\documentclass[11pt,a4paper]{article}
\RequirePackage{amsmath}%
\RequirePackage{amsthm}%
\usepackage{amsmath, amssymb}
\usepackage{amsthm}
\usepackage{url}

\usepackage{color}

\usepackage{amsfonts}%
\usepackage{color}
\usepackage[english]{babel}
\usepackage{verbatim,here}
\usepackage[T1]{fontenc}
\usepackage{graphicx}
\usepackage{epstopdf}
\usepackage{floatflt}
\usepackage{color,xcolor}
\usepackage{enumerate}
\usepackage{titlesec}
\usepackage{a4wide}

\usepackage{graphicx}%
\usepackage{color}%
\usepackage{url}%

\usepackage[autolinebreaks,framed]{mcode}
\lstnewenvironment{matlab}[1][]
{\lstset{numbers=left,xleftmargin=2.5em,framexleftmargin=2em,language=Matlab,float=tp,%
		floatplacement=tp,abovecaptionskip=0pt,#1}}
{}
\titleformat{\subsubsection}[runin]
{\normalfont\bfseries}
{\thesubsubsection}
{1em}
{}
[.]

\renewcommand{\Re}{\mathop{\mathrm{Re}}}
\renewcommand{\Im}{\mathop{\mathrm{Im}}}

\renewcommand{\i}{\mathrm{i}}

\newcommand{\CC}{{\mathbb C}}

\newcommand{\bM}{{\bf M}}
\newcommand{\bN}{{\bf N}}
\newcommand{\bI}{{\bf I}}

\newcounter{minutes}
\divide\time by 60
\newcounter{hours}
\multiply\time by 60 \addtocounter{minutes}{-\time}
\begin{document}

	\title{Numerical computation of Mityuk's function and radius for circular/radial slit domains}
	
	\author{El Mostafa Kalmoun$^{\rm a}$, Mohamed M S Nasser$^{\rm b}$, and Matti Vuorinen$^{\rm c}$}
	
	\date{}
	\maketitle
	
	\vskip-0.8cm %
	\centerline{$^{\rm a,b}$Department of Mathematics, Statistics and Physics, Qatar University,} %
	\centerline{P.O. Box 2713, Doha, Qatar.}%
	\centerline{$^{\rm c}$Department of Mathematics and Statistics, University of Turku,} %
	\centerline{Turku, Finland.}%
	\centerline{E-mail: $^{\rm a}$ekalmoun@qu.edu.qa, $^{\rm b}$mms.nasser@qu.edu.qa, $^{\rm c}$vuorinen@utu.fi}
		
\begin{center}
\begin{quotation}
{\noindent \bf Abstract.}\;\; We consider Mityuk's function and radius which have been proposed in \cite{Mit} as generalizations of the reduced modulus and conformal radius to the cases of multiply connected domains.
We present a numerical method to compute Mityuk's function and radius for canonical domains that consist of the unit disk with circular/radial slits. Our method is based on the boundary integral equation with the generalized Neumann kernel. Special attention is given to the validation of the theoretical results on the existence of critical points and the boundary behavior of Mityuk's radius.
\end{quotation}
\end{center}
\begin{center}
\begin{quotation}
{\noindent {\bf Keywords.\;\;}%
Mityuk's function/radius, reduced modulus, conformal radius, numerical conformal mapping.
}%
\end{quotation}
\end{center}
	
	\begin{center}
		\begin{quotation}
			{\noindent {\bf MSC.\;\;}
				65E05; 30C30; 65R20.
			}
		\end{quotation}
	\end{center}

\maketitle

\footnotetext{\texttt{{\tiny File:~\jobname .tex, printed: \number\year-%
\number\month-\number\day, \thehours.\ifnum\theminutes<10{0}\fi\theminutes}}}
\makeatletter

\makeatother
\section{Introduction}
\label{sc:int}
Let $G\subset {\CC}$ be a bounded simply connected domain and $\alpha\in G$.
Then for $\varepsilon \in (0,d(\alpha,\partial G)/2)$
 the domain $G_\varepsilon=G\backslash\{z\,:\, |z-\alpha|<\varepsilon\}$
 is doubly connected. The quantity
\[
m(G,\alpha)=\lim_{\varepsilon\to0}\left(M(G_\varepsilon)+\frac{1}{2\pi}\log\varepsilon\right)
\]
is called the \emph{reduced modulus} of the domain $G$ with respect to the point $\alpha$ where $M(G_\varepsilon)$ is the modulus of the doubly connected domain $G_\varepsilon$~\cite{Du, Vas02}. If $w=\Phi_\alpha(z)$ is the unique conformal map of $G$ onto the unit disk $|w|<1$ such that
\begin{equation}\label{eq:cm-cond}
\Phi_\alpha(\alpha)=0 \;\text{ and }\; \Phi'_\alpha(\alpha)>0,
\end{equation}
then the conformal radius of $G$ with respect to the point $\alpha$ is defined by~\cite{Mit}
\begin{equation}\label{eq:R}
R(G,\alpha)=\frac{1}{\Phi_\alpha'(\alpha)}.
\end{equation}
The reduced modulus of the domain $G$ with respect to the point $\alpha$ is then given by~\cite[p.~16]{Vas02}
\begin{equation}\label{eq:m}
m(G,\alpha)=\frac{1}{2\pi}\log R(G,\alpha)=-\frac{1}{2\pi}\log\Phi_\alpha'(\alpha).
\end{equation}

Vasil'ev~\cite[Section 2.2.1]{Vas02} assumed that the mapping function $w=\hat\Phi_\alpha(z)$ satisfies the normalization $\hat\Phi_\alpha(\alpha)=0$ and $\hat\Phi'_\alpha(\alpha)=1$. Hence $w=\hat\Phi_\alpha(z)$ maps the domain $G$ onto the disk $|w|<R$ where $R=R(G,\alpha)$ is the conformal radius of $G$ with respect to the point $\alpha$. This is equivalent to the above definition~\eqref{eq:R} since $\Phi_\alpha(z)=\hat\Phi_\alpha(z)/R(G,\alpha)$.

A generalization of the reduced modulus to multiply connected domains has been proposed by Mityuk~\cite{Mit}. Assume that $G\subset\overline{\CC}$ is a bounded multiply connected domain of connectivity $\ell+1$ bordered by
\[
\Gamma=\partial G=\bigcup_{k=0}^{\ell}\Gamma_k
\]
with Jordan curves $\Gamma_0,\Gamma_1,\ldots,\Gamma_m$ such that $\Gamma_0$ encloses all the other curves. Let $\Omega \subset \overline{\CC}$ be the canonical multiply connected domain consisting of the unit disk $|w|<1$ with $p$ circular slits and $\ell-p$ radial slits where $0\le p\le\ell$. Then there exists a unique conformal mapping $w=\Phi_\alpha(z)$ from $G$ onto $\Omega$ holding the same normalization~\eqref{eq:cm-cond}~\cite[Theorem 6, p. 242]{Gol69}.
The mapping function $\Phi_\alpha$ depends also on the piecewise constant real-valued function $\theta$ defined on $\Gamma$ by
\begin{equation*}
\theta(\zeta)=\left\{ \begin{array}{l@{\hspace{0.5cm}}l}
	\pi/2 ,&\zeta\in \Gamma_0,\\
	\theta_1 ,&\zeta\in \Gamma_1,\\
	\hspace{0.3cm}\vdots\\
	\theta_\ell,&\zeta\in \Gamma_\ell,
	\end{array}
	\right.
\end{equation*}
where, for each $k=1,\ldots,\ell$, $\theta_k$ is the ``oblique angle'' of the slit which represents the angle of intersection between the slit and any ray emanating from the origin (see~\cite[p.~109]{Wen92}). This means, for
$k=1,\dots,\ell$
\[
\theta_k=\left\{ \begin{array}{l@{\hspace{0.25cm}}l}
	\pi/2 &\mbox{if $\Gamma_k$ is mapped to a circular slit},\\
	0     &\mbox{if $\Gamma_k$ is mapped to a radial slit}.
	\end{array}
	\right.
\]
So, it is maybe more appropriate to denote the mapping function by $w=\Phi_{\alpha,\theta}(z)$. However, for convenience, we will omit to mention the subscript $\theta$ as it will be clear from the context.

Under this setting, the definitions~\eqref{eq:R} and~\eqref{eq:m} of the functions $R(G,\alpha)$ and $m(G,\alpha)$, respectively, can be extended to the multiply connected domain $G$. In this context, the function $R(G,\alpha)$ is called Mityuk's radius of the domain $G$ with respect to the point $\alpha\in G$ and the canonical domain $\Omega$~\cite{Kin18}. Similarly, the function $m(G,\alpha)$, which is called Mityuk's function in~\cite{Kin18}, is the generalized reduced module of the multiply connected domain $G$ with respect to the point $\alpha\in G$ and the canonical domain $\Omega$.
It is worth mentioning that for the canonical domain of the unit disk with circular slits, Mityuk's function $m(G,\alpha)$ is the negative of the Robin function related to the motion of a single point vortex in the domain $G$~\cite{CM05,NSML}.

One important motivation for the study of Mityuk's function/radius is the connection of their critical points to the solutions of exterior inverse boundary value problems~\cite{Eli17,gak77}.
By rewriting the function $\Phi_\alpha$ as
\begin{equation}\label{eq:Phi-phi}
\Phi_\alpha(z)=(z-\alpha)\phi_\alpha(z),
\end{equation}
where $\phi_\alpha$ is an analytic function in $G$, the critical points of Mityuk's function/radius are the roots of the equation
$\phi_\alpha'(\alpha) = 0$ \cite{Kin18}.
The existence of these critical points for $\ell\ge2$ was proven in the case of circular concentric slits in~\cite{Kin84} and in the case of a mix of circular and radial slits domain in~\cite{Kin18}. In the former case, Kinder~\cite{Kin84} indicated that the nature of critical points can be specified by the equation
\begin{equation}\label{eq:criticalPts}
n_m-n_s=1-\ell ,
\end{equation}
where $n_m\geq 1$ is the number of local maxima and $n_s$ is the number of saddle points.
When $G$ is a doubly connected domain ($\ell=1$), an example was constructed in \cite{Kin18} where Mityuk's radius has no critical points.

Furthermore, given that the domain $G$ has smooth boundary curves, it was proven in~\cite{Kin18} that Mityuk's radius $R(G,\alpha)$ is infinitely differentiable on $G$ and has the following limit values on the boundary $\Gamma$.
First, for the external boundary $\Gamma_0$, we have
\begin{equation}\label{eq:lim-k}
\lim_{\alpha\to\beta\in\Gamma_0}R(G,\alpha)=0.
\end{equation}
Then, for the internal boundary components $\Gamma_k$, $k=1,2,\ldots,\ell$,
\begin{equation}\label{eq:lim-j}
\lim_{\alpha\to\beta\in\Gamma_k}R(G,\alpha)=\left\{ \begin{array}{l@{\hspace{0.25cm}}l}
	0     &\mbox{if $\Gamma_k$ is mapped to a circular slit},\\
	\infty&\mbox{if $\Gamma_k$ is mapped to a radial slit}.
	\end{array}
	\right.
\end{equation}

\begin{figure}[t] %
\centerline{
\scalebox{0.4}{\includegraphics[trim=0.0cm 0.0cm 0.0cm 0.0cm,clip]{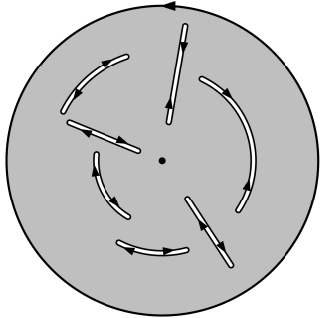}}}
\caption{The canonical domain $\Omega$: The unit disk with circular/radial slits.}
\label{fig:sp}
	\end{figure}

In this paper, we shall present a numerical method to compute Mityuk's radius $R(G,\alpha)$ and Mityuk's function $m(G,\alpha)$ of the domain $G$ with respect to the point $\alpha\in G$ and the canonical domain $\Omega$ consisting of the unit disk with $\ell$ circular/radial slits (see Figure~\ref{fig:sp}). The presented method depends on the boundary integral equation with the generalized Neumann kernel. Through various numerical examples, we aim to validate the theoretical results presented in~\cite{Kin84,Kin18} on the existence of critical points and on the boundary behavior of Mityuk's radius. In particular, we will numerically examine the validity of the limits~\eqref{eq:lim-k} and~\eqref{eq:lim-j} in the case of multiply connected domains with smooth, piecewise smooth and slit boundaries.
At the end of the paper we summarize our experimental discoveries which suggest some problems for theoretical investigation.

	\section{Preliminaries}
	
Let us first describe briefly the boundary parametrization used in this paper; see~\cite{Nas-ETNA} for more details. We parametrize each boundary component $\Gamma_k$, for $k=0,1,\ldots,\ell$, by a $2\pi$-periodic complex function $\eta_k(t)$, $t\in J_k:=[0,2\pi]$. We assume that $\eta_k(t)$ is  twice continuously differentiable with $\eta'_k(t)\ne 0$ for $t\in J_k$.
The whole parameter domain is defined by
	\begin{equation*}
	J = \bigsqcup_{k=0}^{\ell} J_k=\bigcup_{k=0}^{\ell}\{(t,k)\;:\;t\in J_k\}.
	\end{equation*}
	A parametrization of the whole boundary $\Gamma$ is then defined on $J$ by
	\begin{equation}\label{e:eta-1}
	\eta(t,k)=\eta_k(t), \quad t\in J_k,\quad k=0,1,\ldots,\ell.
	\end{equation}
As the value of the auxiliary index $k$ for a given $t$ will be clear from the context, we will drop it in the notation of the function $\eta$ in~(\ref{e:eta-1}) and simply write
	\begin{equation*}
	\eta(t)= \left\{ \begin{array}{l@{\hspace{0.5cm}}l}
	\eta_0(t),&t\in J_0,\\
	\hspace{0.3cm}\vdots\\
	\eta_\ell(t),&t\in J_\ell.
	\end{array}
	\right.
	\end{equation*}

On the other hand, the smoothness of $\eta(t)$ allows us to identify any function $\phi $ in the  space of all real-valued H\"older continuous functions $H$ on $\Gamma$ with a $2\pi$-periodic H\"older continuous real function $\hat\phi$ of the parameter $t$ in $J_k$ by writing $\hat\phi(t):=\phi(\eta(t))$. 
Therefore, we shall not distinguish between $\phi(\eta(t))$ and $\phi(t)$ in the sequel. 

In order to introduce the boundary integral equation that will be used in the next section to compute Mityuk's function and radius, we recall the definition of the generalized Neumann kernel formed with a continuously differentiable complex function $A$ on $J$ and the parametrization function $\eta$~\cite{Mur-Bul,Weg-Mur-Nas,Weg-Nas}
\begin{equation*}
N(s,t) =  \frac{1}{\pi}\Im\left(
\frac{A(s)}{A(t)}\frac{\eta'(t)}{\eta(t)-\eta(s)}\right).
\end{equation*}
We recall also the following associated kernel $M(s,t)$ on $J\times J$~\cite{Weg-Nas}
\begin{equation*}
M(s,t) =  \frac{1}{\pi}\Re\left(
\frac{A(s)}{A(t)}\frac{\eta'(t)}{\eta(t)-\eta(s)}\right).
\end{equation*} 
The first kernel is continuous and the second has a singularity of cotangent type~\cite{Weg-Nas}. 
Using these two kernels, the integral operators $\bN, \bM: H\to H$ defined by
	\begin{eqnarray*}
	\bN\mu(s) &=& \int_J N(s,t) \mu(t) dt, \quad s\in J, \\
	\bM\mu(s) &=& \int_J  M(s,t) \mu(t) dt, \quad s\in J,
	\end{eqnarray*}
are both bounded with the first one being compact and the second being singular. Further details can be found in~\cite{Weg-Mur-Nas,Weg-Nas}.
	
\section{Computing Mityuk's radius and Mityuk's function}
We use a boundary integral equation involving the integral operators $\bN$ and $\bM$ to compute Mityuk's function/radius with respect to the canonical domain $\Omega$ in Figure~\ref{fig:sp}.
To this aim, we first review the following method from~\cite[Section 4.2]{Nas-jmaa11} for computing the unique conformal mapping $w=\Phi_\alpha(z)$ satisfying the normalization~\eqref{eq:cm-cond} from the multiply connected domain $G$ in $z$-plane onto the canonical circular/radial slit domain $\Omega$ in $w$-plane. We assume that the conformal mapping $w=\Phi_\alpha(z)$ maps the curve $\Gamma_0$ onto the unit circle $|w|=1$ and maps the curve $\Gamma_k$, for each $k=1,2,\ldots,m$, onto a circular or radial slit, i.e., for $\eta_k(t)\in\Gamma_k$,
\begin{equation*}
\Im\left[e^{-\i\theta_k}\log\Phi_\alpha(\eta_k(t))\right]=R_k, \quad k=1,2,\ldots,\ell,
\end{equation*}
where $\theta_k=0$ for the radial slit case and $\theta_k=\pi/2$ for the circular slit case. The undetermined real constants $R_1,\ldots,R_\ell$ should be computed alongside the conformal mapping $\Phi_\alpha$.  Thus, the boundary values of the conformal mapping $w=\Phi_\alpha(z)$ satisfy
\begin{equation}\label{eq:d-cd1}
\Im\left[e^{-\i\theta(t)}\log\Phi_\alpha(\eta(t))\right]
=R(t),\quad \eta(t)\in\Gamma,
\end{equation}
where $R(t)=(0,R_1,\ldots,R_\ell)$.	Here and in what follows, the notation $h(t)=(h_0,\ldots,h_\ell)$ stands for the real piecewise constant function
\[
h(t) = \left\{
\begin{array}{l@{\hspace{0.5cm}}l}
h_0,     & t\in J_0, \\
\vdots  & \\
h_\ell,  & t\in J_\ell, \\
\end{array}%
\right.
\]
where $h_0,\ldots,h_\ell$ are real constants.

The mapping function $\Phi_\alpha$ can be written as
\begin{equation*}
w=\Phi_\alpha(z)=c(z-\alpha)e^{(z-\alpha)f(z)}, \quad z\in G\cup\Gamma,
\end{equation*}	
where
\begin{equation*}
c=\Phi_\alpha'(\alpha)>0
\end{equation*}
and $f(z)$ is an analytic function on $G$.
Therefore, for $\eta(t)\in\Gamma$,
\[
\log\Phi_\alpha(\eta(t))=\log c+\log(\eta(t)-\alpha)+(\eta(t)-\alpha)f(\eta(t)).
\]
Multiplying both sides by $e^{-\theta(t)}$, taking the imaginary of both sides, and using~\eqref{eq:d-cd1}, we conclude that the boundary values of the function $f$ are given by
\begin{equation}\label{eq:Af}
A(t)f(\eta(t))=\gamma(t)+h(t)+\i\mu(t),\quad\eta(t)\in\Gamma,
\end{equation}
where
\begin{equation*}
A(t) = e^{\i\left(\frac{\pi}{2}-\theta(t)\right)}\,(\eta(t)-\alpha), \quad \gamma(t)=\Im\left[e^{-\i\theta(t)}\log\left(\eta(t)-\alpha\right)\right]
\end{equation*}
and
\begin{equation}\label{eq:h2}
h(t)=(-\log c,-R_1-\log c\sin\theta_1,\ldots,-R_\ell-\log c\sin\theta_\ell).
\end{equation}
In~\eqref{eq:Af}, only $\gamma$ is known whereas $h$ and $\mu$ are unknown. The function $\mu$ is the unique solution of the integral equation
\begin{equation}\label{eq:ie}
(\bI-\bN)\mu=-\bM\gamma,
\end{equation}
and the piecewise constant real-valued function
\begin{equation}\label{eq:h-k}
h=(h_0,h_1,\ldots,h_\ell)
\end{equation}
is given by
\begin{equation}\label{eq:h}
h=[\bM\mu-(\bI-\bN)\gamma]/2.
\end{equation}
By computing $\mu$ and $h$, we obtain the boundary values of the auxiliary analytic function $f$ from~\eqref{eq:Af}. Moreover, comparing~\eqref{eq:h-k} with~\eqref{eq:h2} implies that
\begin{equation}\label{eq:c}
\Phi_\alpha'(\alpha)=c=e^{-h_0}.
\end{equation}

Now, we are in position to compute Mityuk's radius and Mityuk's function. In view of~\eqref{eq:c}, it follows from~\eqref{eq:R} and~\eqref{eq:m} that Mityuk's radius and Mityuk's function of the domain $G$ with respect to the point $\alpha\in G$ and the canonical domain $\Omega$ can be computed through
\begin{equation}\label{eq:R-m}
R(G,\alpha)=e^{h_0}\quad{\rm and}\quad m(G,\alpha)=\frac{h_0}{2\pi}.
\end{equation}

A MATLAB function \verb|fbie| for solving numerically the integral equation~\eqref{eq:ie} has been presented in~\cite{Nas-ETNA}. 
The function \verb|fbie| is based on discretizing the integral equation~(\ref{eq:ie}) by the Nystr\"om method with the trapezoidal rule~\cite{Atk97,Tre-Trap1,Tre-Trap} to obtain a linear system. For domains with corners, the trapezoidal rule with a graded mesh is used, see~\cite{Kre,LSN} for more details.
The MATLAB function $\mathtt{gmres}$ is then used to solve iteratively this linear system. The matrix-vector multiplication required by the GMRES method is computed efficiently using the MATLAB function $\mathtt{zfmm2dpart}$ from the toolbox $\mathtt{FMMLIB2D}$~\cite{Gre-Gim12}. 
The function \verb|fbie| is also used to compute approximations of the function $h$ in~\eqref{eq:h}.
The total computational cost of the overall method is $O((\ell+1)n\ln n)$ operations where $n$ is the number of nodes in each of the intervals $J_0,J_1,\ldots,J_\ell$. For more details, see~\cite{Nas-ETNA}.

If the discretized versions of the functions $\eta(t)$, $\eta'(t)$, $A(t)$, $\gamma(t)$, $\mu(t)$ and $h(t)$ are denoted by the vectors \verb|et|, \verb|etp|, \verb|A|, \verb|gam|, \verb|mu| and \verb|h|, respectively, then the linear systems \eqref{eq:ie} and \eqref{eq:h} are numerically solved by calling
\[
[\verb|mu|,\verb|h|] = \verb|fbie|(\verb|et|,\verb|etp|,\verb|A|,\verb|gam|,\verb|n|,\verb|iprec|,\verb|restart|,\verb|gmrestol|,\verb|maxit|).
\]
In the numerical experiments presented in this paper, we choose $\mathtt{iprec}=5$ which means the tolerance of the FMM is $0.5\times 10^{-15}$. For the GMRES method, we choose $\mathtt{gmrestol}=10^{-14}$ so that the tolerance of the GMRES method is $10^{-14}$. The GMRES is used without restart by choosing $\mathtt{restart}=[\,]$ and the maximum number of GMRES iterations is $\mathtt{maxit}=100$.
Choosing the value of $n$ depends on the geometry of the domain $G$ and the location of the point $\alpha$. For domains with smooth boundaries, accurate numerical results can be obtained for moderate values of $n$ if $\alpha$ is sufficiently far from the boundary of $G$. If $\alpha$ is close to the boundary of $G$ or if the boundary of $G$ has corners, a sufficiently large value of $n$ is required to obtain accurate results. In the numerical examples below, $\alpha$ could be very close to the boundary of $G$ so we shall choose $n=2^{15}$ for all examples.
By obtaining the vector \verb|h|, an approximate value for the constant $h_0$ is computed as the average of the values of \verb|h| over the interval $J_0$. Finally, from $h_0$ we can compute Mityuk's radius and Mityuk's function through~\eqref{eq:R-m}. See Listing~\ref{MityukMatlab} for a MATLAB implementation. All the computer codes of our computations are available through the internet link at
{\centering \url{https://github.com/mmsnasser/mityuk-radius}}.

\begin{matlab}[caption={MATLAB function for computing Mityuk's function/radius for bounded multiply connected domains}, label=MityukMatlab]
function [R,m] = Mityuk(et,etp,n,thetak,alpha)
% Compute Mityuk's radius R(G,alpha) and Mityuk's function m(G,alpha)
% where alpha is a point in G, et, etp: parametrization of the boundary
% and its first derivative, thetk=[theta_1,...,theta_L] are the oblique
% angles of the spirals
for k=1:length(thetak)
    J       = 1+(k-1)*n:k*n;
    A(J,1)  = exp(i.*(pi/2-thetak(k))).*(et(J)-alpha);
    gam(J,1)=imag(exp(-i.*thetak(k)).*clog(et(J)-alpha));
end
[~,h]  =  fbie(et,etp,A,gam,n,5,[],1e-14,100);
h0     =  mean(h(1:n)); R = exp(h0); m = h0/(2*pi);
end
\end{matlab}

\section{Numerical Examples}
In this section, we use the MATLAB function \verb|Mityuk| given in Listing~\ref{MityukMatlab} to compute the values of Mityuk's radius of several multiply connected domains. We examine especially the existence of critical points of Mityuk's radius for different cases of the canonical domain, and check the validity of the limit values in equations~\eqref{eq:lim-k} and~\eqref{eq:lim-j}.
Note that a considerable computational effort is needed for each example. In fact, computing the values of Mityuk's radius $R(G,\alpha)$ at a vector of points $\alpha_k$, $k=1,2,\ldots,p$, requires calling the MATLAB function \verb|Mityuk| for each point $\alpha_k$, i.e., $p$ times. In the numerical examples below, to plot the contour maps of Mityuk's radius $R(G,\alpha)$, we first discretize the domain $G$ to get a mesh-grid and then compute values of $R(G,\alpha)$ at the mesh-grid points. Although the function \verb|Mityuk| takes only 0.5-3 seconds to execute, the overall running time for all mesh-grid points could be several hours.

\subsection{Doubly connected domains}

In the following examples, we consider $G$ to be a doubly connected domain and compute $R(G,\alpha)$  for two cases of the canonical domain: the unit disk with a circular slit and the unit disk with a radial slit.

\subsubsection{Annulus}\label{sec:ann} 

{\color{blue}First, as a special case, we consider an annulus $G=\{z\,:\,q<|z|<1\}$ when $0<q<1$. In~\cite{Kin18}, it was proven that Mityuk's radius with respect to the canonical domain consisting of the unit disk with a radial slit has no critical points in $G$. On the contrary, it was shown in~\cite{Kin1} that Mityuk's radius with respect to the canonical domain consisting of the unit disk with a circular slit has an infinite number of critical points in $G$, which constitute the circle $|\alpha|=\sqrt{q}$.

Furthermore, an analytical formula for Mityuk's radius can be derived in terms of the inner radius $q$ of the annulus (see~\cite{Kin1,Kin18} for more details).
In the case of the canonical domain consisting of the unit disk with a circular slit, Mityuk's radius can be written as
\[
R(G,\alpha)=(1-|\alpha|^2)\prod_{j=1}^{\infty}\frac{(1-q^{2j}|\alpha|^2)(1-q^{2j}/|\alpha|^2)}{(1-q^{2j})^2}.
\]
Similarly, the exact form of Mityuk's radius with respect to the canonical domain consisting of the unit disk with a radial slit is given by
\[
R(G,\alpha)=(1-|\alpha|^2)\prod_{j=1}^{\infty}\left[\frac{(1-q^{2j}|\alpha|^2)(1-q^{2j}/|\alpha|^2)}{(1-q^{2j})^2}\right]^{(-1)^j}.
\]
}

\begin{figure}[t] %
		\centerline{
			\includegraphics[height=5cm,trim=0cm 0cm 0cm 0cm,clip]{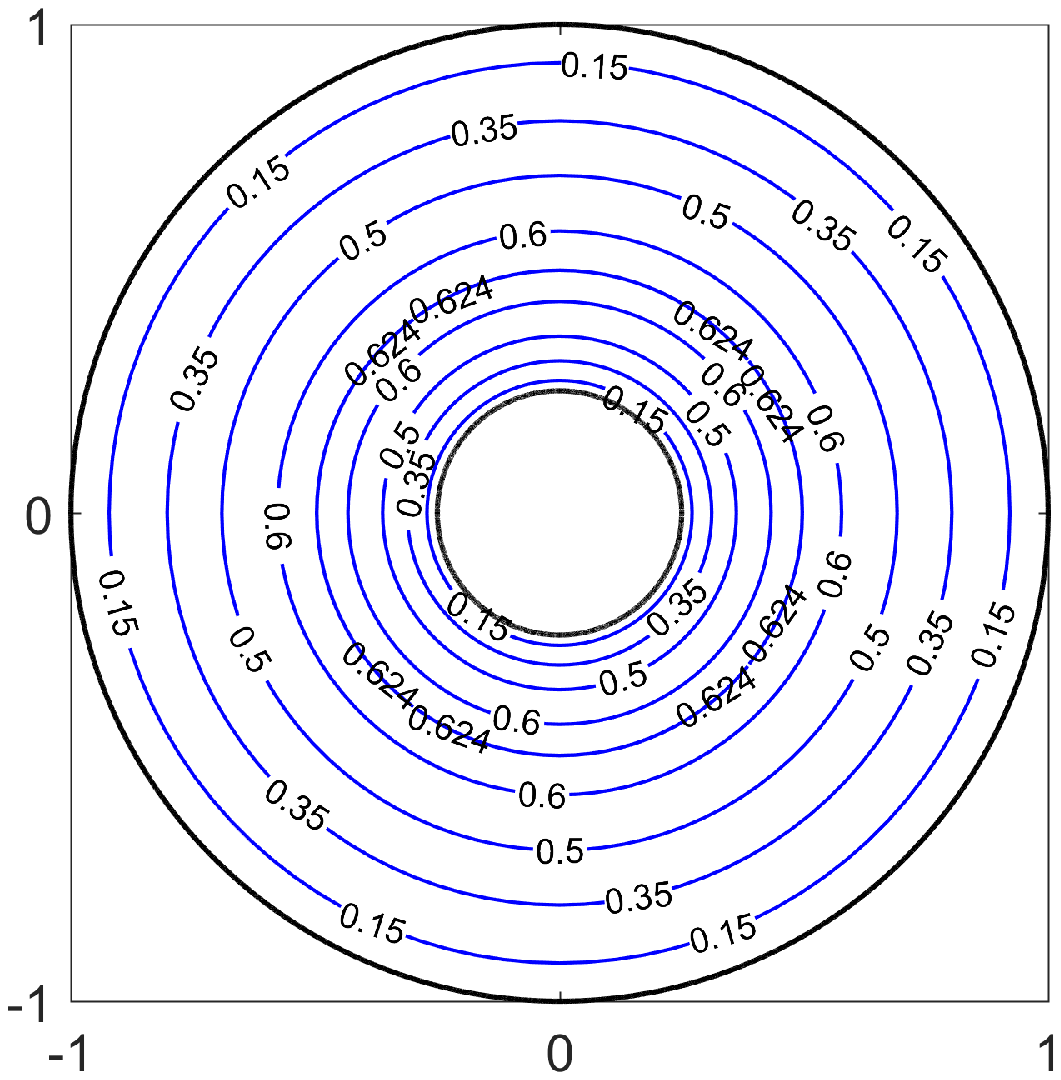}
			\qquad \includegraphics[height=5cm,trim=0cm 0cm 0cm 0cm,clip]{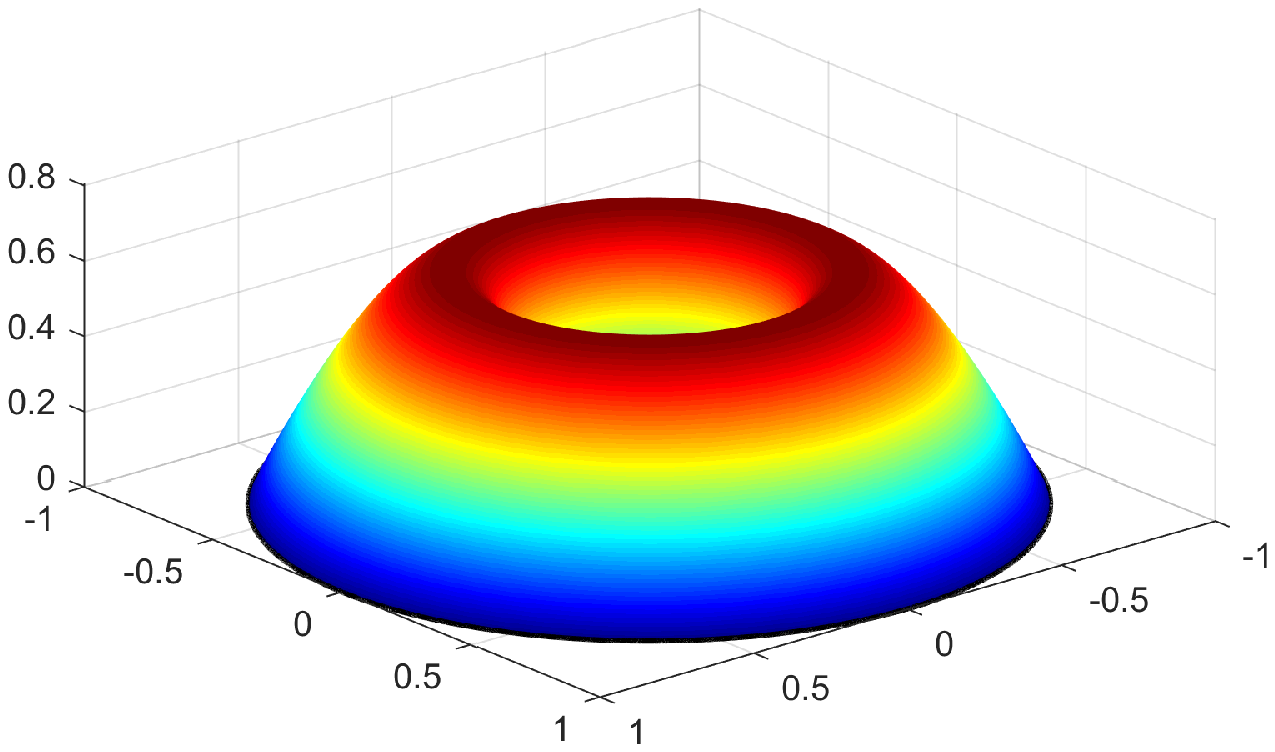}
		}
		\caption{The contour maps and the surface plots of the function $R(G,\alpha)$ for $\theta_1=\pi/2$.}
		\label{fig:an-cir}
	\end{figure}

\begin{figure}[t] %
		\centerline{
			\includegraphics[height=5cm,trim=0cm 0cm 0cm 0cm,clip]{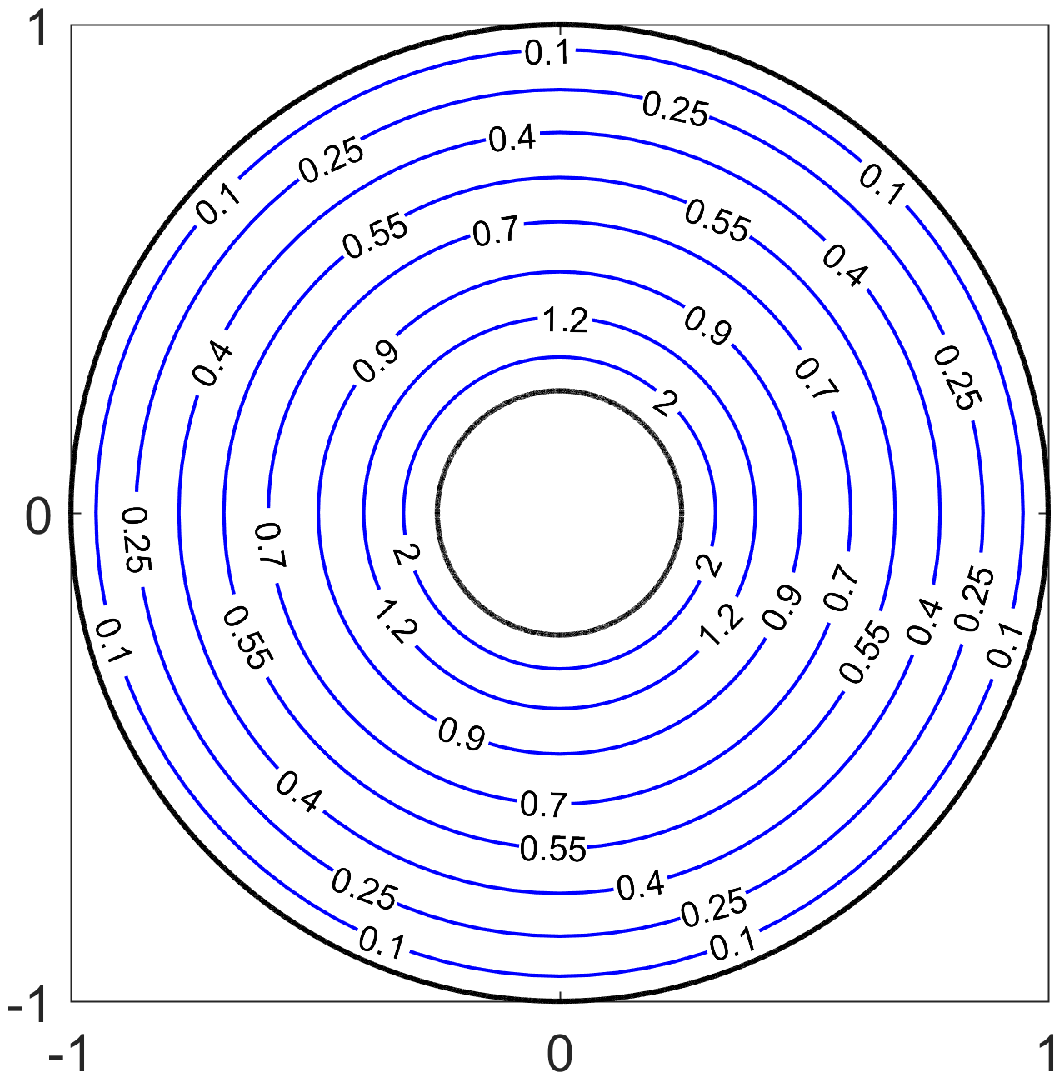}
			\qquad \includegraphics[height=5cm,trim=0cm 0cm 0cm 0cm,clip]{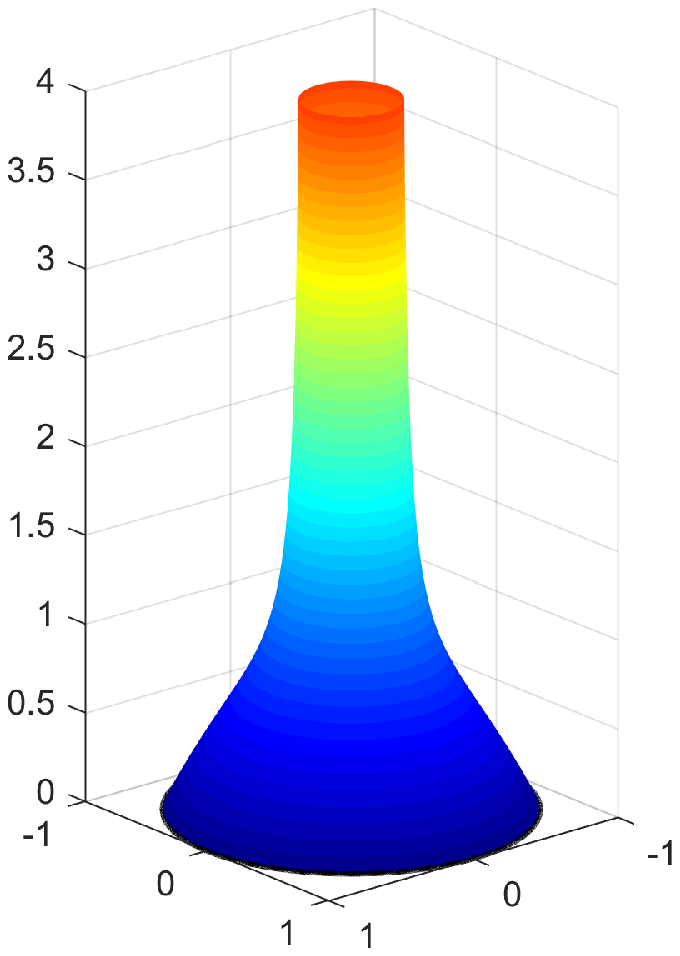}
		}
		\caption{The contour maps and the surface plots of the function $R(G,\alpha)$ for $\theta_1=0$.}
		\label{fig:an-rad}
	\end{figure}

\begin{figure}[t] %
		\centerline{
			\scalebox{0.5}{\includegraphics[trim=0cm 0cm 0cm 0cm,clip]{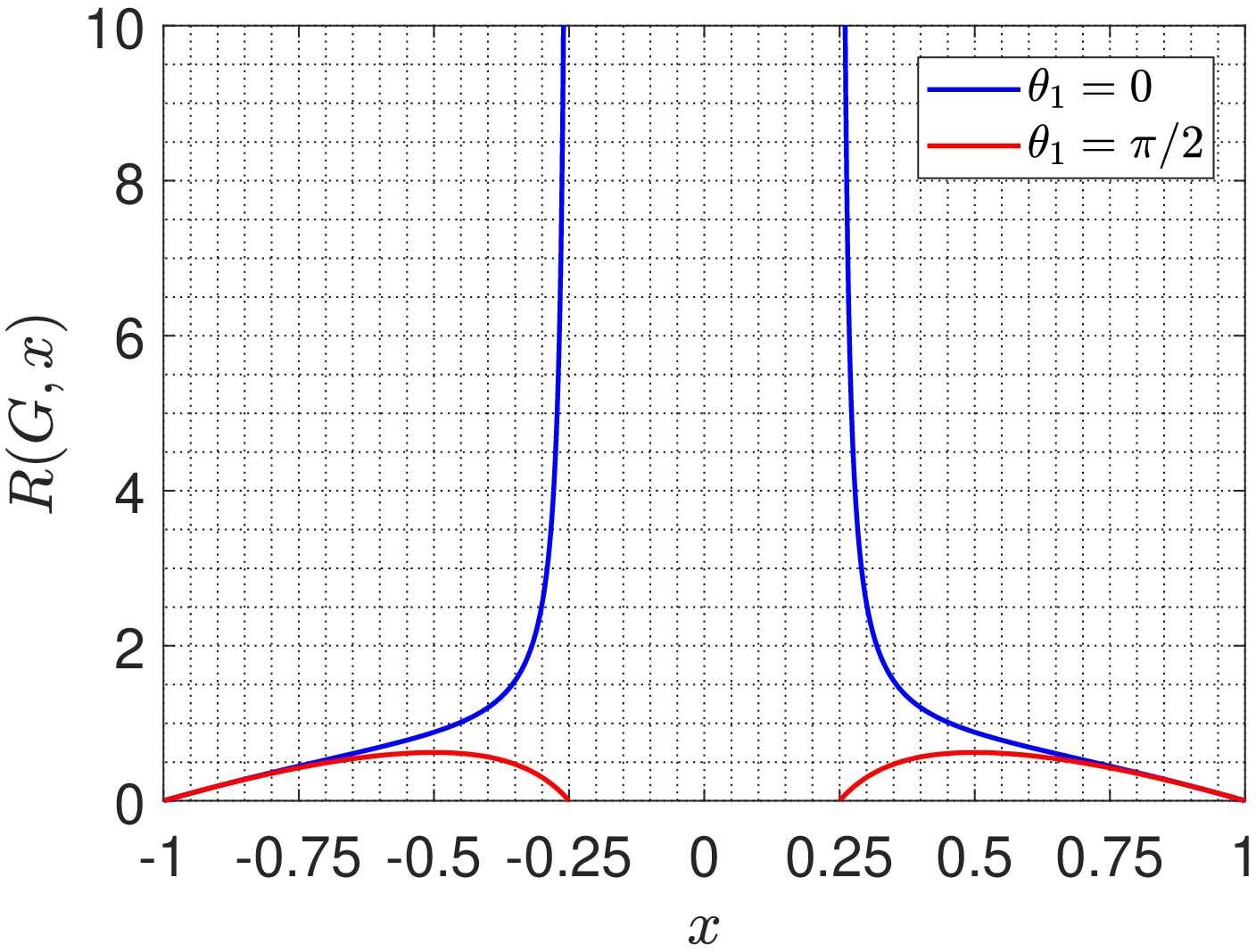}}
		}
		\caption{The values of the radius function $R(G,x)$ for $x\in G$ with $-1< x <1$. The domain $G$ is the same as in Figures~\ref{fig:an-cir} and~\ref{fig:an-rad}.}
		\label{fig:an-x}
	\end{figure}

In order to validate these analytic results, we compute numerically the values of Mityuk's radius for $G=\{z\,:\,q<|z|<1\}$ with $q=0.25$. The contour maps and the surface plots of the function $R(G,\alpha)$ are shown in Figures~\ref{fig:an-cir} and~\ref{fig:an-rad}. Figure~\ref{fig:an-cir} reveals that all the points on the circle $|\alpha|=0.5$ are critical points of the function $R(G,\alpha)$ for the canonical domain consisting of the unit disk with a circular slit, which is an agreement with the {\color{blue} analytical results~\cite{Kin1}.} 
On the other hand, Figure~\ref{fig:an-rad} illustrates the non-existence of critical points of $R(G,\alpha)$ in the case of a radial slit, which confirms the result presented in~\cite{Kin18}.
We also numerically validate the limits~\eqref{eq:lim-k} and~\eqref{eq:lim-j} by computing $R(G,x)$ for several values of $x\in G$ with $-1< x<1$. The graph of $R(G,x)$ as a function of $x$ is shown in Figure~\ref{fig:an-x} for the two cases of the canonical domain. We can see that
	$$\lim_{|x|\to 1^-} R(G,x) = 0$$ for the two cases of the canonical domain as follows from~\eqref{eq:lim-k}, and
	$$
	\lim_{|x|\to 0.25^+}R(G,x)=\left\{ \begin{array}{l@{\hspace{0.25cm}}l}
	0     &\mbox{if $\theta_1 = \frac{\pi}{2}$},\\
	\infty&\mbox{if $\theta_1 = 0$},
	\end{array}
	\right.
	$$
	which is consistent with~\eqref{eq:lim-j}.
	
\subsubsection{Two circles}\label{sec:2cir}
\begin{figure}[t] %
	\centerline{
    	 \scalebox{0.4}{\includegraphics[trim=0cm 0cm 0cm 0cm,clip]{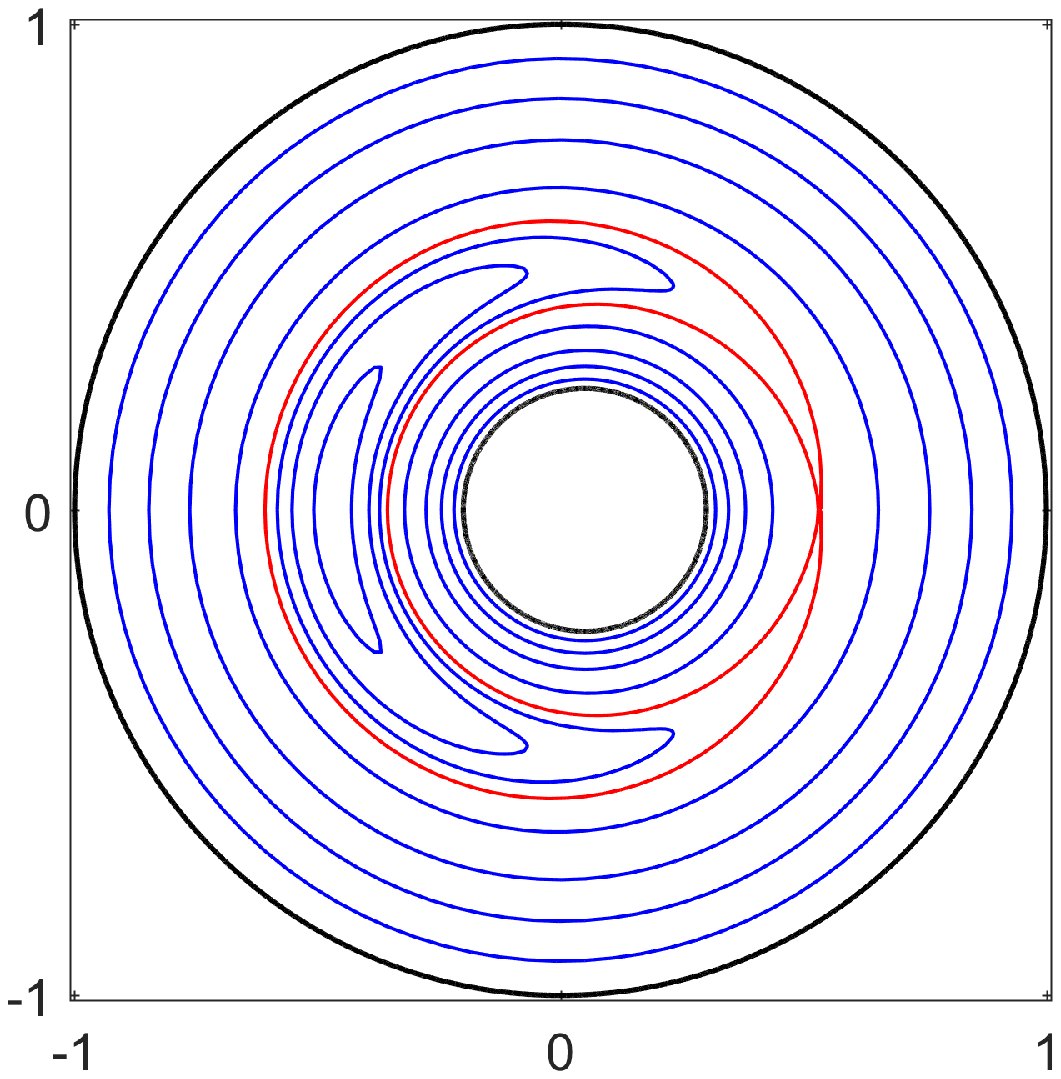}}
\hfill \scalebox{0.4}{\includegraphics[trim=0cm 0cm 0cm 0cm,clip]{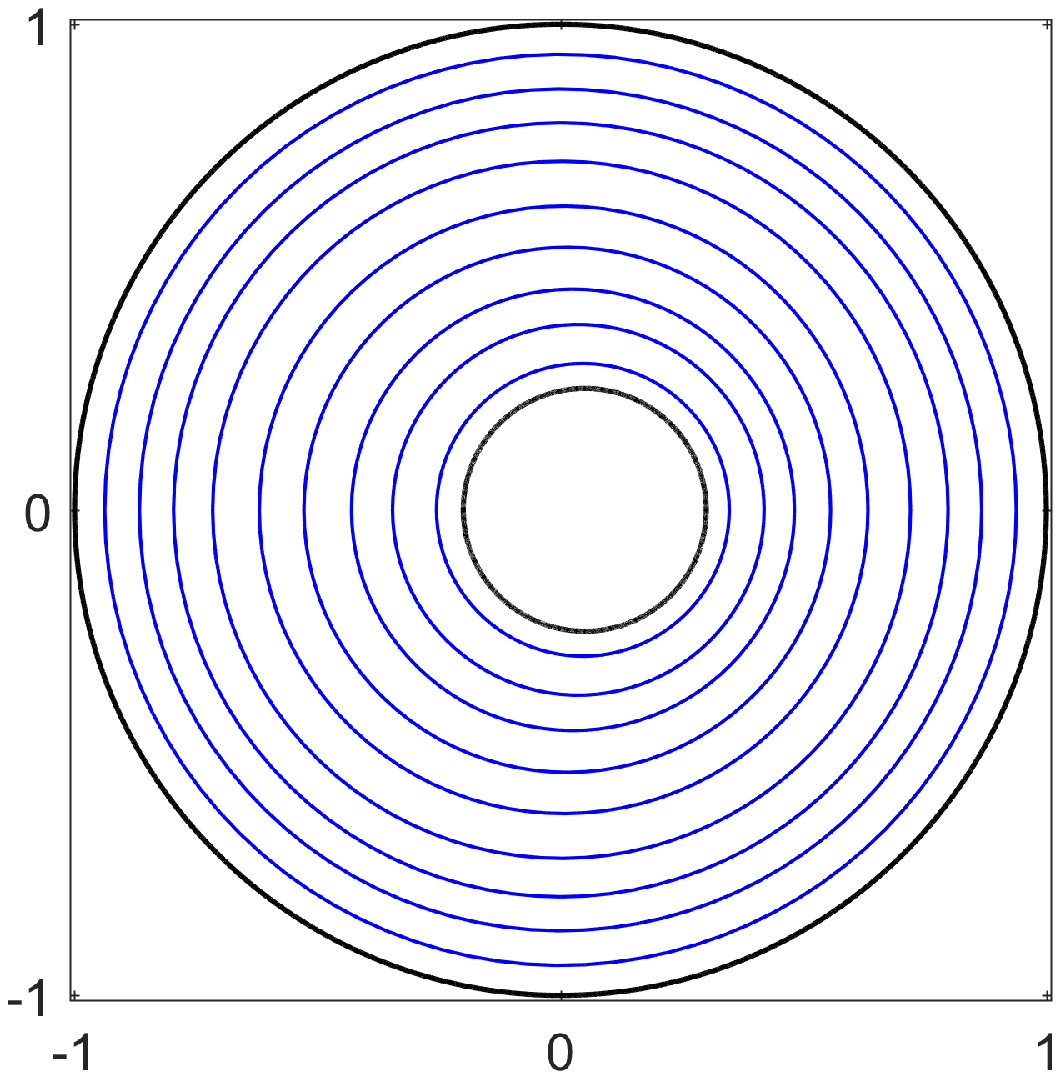}}
\hfill \scalebox{0.4}{\includegraphics[trim=0cm 0cm 0cm 0cm,clip]{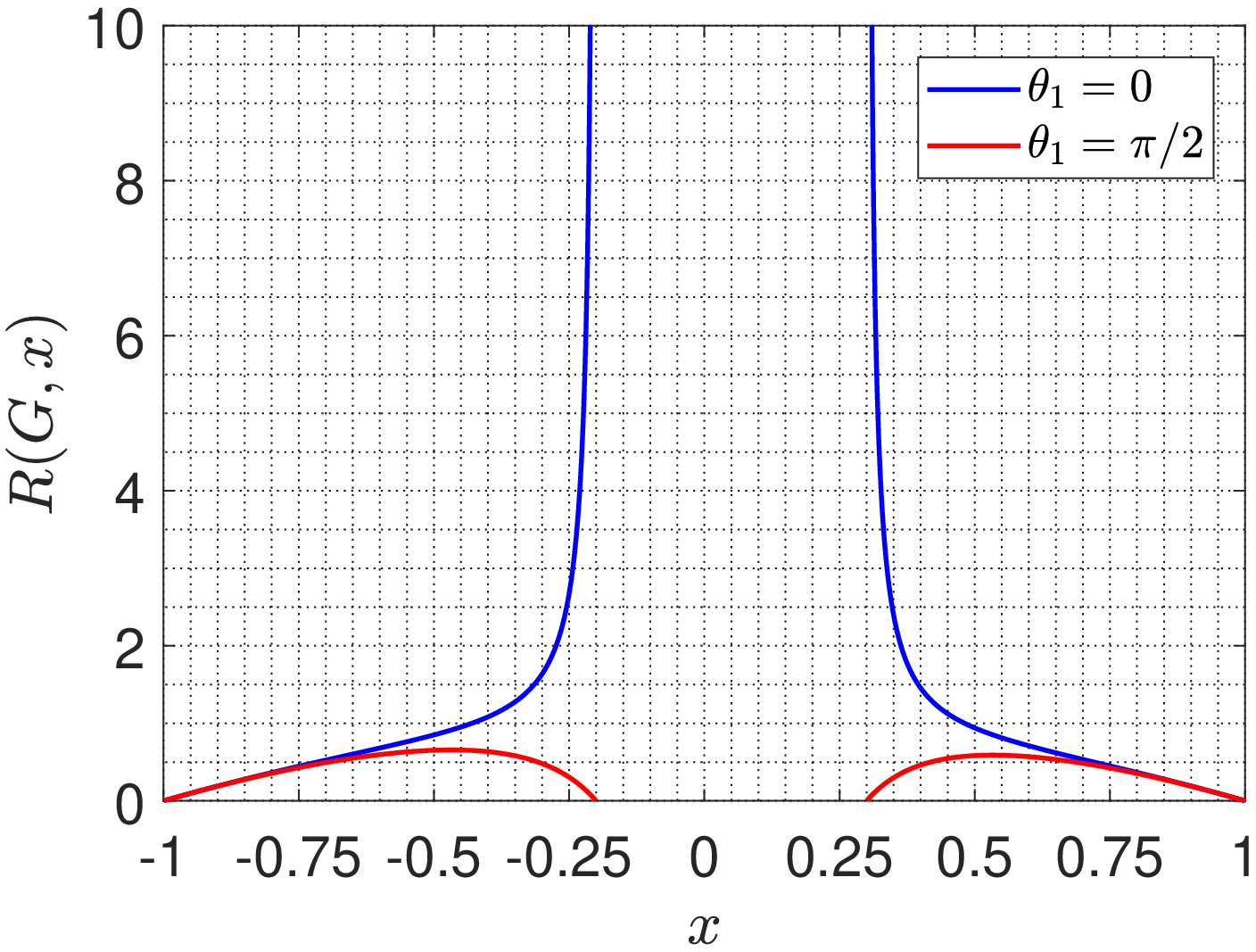}}
	}
	\centerline{
    	 \scalebox{0.4}{\includegraphics[trim=0cm 0cm 0cm 0cm,clip]{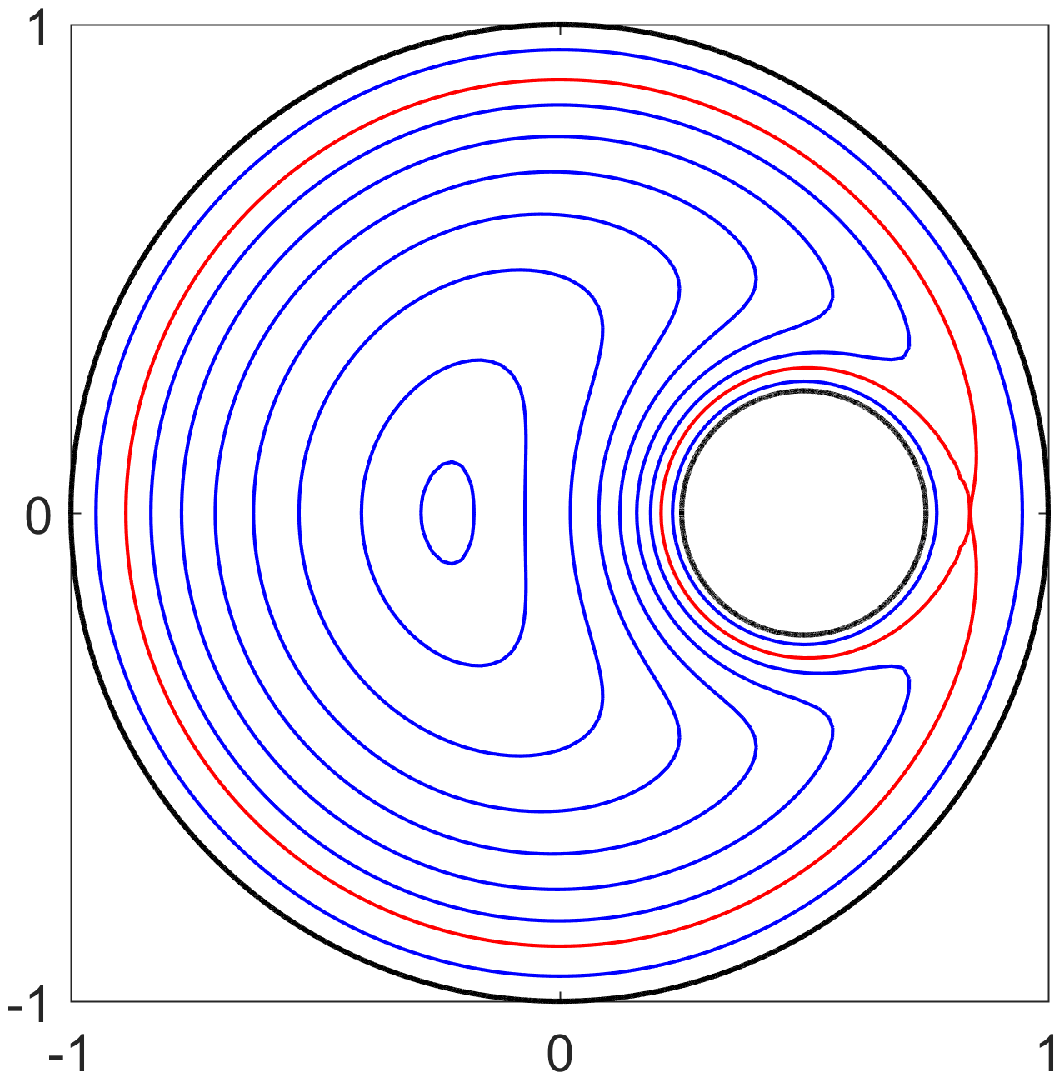}}
\hfill \scalebox{0.4}{\includegraphics[trim=0cm 0cm 0cm 0cm,clip]{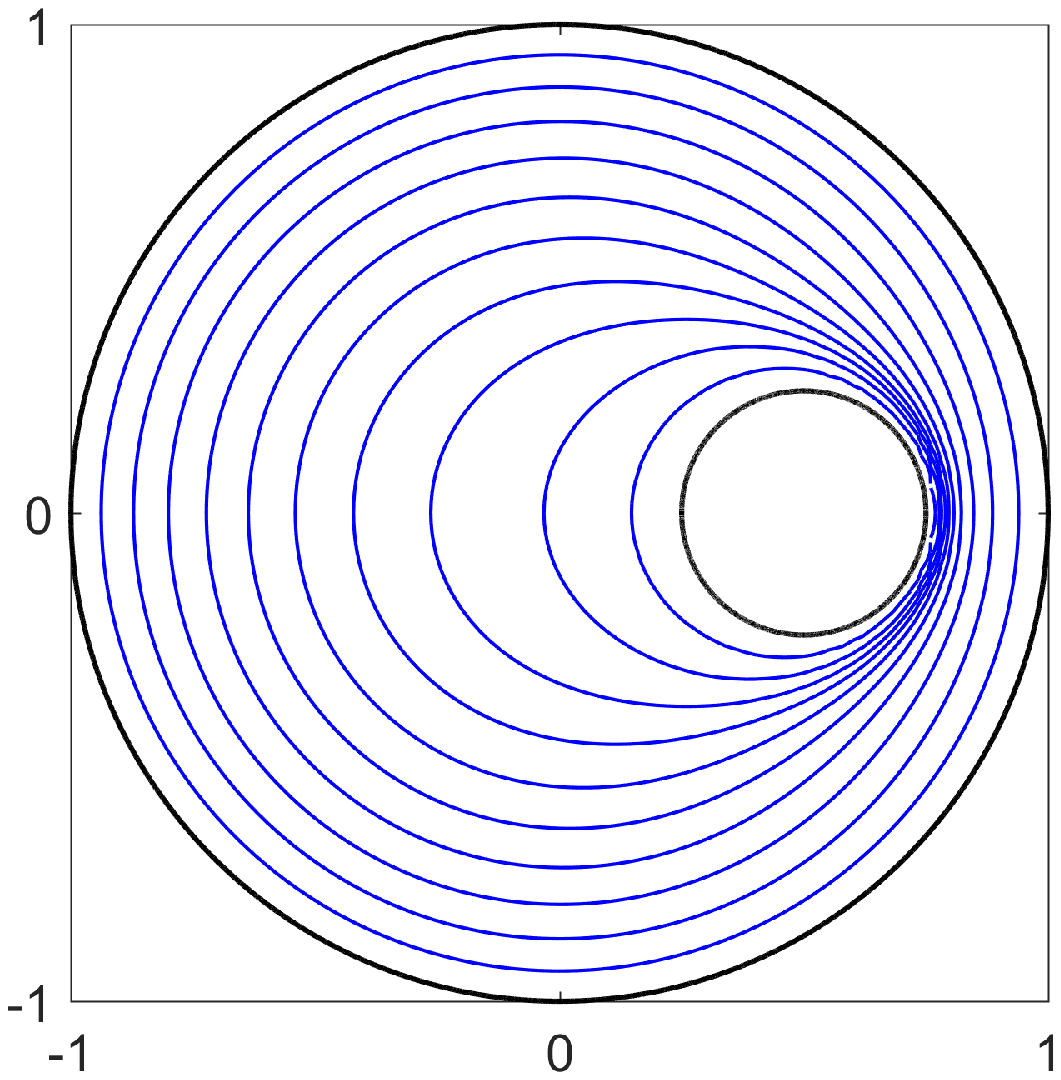}}
\hfill \scalebox{0.4}{\includegraphics[trim=0cm 0cm 0cm 0cm,clip]{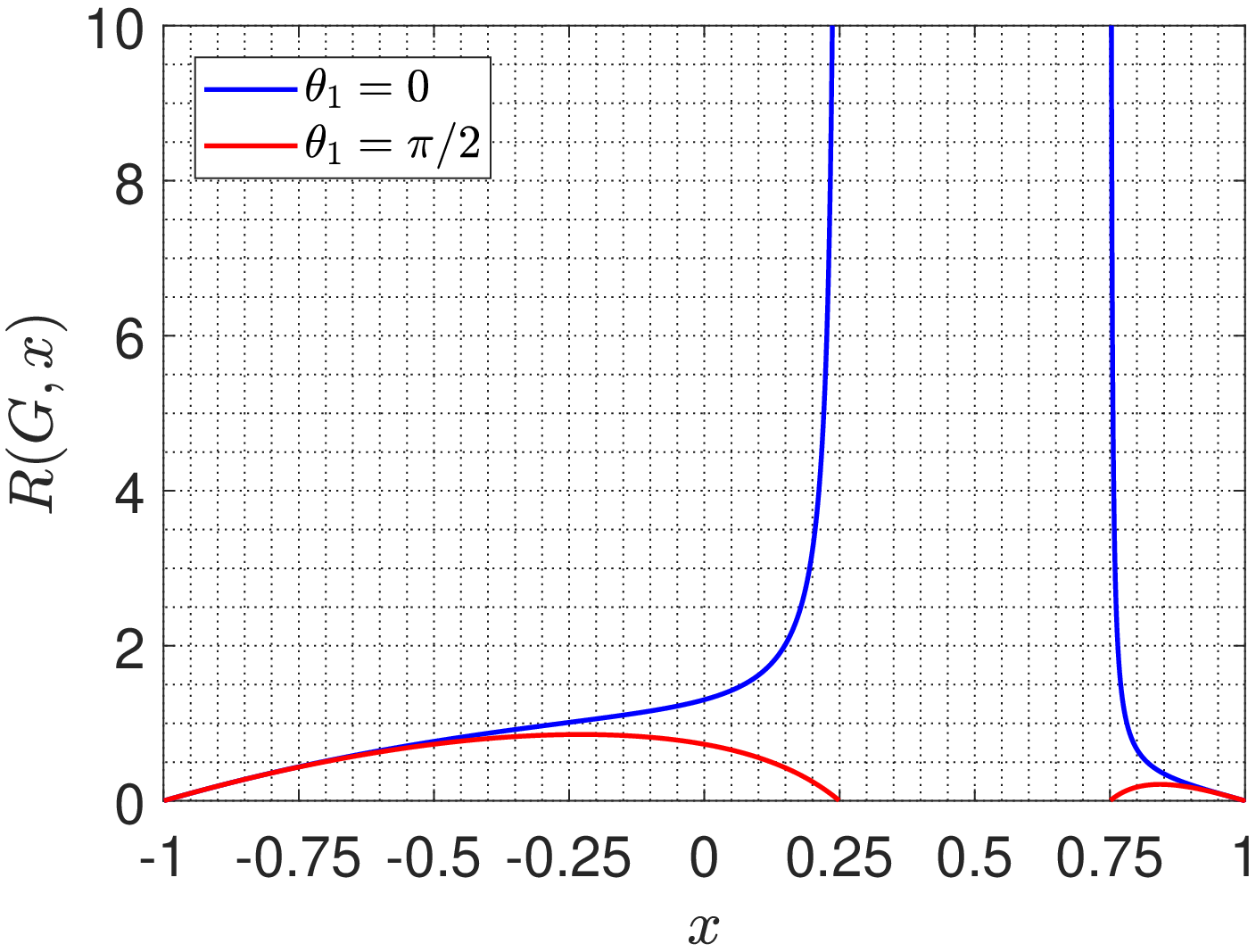}}
	}
	\caption{The contour maps of the function $R(G,\alpha)$ for $\theta_1=\pi/2$ (left), $\theta_1=0$ (center), and the values of $R(G,x)$ for $x\in G$ with $-1< x <1$ (right) where the center of the inner circle is $a=0.05$ for the first row and $a=0.5$ for the second row. Critical streamlines are shown in red color.}
	\label{fig:cir2}
\end{figure}

Let $G$ be the domain interior to the circle $|z|=1$ and exterior to the circle $|z-a|=0.25$ for $a=0.05$ and $a=0.5$.
The contour maps of the function $R(G,\alpha)$ are shown in Figure~\ref{fig:cir2} for $a=0.05$ (first row) and $a=0.5$ (second row). As we can see from this figure, we have two critical points ($n_m=n_s=1$) for the unit disk with a circular slit and no critical point in the case of a radial slit. Interestingly, shifting the origin of the inner circle in the annulus changes dramatically the number of critical points in the case of a circular slit.
On the other hand, the graph of $R(G,x)$ as a function of $x\in G$ such that  $-1<x<1$ is shown in Figure~\ref{fig:cir2} for the two cases of the canonical domain. Notice that limit values of Mityuk's radius at domain boundaries as given in equations~\eqref{eq:lim-k} and~\eqref{eq:lim-j} are again confirmed for this example.

\subsubsection{Rectangle with a slit}\label{sec:rs}
\begin{figure}[t] %
\centerline{
\scalebox{0.4}{      \includegraphics[trim=0cm 0cm 0cm 0cm,clip]{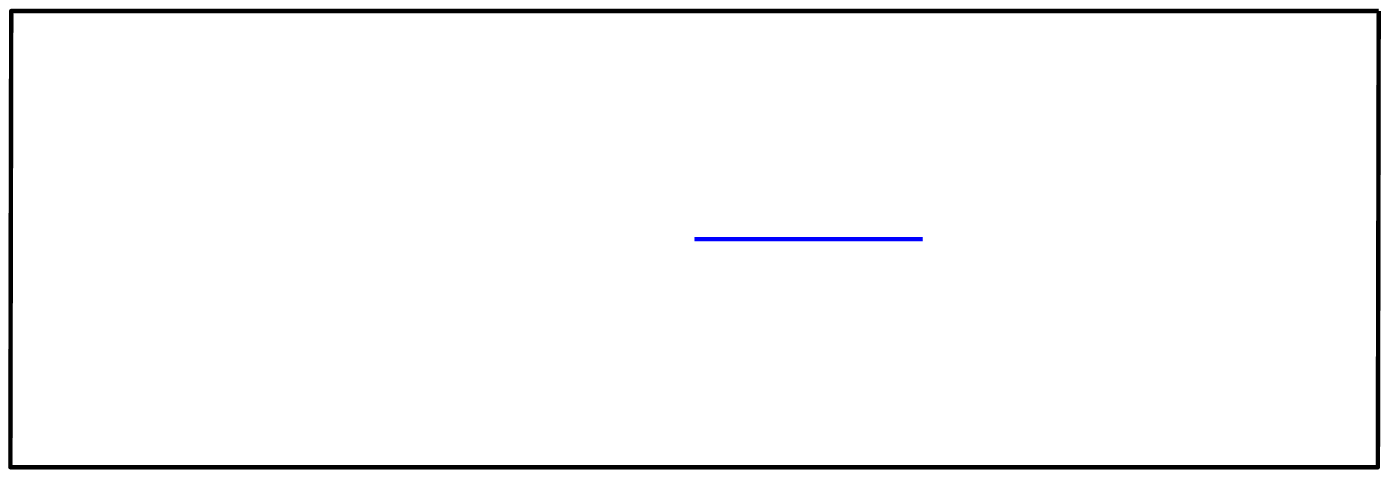}}
\hfill\scalebox{0.4}{\includegraphics[trim=0cm 0cm 0cm 0cm,clip]{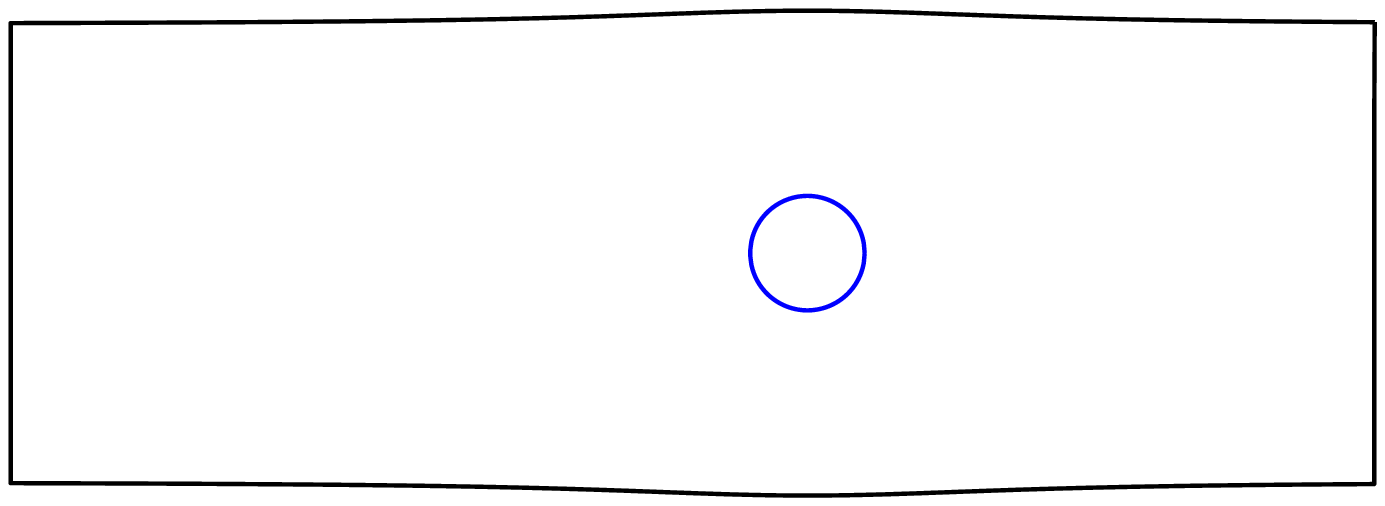}}
}
\caption{The rectangle with a slit domain $G$ (left) and its image $\Psi_2(G)$ (right) as described in~\S\ref{sec:rs}.}
\label{fig:rec}
\end{figure}

Consider the domain (see Figure~\ref{fig:rec} (left))
\[
G=\{z\,:\, -3<\Re z<3,\; -1<\Im z<1\}\backslash[-1,0].
\]
The MATLAB function \verb|Mityuk| is not directly applicable to such a domain as it is not bounded by Jordan curves. So we shall use an elementary conformal mapping to ``open up'' the intervals and obtain a domain bounded by piecewise smooth Jordan curves. Since the elementary mapping
\[
\Psi_1(z)=\frac{1}{4}\left(z+\frac{1}{z}\right)+\frac{1}{2}
\]
maps conformally the exterior of the unit disk onto the exterior of the segment $[0,1]$, its inverse mapping
\[
\Psi_2(z)=\Psi_1^{-1}(z)
=(2z-1)\left(1+\sqrt{1-\frac{1}{(2z-1)^2}}\right)
\]
maps the exterior of the segment $[0,1]$ onto the exterior of the unit disk where we choose the branch for which $\sqrt{1}=1$. Thus, the mapping $\zeta=\Psi_1^{-1}(z)$ maps the domain $G$ onto the domain $\Psi_2(G)$ exterior to the unit disk and interior to the piecewise smooth Jordan curve $\Psi_2(\Gamma_0)$ where $\Gamma_0$ is the external boundary component of $G$ (see Figure~\ref{fig:rec} (right)).

 Let $w=\Psi(\zeta)$ be the conformal mapping from the domain $\Psi_2(G)$ onto the unit disk with a circular/radial slit such that $\Psi(\Psi_2(\alpha))=0$ and $\Psi'(\Psi_2(\alpha))>0$. Hence, the function
\[
w=\Phi_\alpha(z)=e^{-\i\arg\Psi_2(\alpha)}\Psi(\Psi_2(z))
\]
conformally maps the domain $G$ onto the unit disk with a circular/radial slit such that
\[
\Phi_\alpha(\alpha)=0, \quad \Phi_\alpha'(\alpha)=\Psi'(\Psi_2(\alpha))|\Psi'_2(\alpha)|>0,
\]
where
\[
\Psi'_2(\alpha)=\frac{1}{\Psi'_1(\Psi_2(\alpha))}=\frac{4\Psi_2(\alpha)^2}{\Psi_2(\alpha)^2-1}.
\]
Henceforth, in view of~\eqref{eq:R}, Mityuk's radius of the domain $G$ with respect to the point $\alpha$ is given by
$$
R(G,\alpha)=\frac{1}{\Psi'(\Psi_2(\alpha))|\Psi'_2(\alpha)|}.
$$
Note that $R(\Psi_2(G),\Psi_2(\alpha))=1/\Psi'(\Psi_2(\alpha))$, and hence the previous equation becomes (see~\cite[Corollary 2.2.1]{Vas02})
\begin{equation}\label{eq:RG-2}
R(G,\alpha)=\frac{R(\Psi_2(G),\Psi_2(\alpha))}{|\Psi'_2(\alpha)|}.
\end{equation}

\begin{figure}[t] %
	\centerline{
		\scalebox{0.5}{      \includegraphics[trim=0cm 0cm 0cm 0cm,clip]{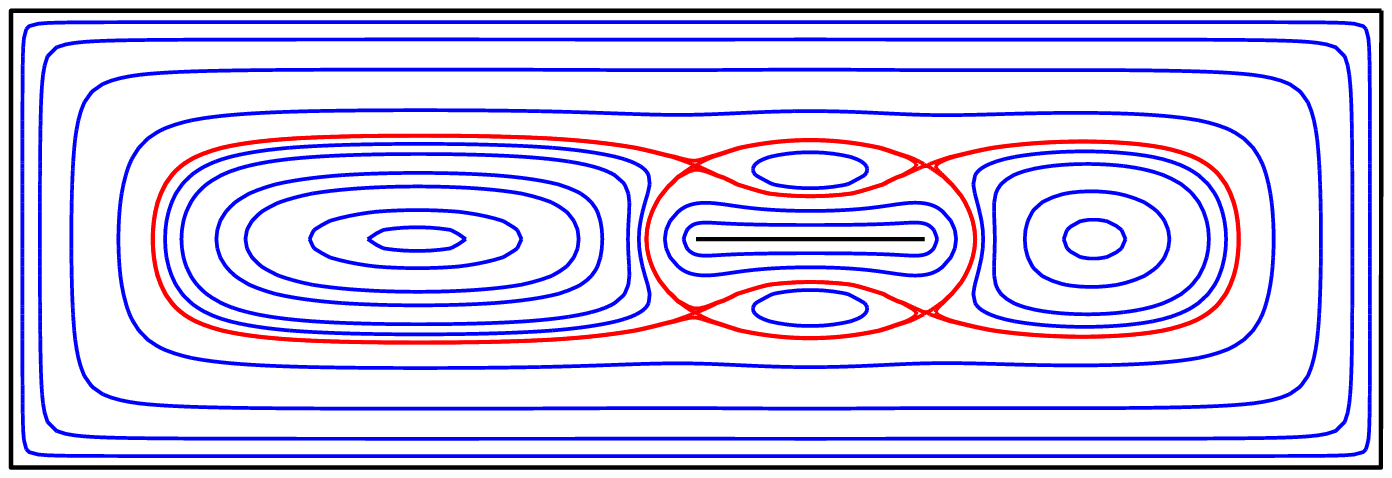}}
		\hfill\scalebox{0.5}{\includegraphics[trim=0cm 0cm 0cm 0cm,clip]{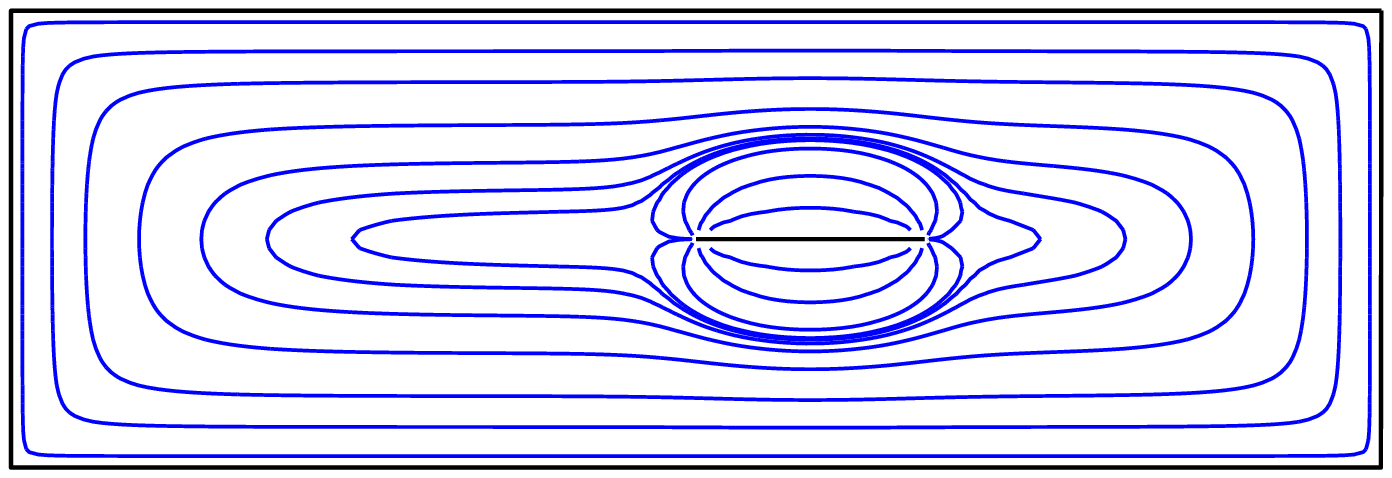}}
	}
	\caption{The contour maps of the function $R(G,\alpha)$ for $\theta_1=\pi/2$ (left) and $\theta_1=0$ (right).  Critical streamlines are shown in red color.}
	\label{fig:rec-L}
\end{figure}

\begin{figure}[t] %
	\centerline{
		\scalebox{0.35}{      \includegraphics[trim=0cm 0cm 0cm 0cm,clip]{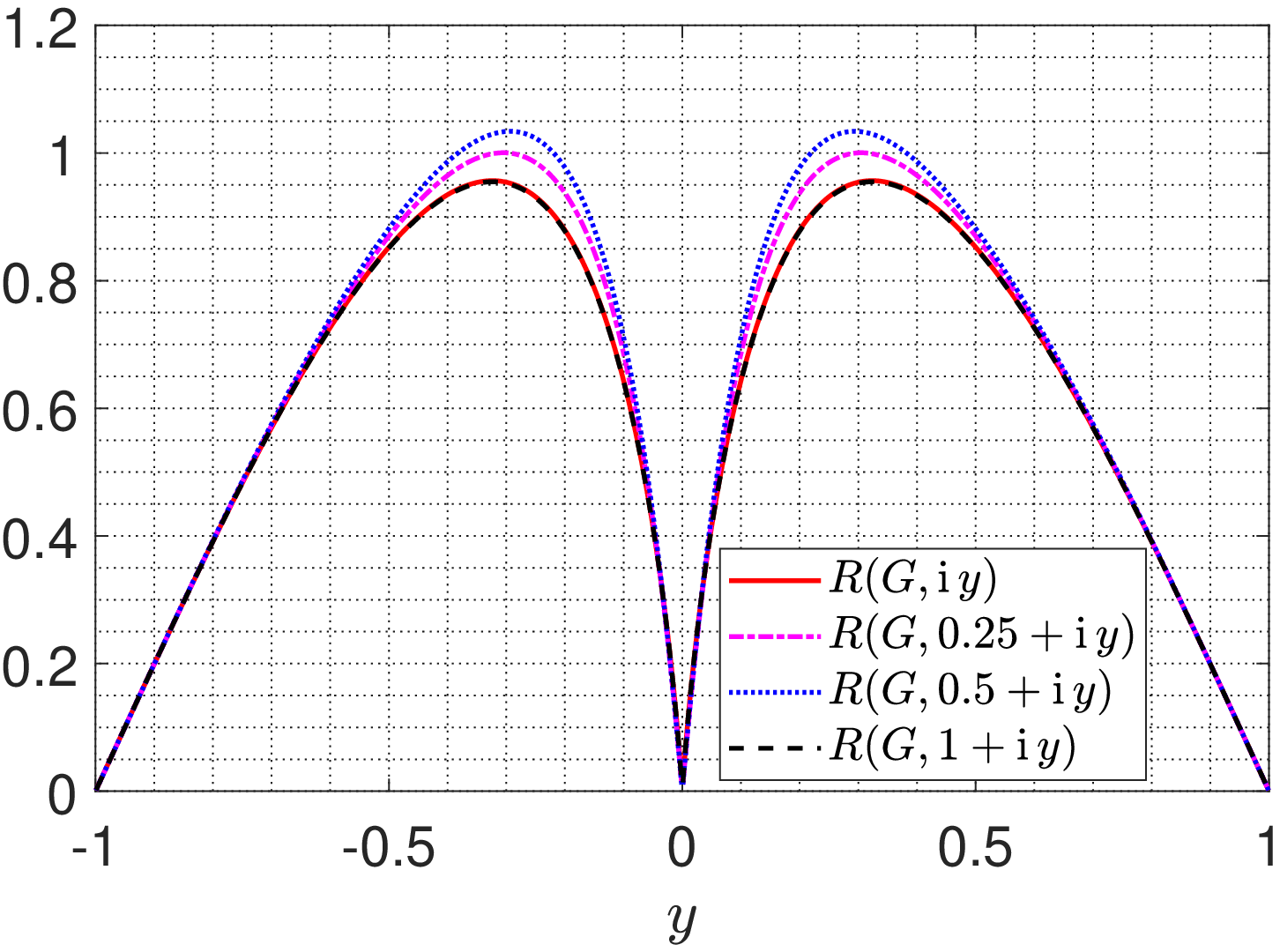}}
		\hfill\scalebox{0.35}{\includegraphics[trim=0cm 0cm 0cm 0cm,clip]{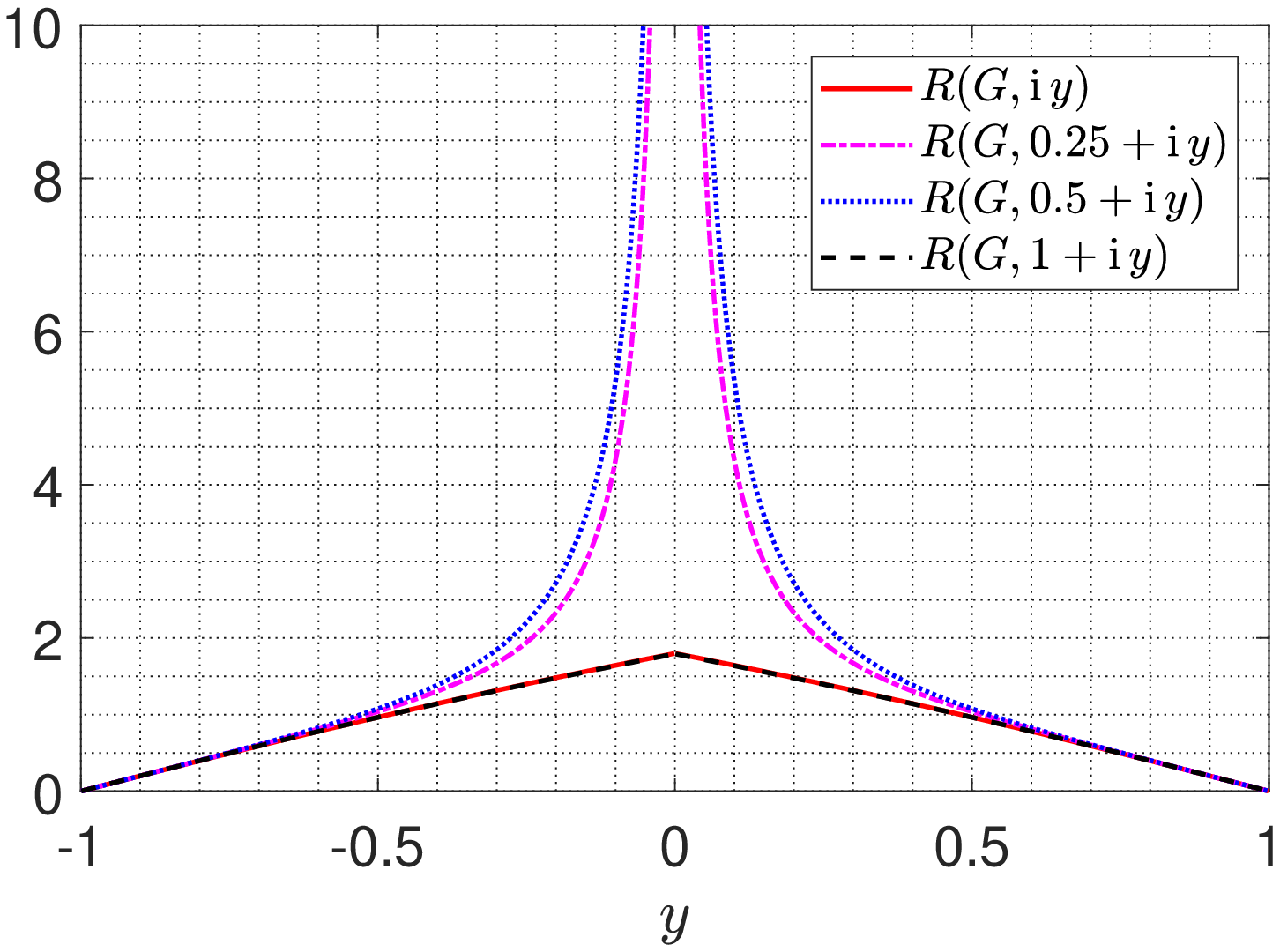}}
		\hfill\scalebox{0.35}{\includegraphics[trim=0cm 0cm 0cm 0cm,clip]{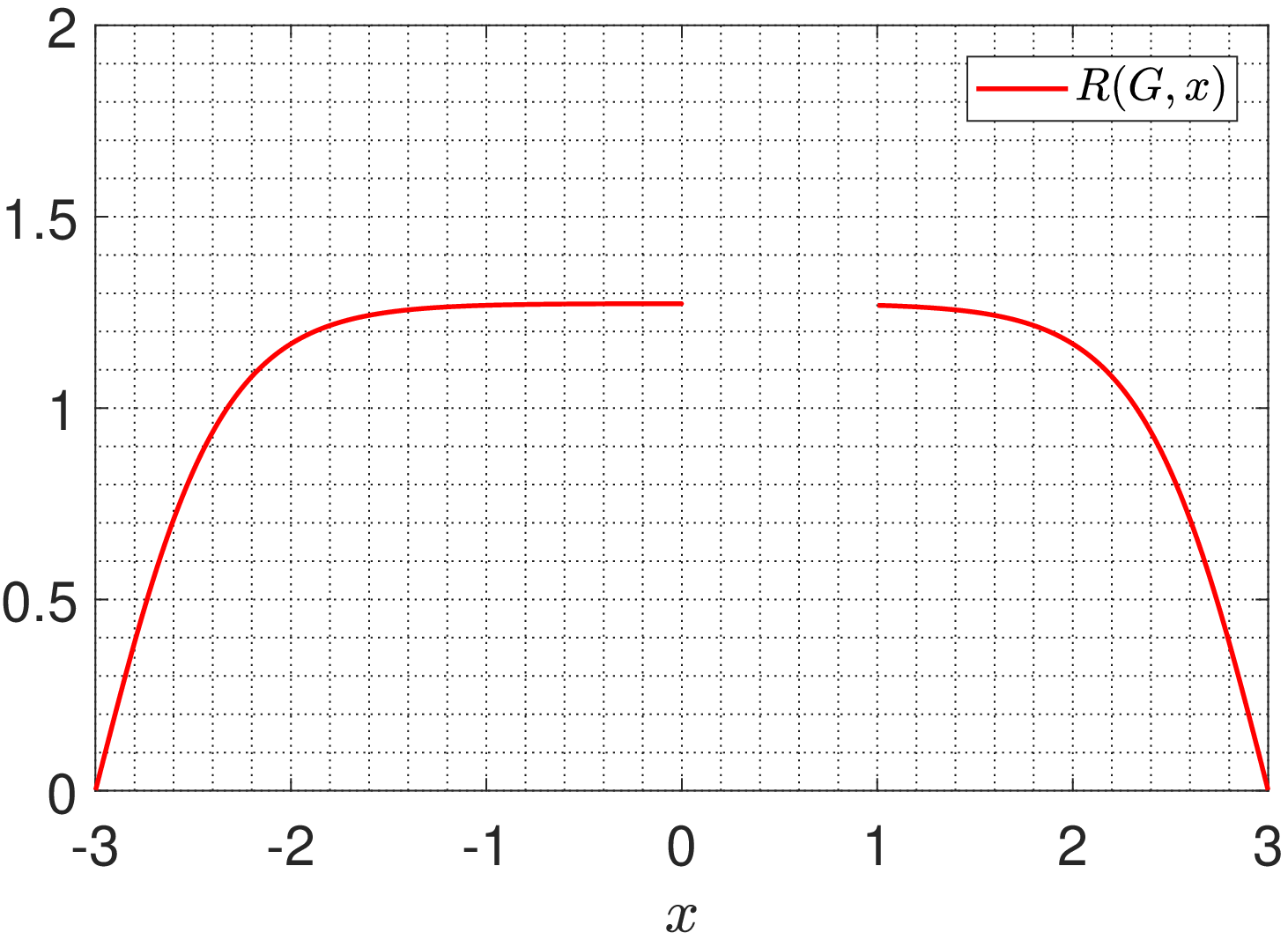}}
	}
	\caption{The values of $R(G,\i y)$, $R(G,0.25+\i y)$, $R(G,0.5+\i y)$, and $R(G,1+\i y)$, $x\in(-1,0)\cup(0,1)$, for $\theta_1=\pi/2$ (left) and $\theta_1=0$ (center), and the values of $R(G,x)$, $x\in(-3,0)\cup(1,3)$,  for $\theta_1=0$ (right). The domain $G$ is the same as in Figure~\ref{fig:rec-L}.}
	\label{fig:rec-x}
\end{figure}

The external boundary component of the domain $\Psi_2(G)$ is a piecewise smooth curve. With suitable parametrization of such a boundary component (see~\cite{LSN} for more details), we can use the MATLAB function \verb|Mityuk| to compute the values of Mityuk's radius $R(\Psi_2(G),\Psi_2(\alpha))$, and then get the values of $R(G,\alpha)$ by~\eqref{eq:RG-2}.
The contour maps of the function $R(G,\alpha)$ are shown in Figure~\ref{fig:rec-L}, and unveil the existence of eight critical points ($n_m=n_s=4$) in the case of a circular slit and no critical point for a radial slit.

The function \verb|Mityuk| is also used to compute $R(G,\i y)$, $R(G,0.25+\i y)$, $R(G,0.5+\i y)$, and $R(G,1+\i y)$ for $x\in(-1,1)\setminus\{0\}$ in the case of the two canonical domains, and $R(G,x)$ for $x\in(-3,0)\cup(1,3)$ in the case of a radial slit. 
It should be pointed out that the limits in~\eqref{eq:lim-k}
and~\eqref{eq:lim-j} were proved in~\cite{Eli17} only for domains with smooth boundaries. The graphs in Figure~\ref{fig:rec-x} presents a numerical counter example that shows these limits are not valid if the inner slit of $G$ is mapped onto a radial slit. More precisely, for $\theta_1=0$ (i.e., when the canonical domain is the unit disk with a radial slit), if $\alpha\in G$ approaches the slit in the domain $G$, then the values of $\Psi_2(\alpha)$ move toward the unit circle in the domain $\Psi_2(G)$. Hence, by~\eqref{eq:lim-j}, the values of $R(\Psi_2(G),\Psi_2(\alpha))$ converge to $+\infty$. However, when $\alpha\in G$ goes to any of the two end-points of the slit, the values of $|\Psi'_2(\alpha)|$ converge to $+\infty$. Hence, by~\eqref{eq:RG-2}, it is not necessary that the function $R(G,\alpha)$ has an infinite limit in this case. In fact, as $\alpha$ moves vertically toward the slit, Figure~\ref{fig:rec-x} (center) shows that the values of $R(G,\alpha)$ converge to a finite number when $\alpha$ approaches particularly one of the two end-points of the slit and the values of $R(G,\alpha)$ become large when $\alpha$ approaches points on the slit away of its two end-points. In the case when $\alpha$ comes close to the two end-points along the real axis, Figure~\ref{fig:rec-x} (right) shows that the values of $R(G,\alpha)$ converge to finite numbers which are different from the limits in the vertical direction. This means the limit of Mityuk's radius $R(G,\alpha )$ does not exists as $\alpha$ approaches the two end-points of the slit.

On the other hand, for $\theta_1=\pi/2$ (i.e., when the canonical domain is the unit disk with a circular slit), if $\alpha\in G$ approaches the slit in the domain $G$, then the values of $\Psi_2(\alpha)$ come close to the unit circle in the domain $\Psi_2(G)$ and, by~\eqref{eq:lim-j}, the values of $R(\Psi_2(G),\Psi_2(\alpha))$ converge to $0$. Since $\Psi'_2(\alpha)\neq 0$ even for points on the slit, then, by~\eqref{eq:RG-2}, the values of $R(G,\alpha)$ go to $0$ as well if $\alpha\in G$ approaches any point on the slit. This is confirmed by the numerical results presented in Figure~\ref{fig:rec-x} (left).
Finally, we note that the values of $R(G,\alpha)$ converge to $0$ if $\alpha\in G$ approaches any point on the outer rectangle for both the two cases of the canonical domain.

\subsubsection{Rectangle in rectangle}\label{sec:rr}

\begin{figure}[t] %
\centerline{
\scalebox{0.5}{      \includegraphics[trim=0cm 0cm 0cm 0cm,clip]{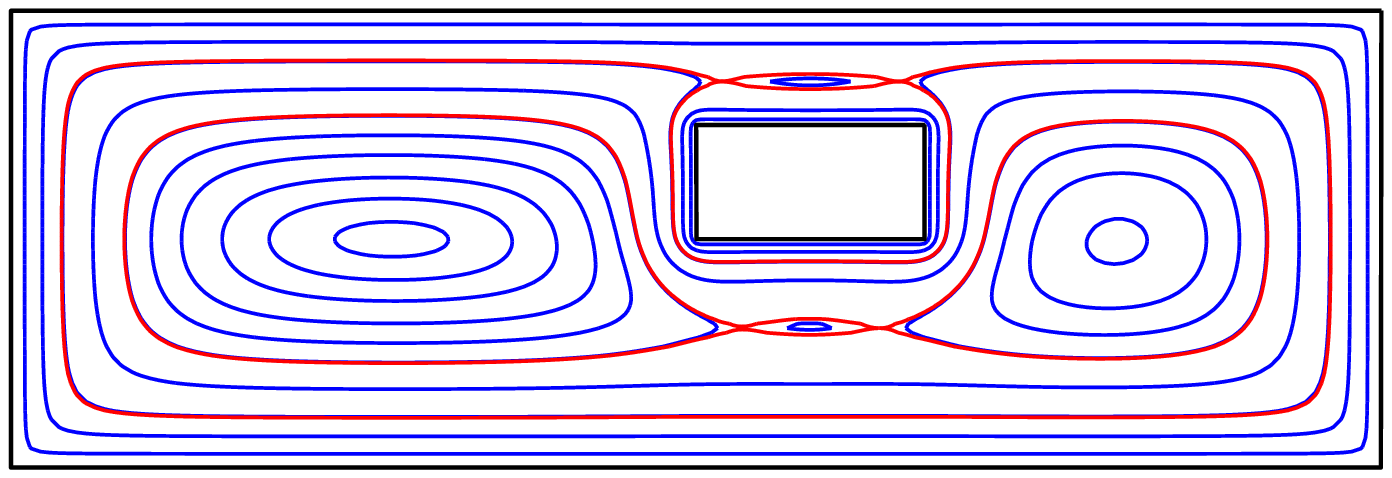}}
\hfill\scalebox{0.5}{\includegraphics[trim=0cm 0cm 0cm 0cm,clip]{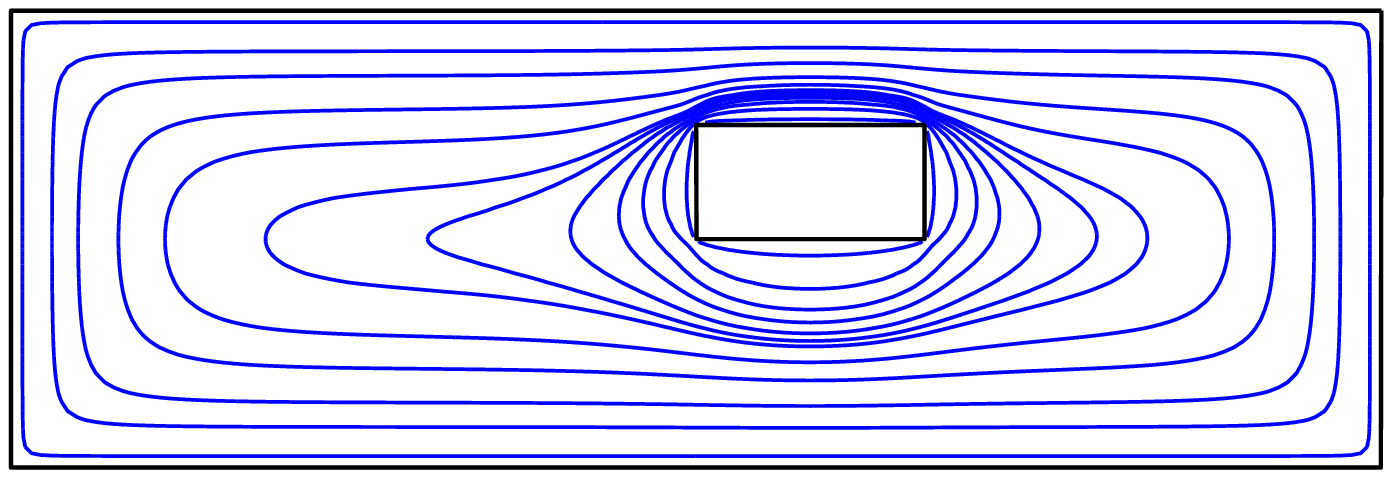}}
}
\caption{The contour maps of the function $R(G,\alpha)$ for $\theta_1=\pi/2$ (left) and $\theta_1=0$ (right).  Critical streamlines are shown in red color.}
\label{fig:recrec-L}
\end{figure}

\begin{figure}[t] %
\centerline{
\scalebox{0.35}{      \includegraphics[trim=0cm 0cm 0cm 0cm,clip]{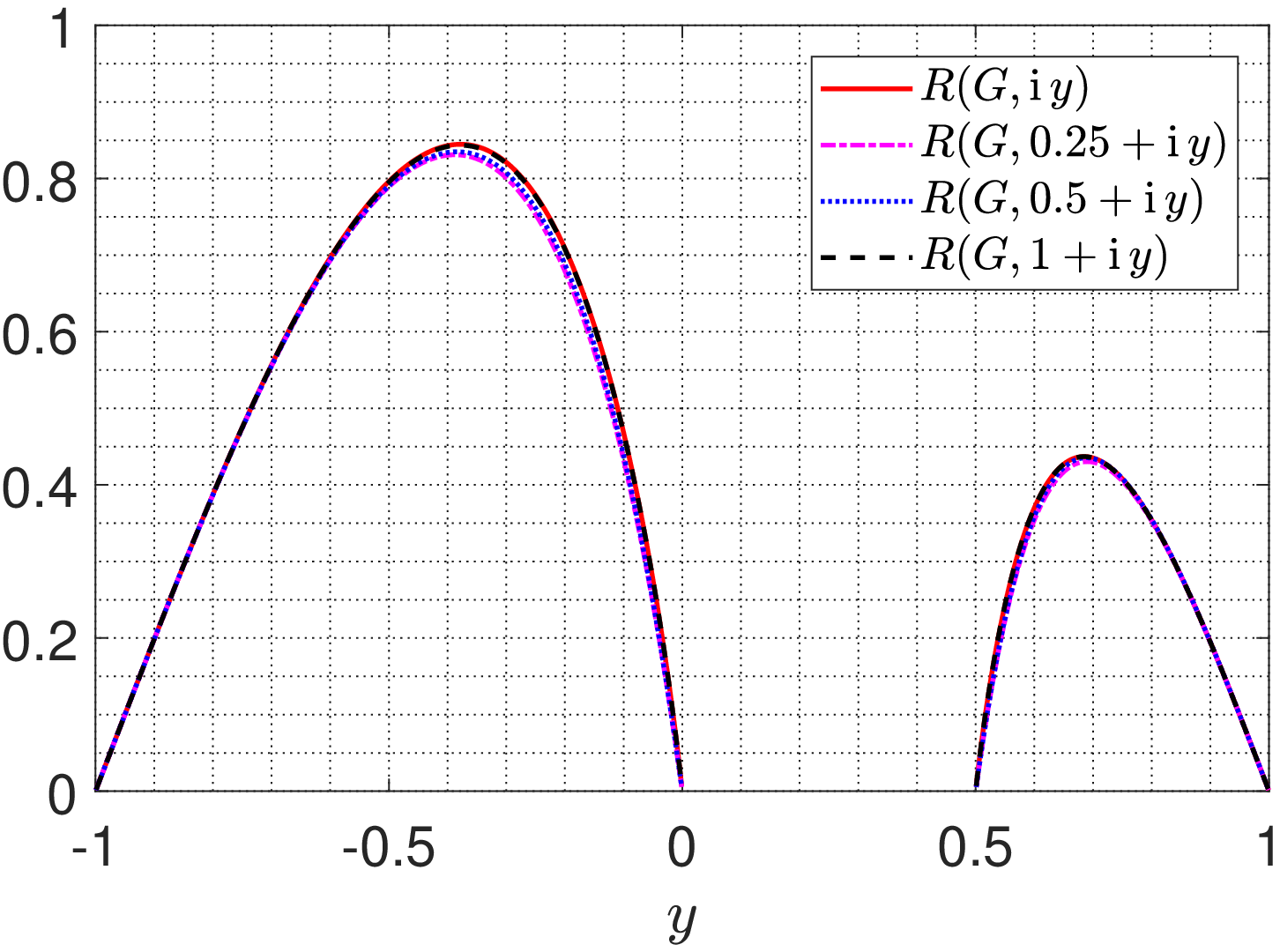}}
\hfill\scalebox{0.35}{\includegraphics[trim=0cm 0cm 0cm 0cm,clip]{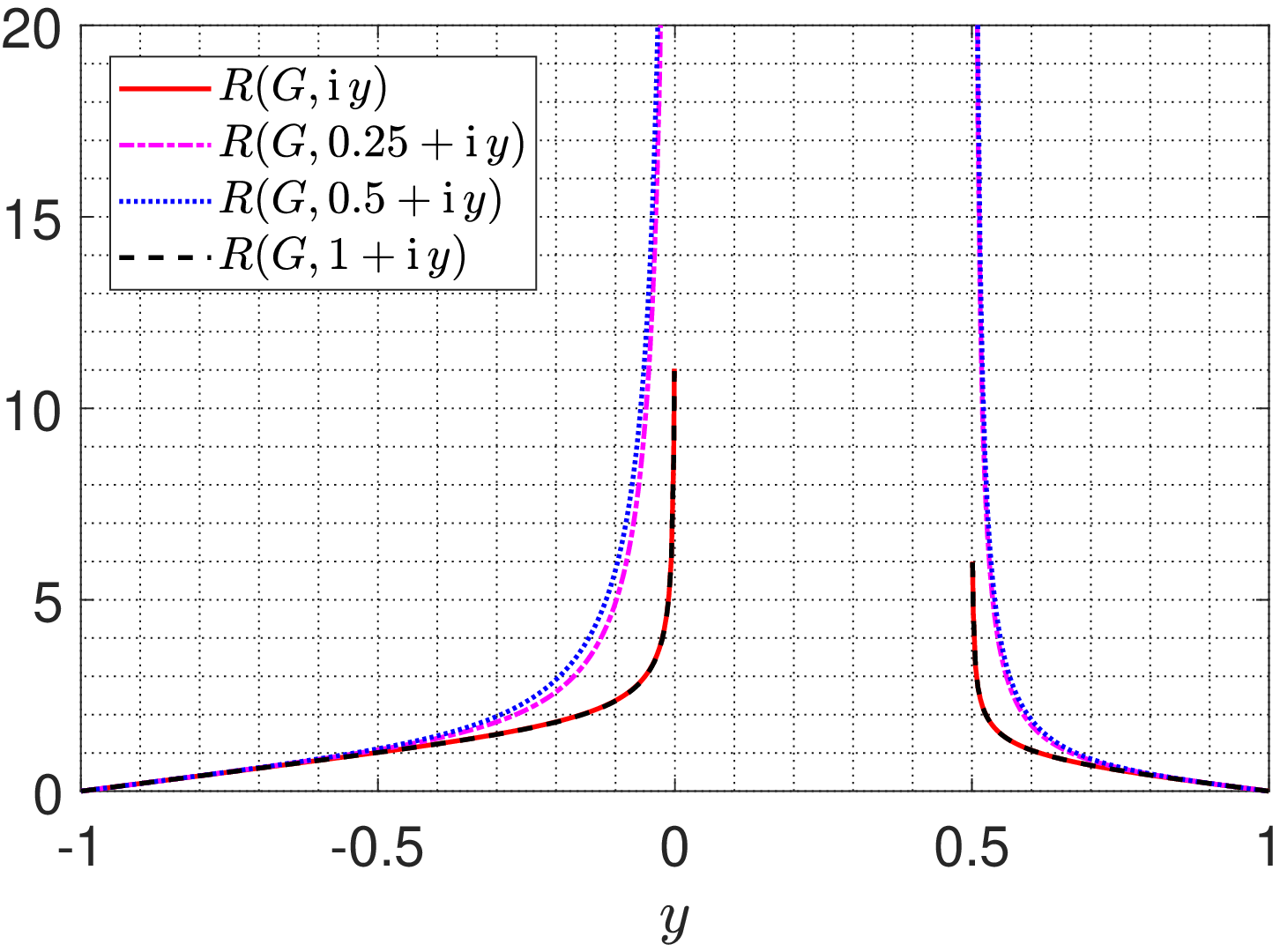}}
\hfill\scalebox{0.35}{\includegraphics[trim=0cm 0cm 0cm 0cm,clip]{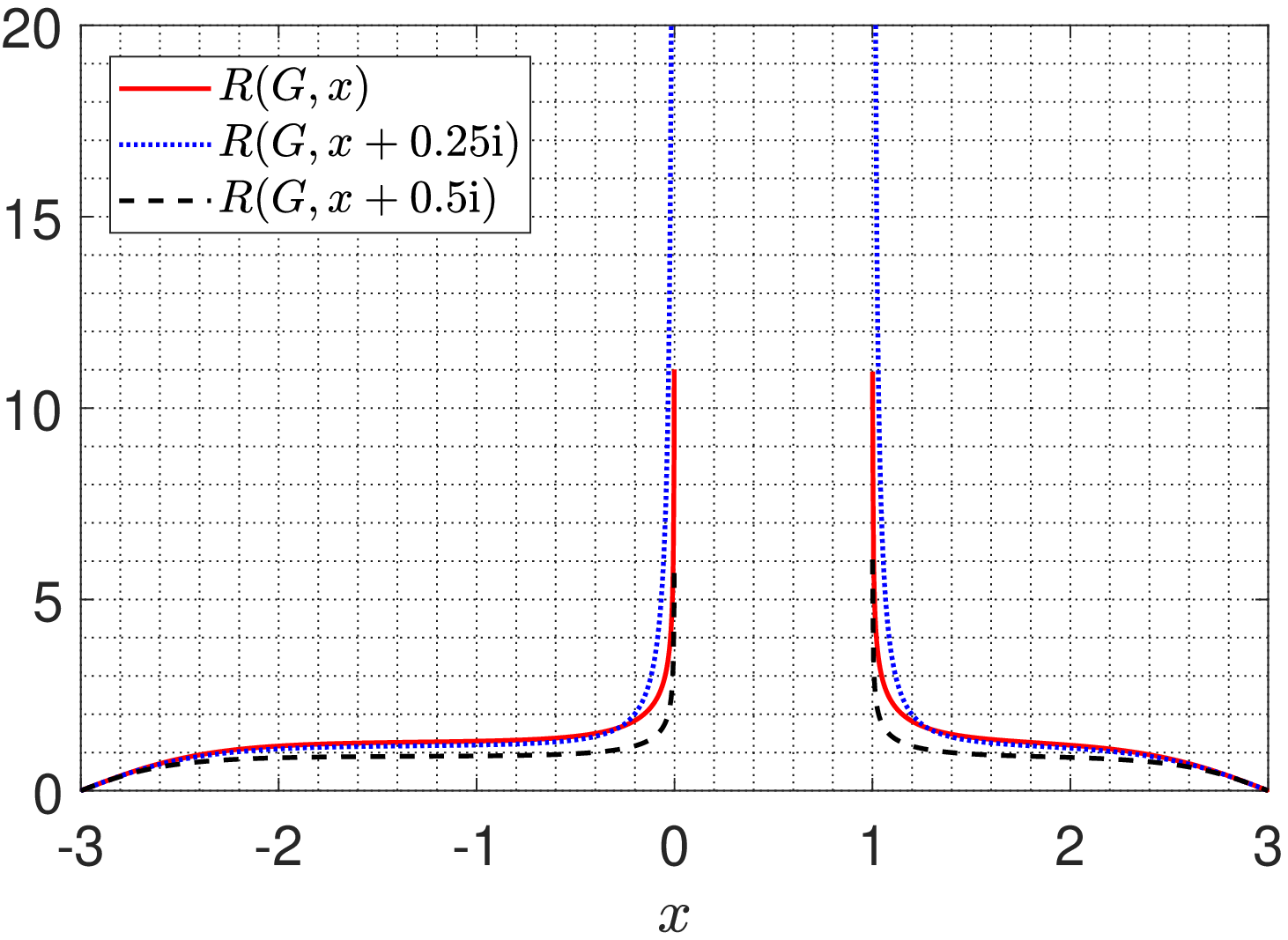}}
}
\caption{The values of $R(G,x)$, $R(G,0.25+\i y)$, $R(G,0.5+\i y)$, and $R(G,1+\i y)$, $y\in(-1,0)\cup(0.5,1)$, for $\theta_1=\pi/2$ (left) and $\theta_1=0$ (center), and the values of $R(G,x)$, $R(x+0.25\i)$, and $R(x+0.5\i)$, $x\in(-3,0)\cup(1,3)$,  for $\theta_1=0$ (right). The domain $G$ is the same as in Figure~\ref{fig:recrec-L}.}
\label{fig:recrec-x}
\end{figure}

Let $$
G=\{z\,:\, -3<\Re z<3,\; -1<\Im z<1\}\setminus\{z\,:\, 0<\Re z<1,\; 0<\Im z<0.5\}.
$$
The function \verb|Mityuk| can be directly used (with suitable parametrization of the boundary components) to compute the values of Mityuk's radius $R(G,\alpha)$.
The contour maps of the function $R(G,\alpha)$  are shown in Figure~\ref{fig:recrec-L}. Similarly to the previous example, we have the existence of eight critical points ($n_m=n_s=4$) for a circular slit and no critical point in the case of a radial slit. We compute also $R(G,\i y)$, $R(G,0.25+\i y)$, $R(G,0.5+\i y)$, and $R(G,1+\i y)$ for $y\in(-1,0)\cup(0.5,1)$. The results are shown in Figure~\ref{fig:recrec-x}.

Although, the boundary components of the domain are only piecewise smooth, the numerical results of this example reveal that the limits in~\eqref{eq:lim-k} and~\eqref{eq:lim-j} are still valid for the two cases of the canonical domain. Figure~\ref{fig:recrec-x} shows the values of $R(G,\alpha)$ converge to $0$ if $\alpha\in G$ approaches any point on the outer rectangle. The results shown in Figure~\ref{fig:recrec-x} (left) suggest also that the limit in~\eqref{eq:lim-j} is valid if the inner rectangle of $G$ is mapped onto a circular slit; i.e., when the canonical domain is the unit disk with a circular slit ($\theta_1=\pi/2$). Indeed, if $\alpha\in G$ approaches any point on the inner rectangle, then the values of $R(G,\alpha)$ converge to $0$. When the canonical domain is the unit disk with a radial slit ($\theta_1=0$), the numerical results presented in Figure~\ref{fig:recrec-x} (center and right) illustrate that if $\alpha\in G$ approaches points on the inner rectangle away of its corners, the values of $R(G,\alpha)$ become very large. However, when $\alpha$ comes vertically or horizontally toward any of the inner rectangle corners, the values of $R(G,\alpha)$ become large but not as large as for off-corner points.

\subsubsection{Triangle in triangle}

\begin{figure}[t] %
\centerline{
\scalebox{0.5}{      \includegraphics[trim=0cm 0cm 0cm 0cm,clip]{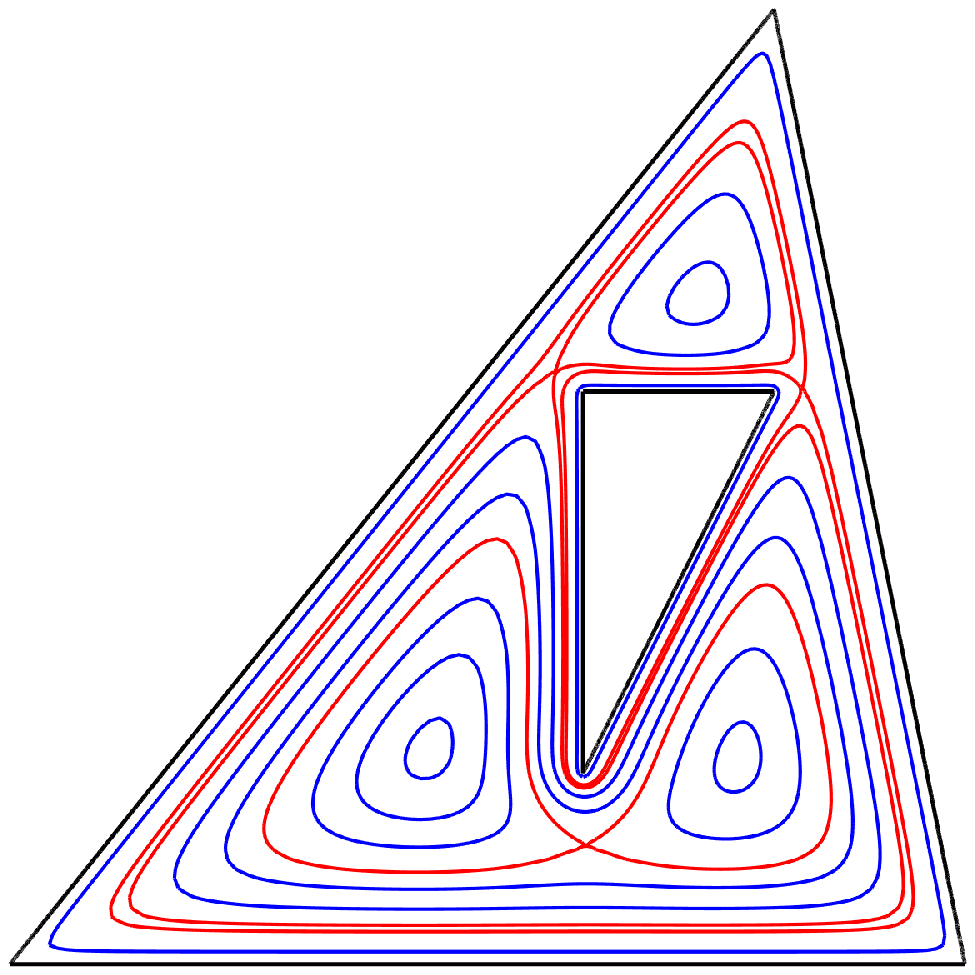}}
\hfill\scalebox{0.5}{\includegraphics[trim=0cm 0cm 0cm 0cm,clip]{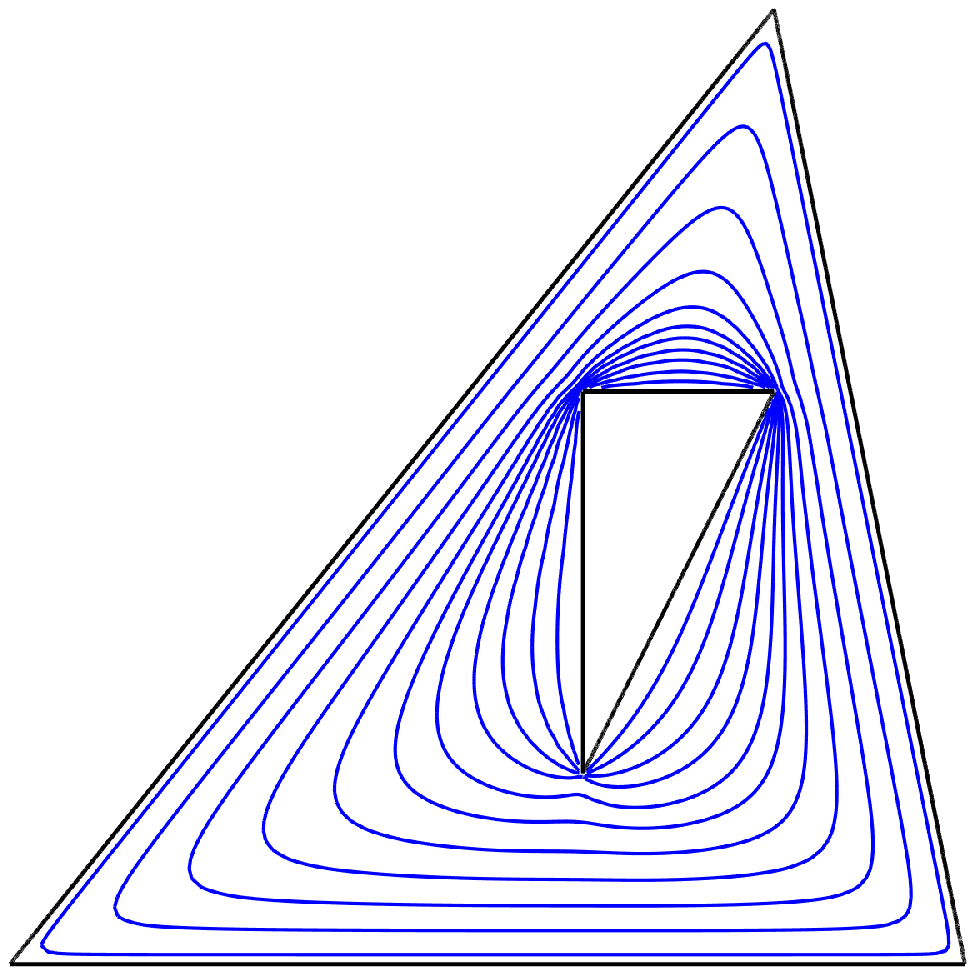}}
}
\caption{The contour maps of the function $R(G,\alpha)$ for $\theta_1=\pi/2$ (left) and $\theta_1=0$ (right).  Critical streamlines are shown in red color.}
\label{fig:tritri-L}
\end{figure}

\begin{figure}[t] %
\centerline{
\scalebox{0.35}{      \includegraphics[trim=0cm 0cm 0cm 0cm,clip]{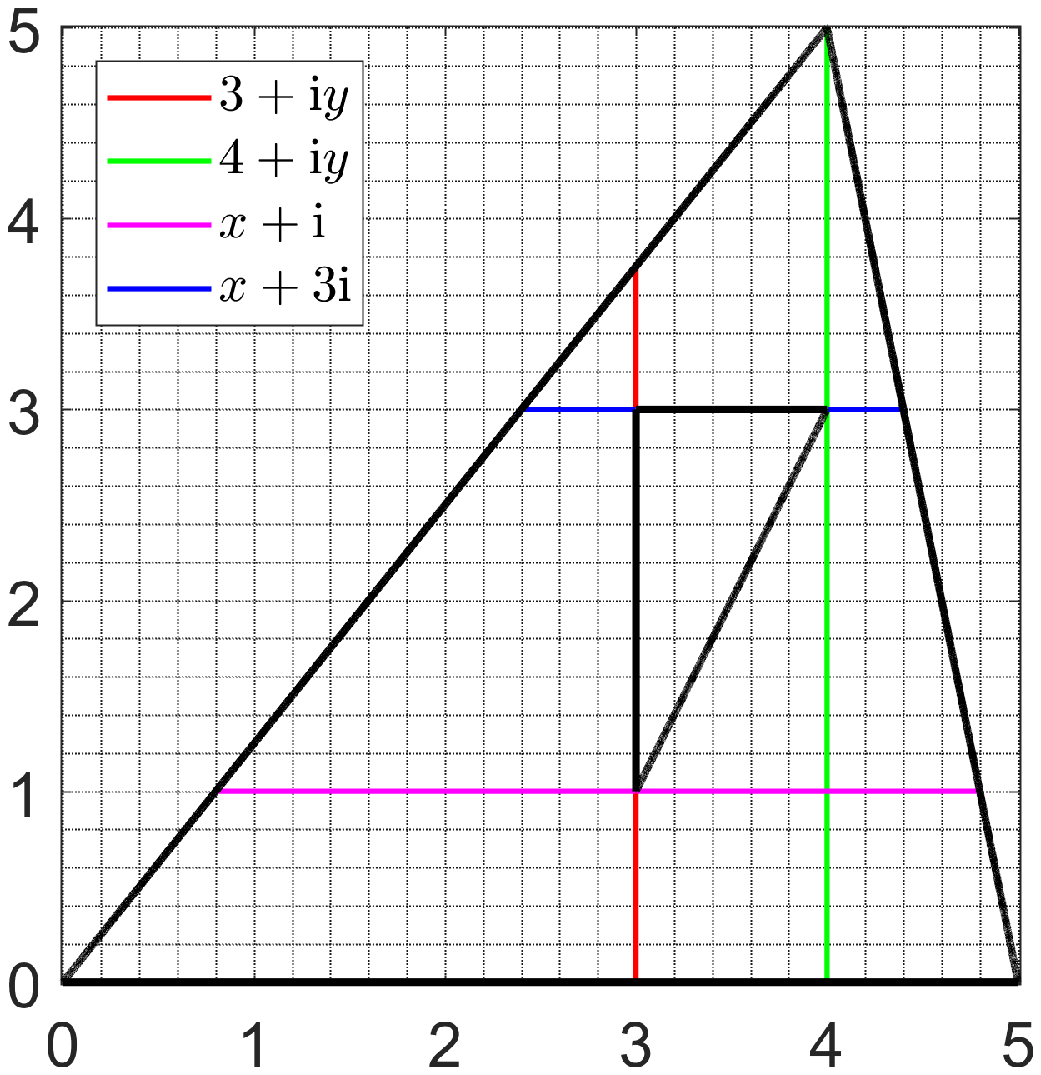}}
\hfill\scalebox{0.35}{\includegraphics[trim=0cm 0cm 0cm 0cm,clip]{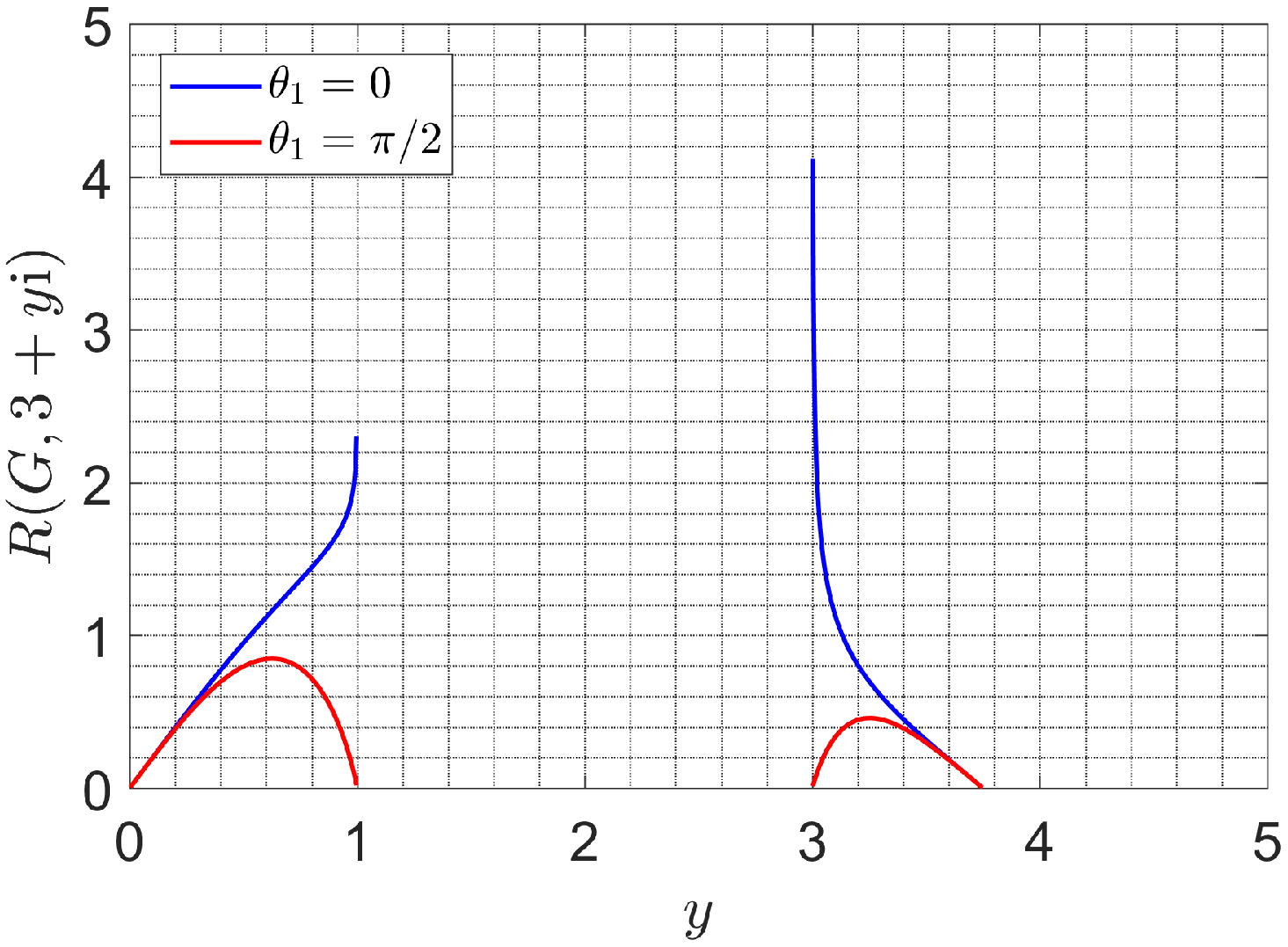}}
\hfill\scalebox{0.35}{\includegraphics[trim=0cm 0cm 0cm 0cm,clip]{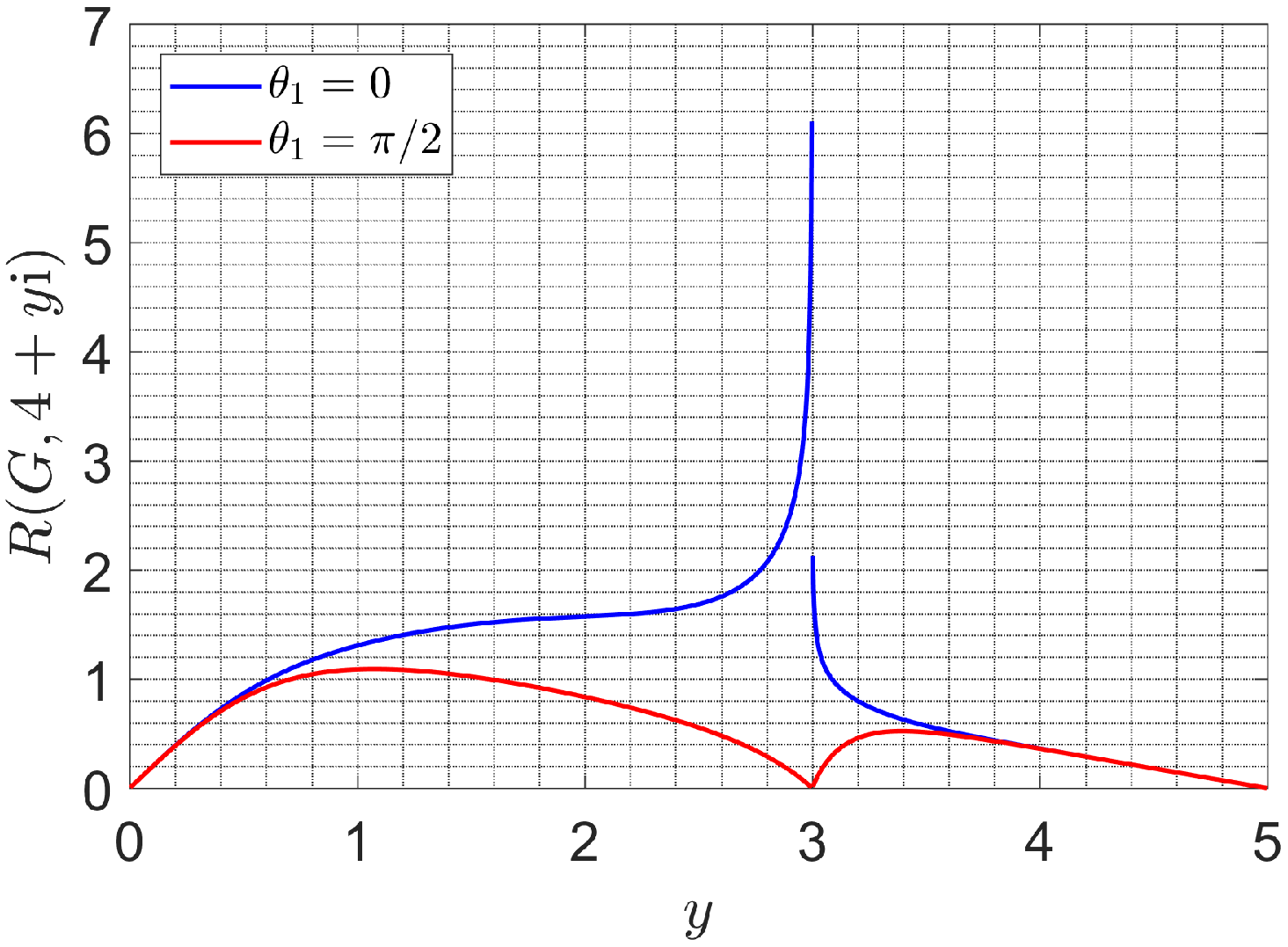}}
}
\centerline{
\hfill\scalebox{0.35}{\includegraphics[trim=0cm 0cm 0cm 0cm,clip]{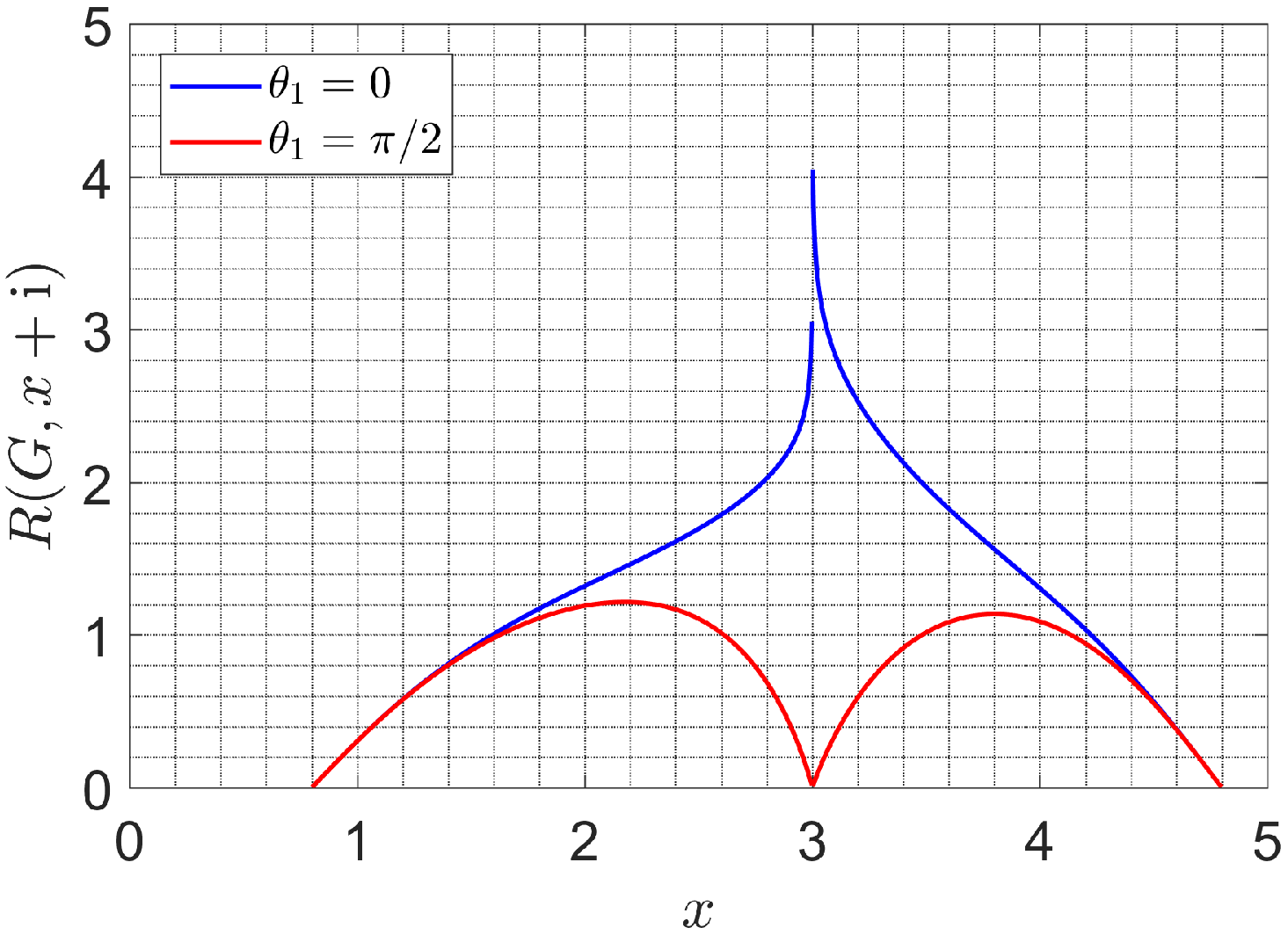}}
\hfill\scalebox{0.35}{\includegraphics[trim=0cm 0cm 0cm 0cm,clip]{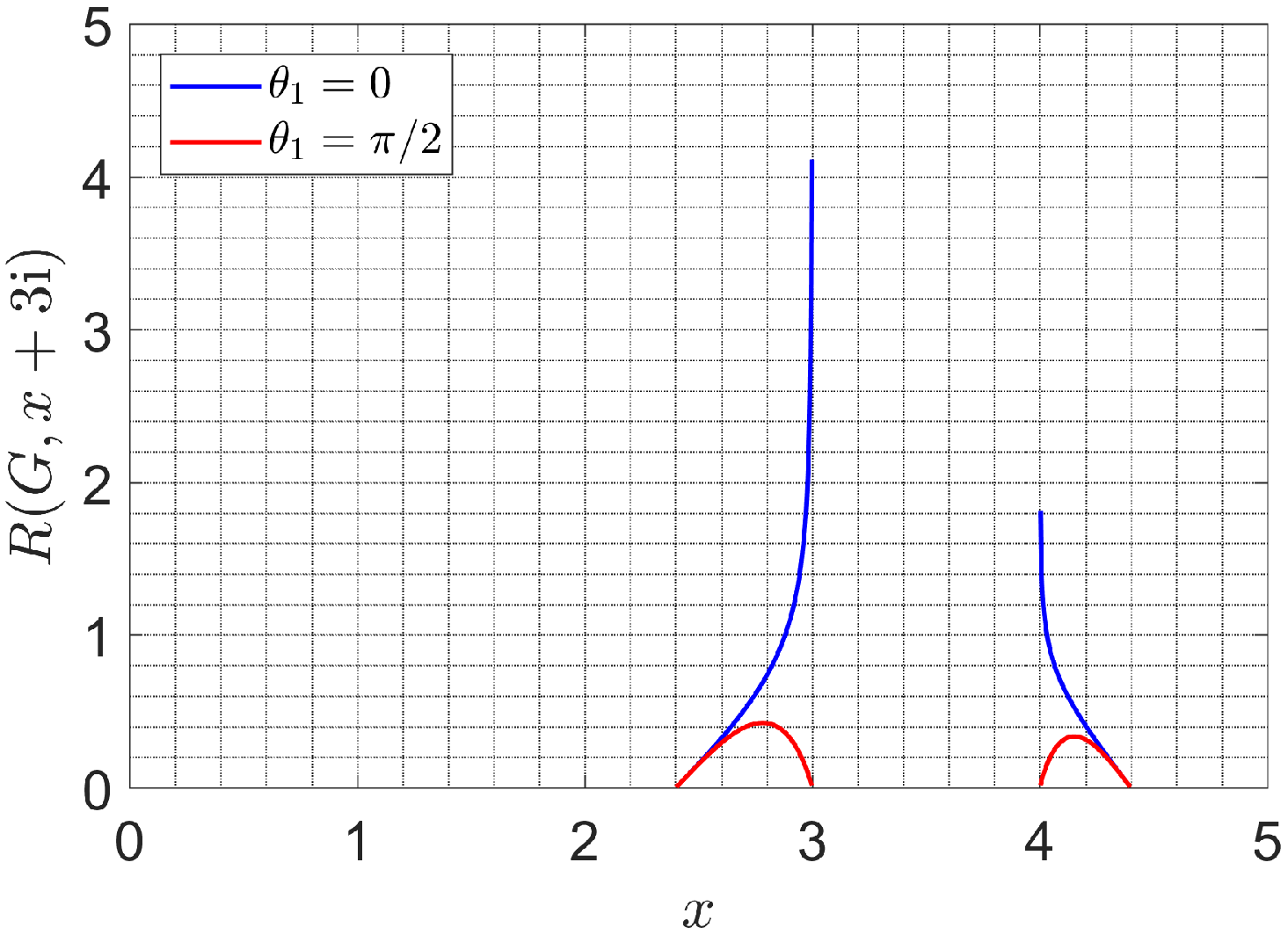}}
\hfill}
\caption{The points $3+\i y$, $4+\i y$, $x+\i$, and $x+3\i$ in the domain $G$ (top, left) and the values of Mityuk's radius $R(G,\alpha)$ at these points.}
\label{fig:tritri-x}
\end{figure}

Consider the domain $G$ bordered externally by the triangle with vertices $0$, $5$, and $4+5\i$, and internally by the triangle with vertices $3+3\i$, $4+3\i$, and $3+\i$ (see Figure~\ref{fig:tritri-x} (top, left)). The contour maps of the function $R(G,\alpha)$  are shown in Figure~\ref{fig:tritri-L}, which unveil the existence of six critical points ($n_m=n_s=3$) in the case of a circular slit and no critical point in the case of a radial slit.
To check the validity of limit values in equations~\eqref{eq:lim-k} and~\eqref{eq:lim-j}, we display in Figure~\ref{fig:tritri-x} the curve of $R(G,\alpha)$ along the vertical and horizontal paths $\alpha=3+\i y$, $\alpha=4+\i y$, $\alpha=x+\i$, and $\alpha=x+3\i$ for real numbers $x$ and $y$ such that $\alpha\in G$. As can be seen from this figure the behavior of the function $R(G,\alpha)$, as $\alpha\in G$ approaches any of the corner points of the inner triangle, is similar to the previous example.

\subsubsection{Square in square / Square in circle / Circle in square}

\begin{figure}[t] %
\centerline{
\scalebox{0.5}{      \includegraphics[trim=0cm 0cm 0cm 0cm,clip]{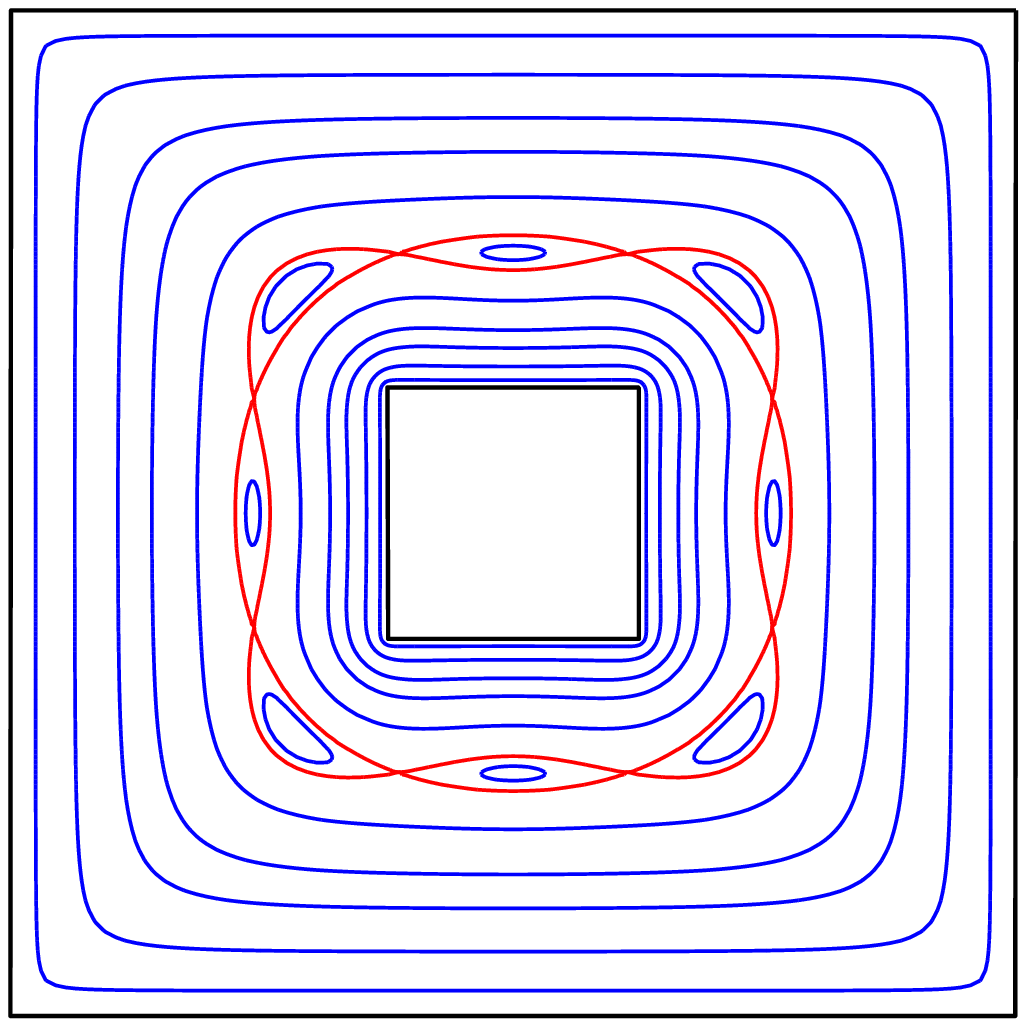}}
\hfill\scalebox{0.5}{\includegraphics[trim=0cm 0cm 0cm 0cm,clip]{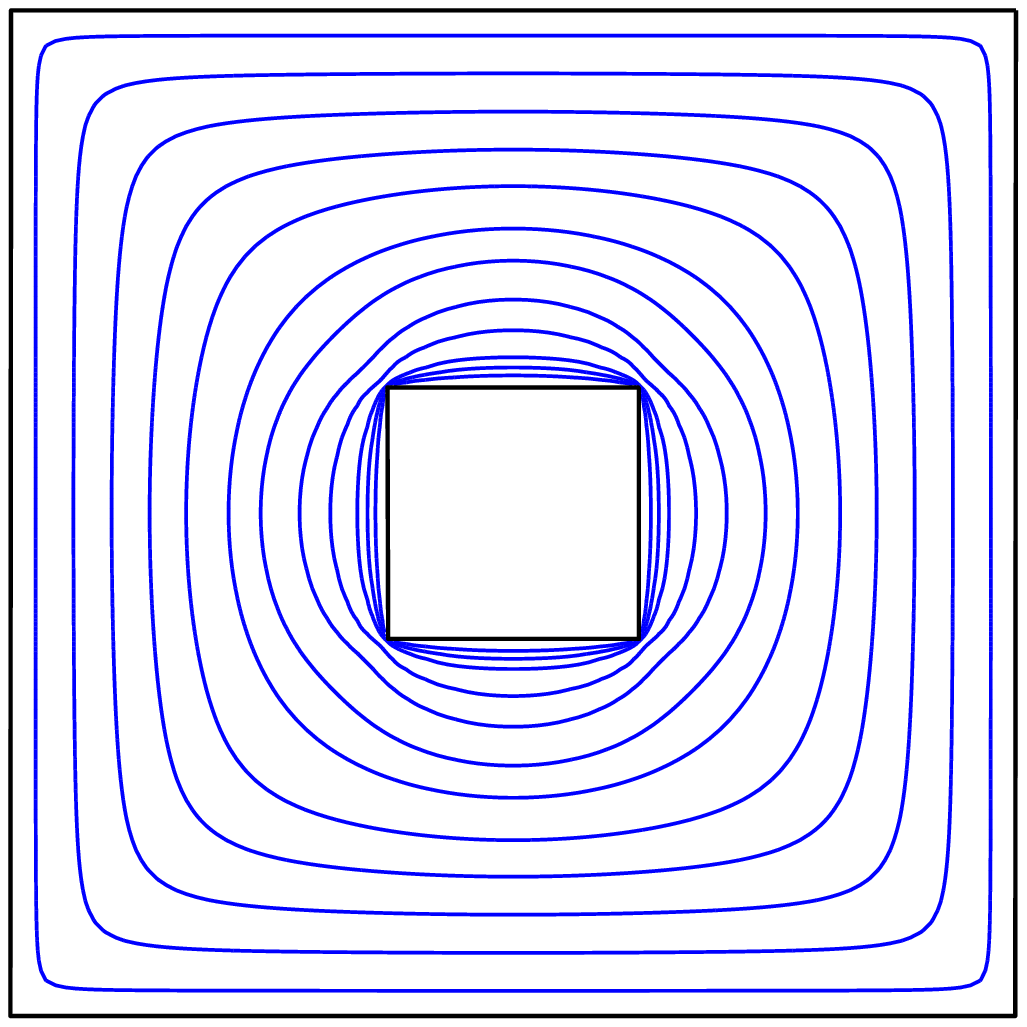}}
}
\centerline{
\scalebox{0.5}{      \includegraphics[trim=0cm 0cm 0cm 0cm,clip]{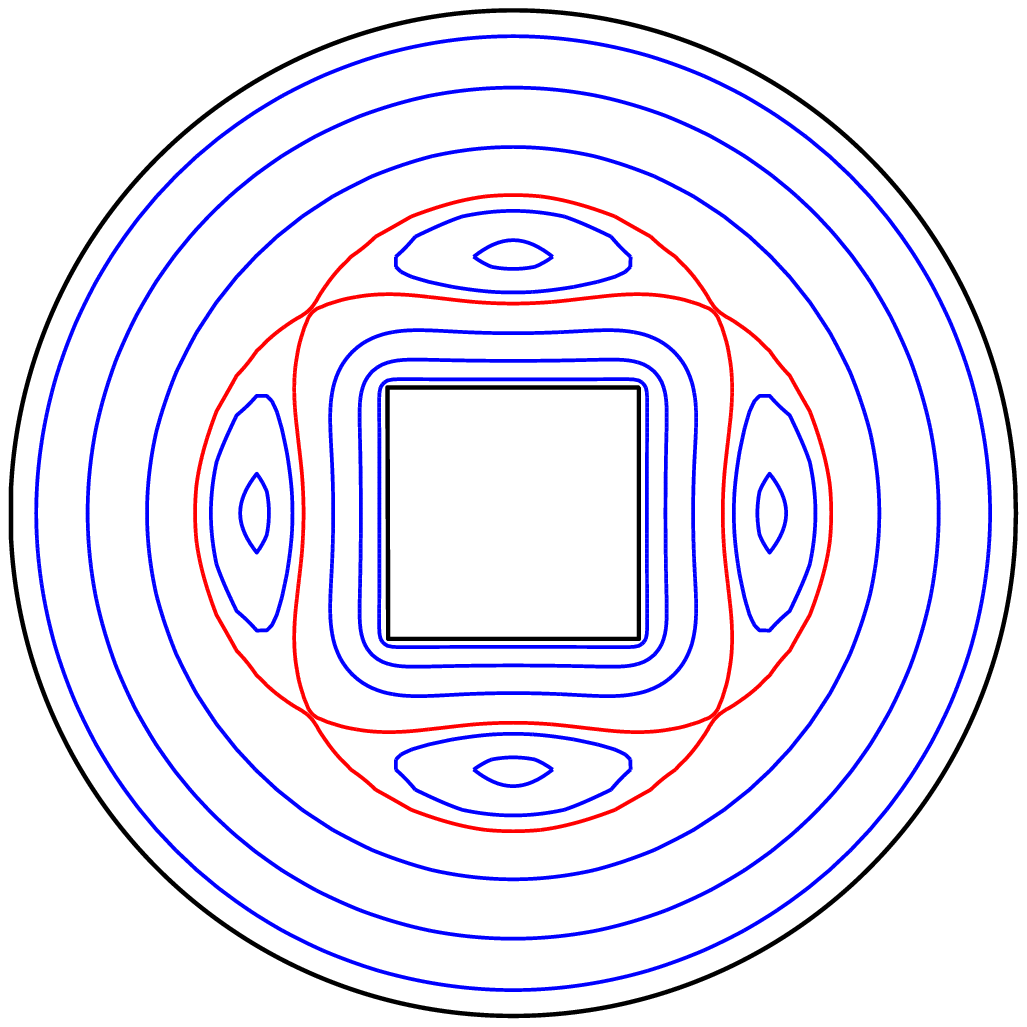}}
\hfill\scalebox{0.5}{\includegraphics[trim=0cm 0cm 0cm 0cm,clip]{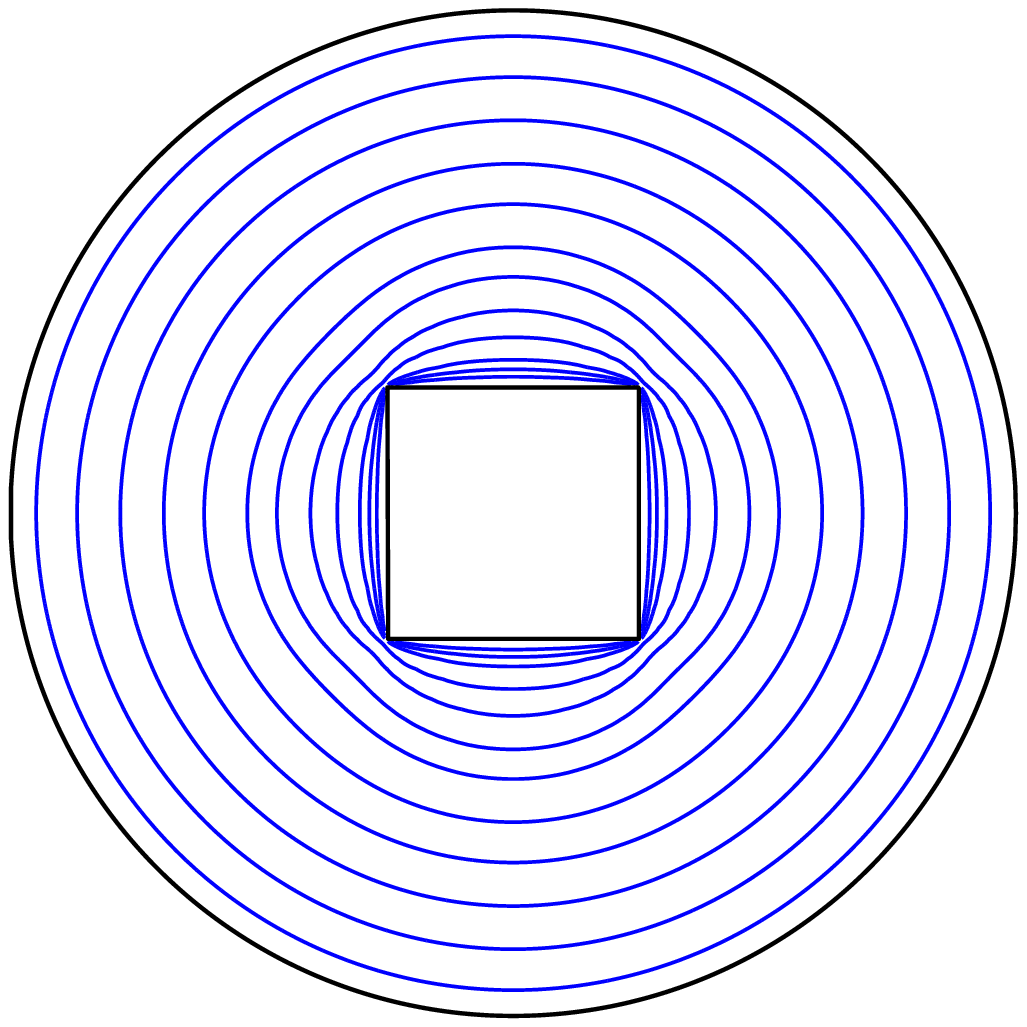}}
}
\centerline{
\scalebox{0.5}{      \includegraphics[trim=0cm 0cm 0cm 0cm,clip]{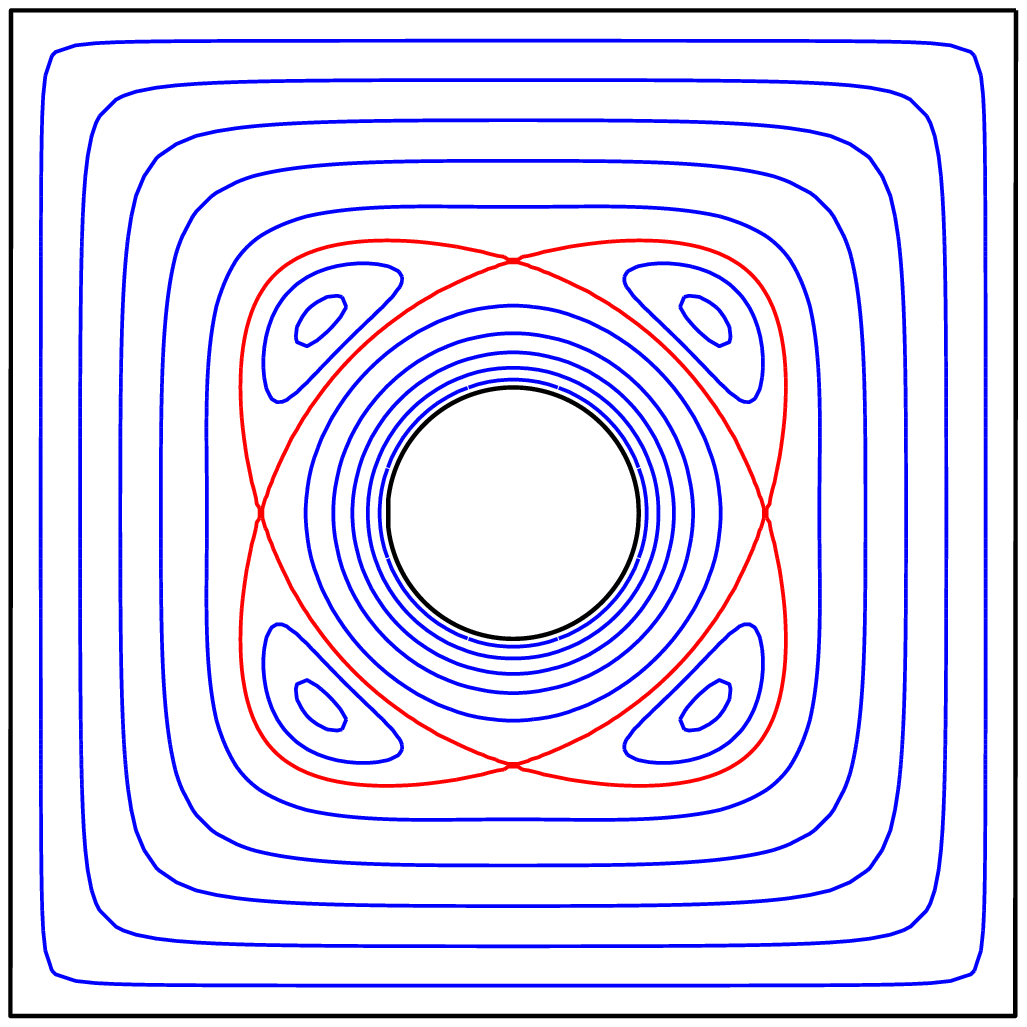}}
\hfill\scalebox{0.5}{\includegraphics[trim=0cm 0cm 0cm 0cm,clip]{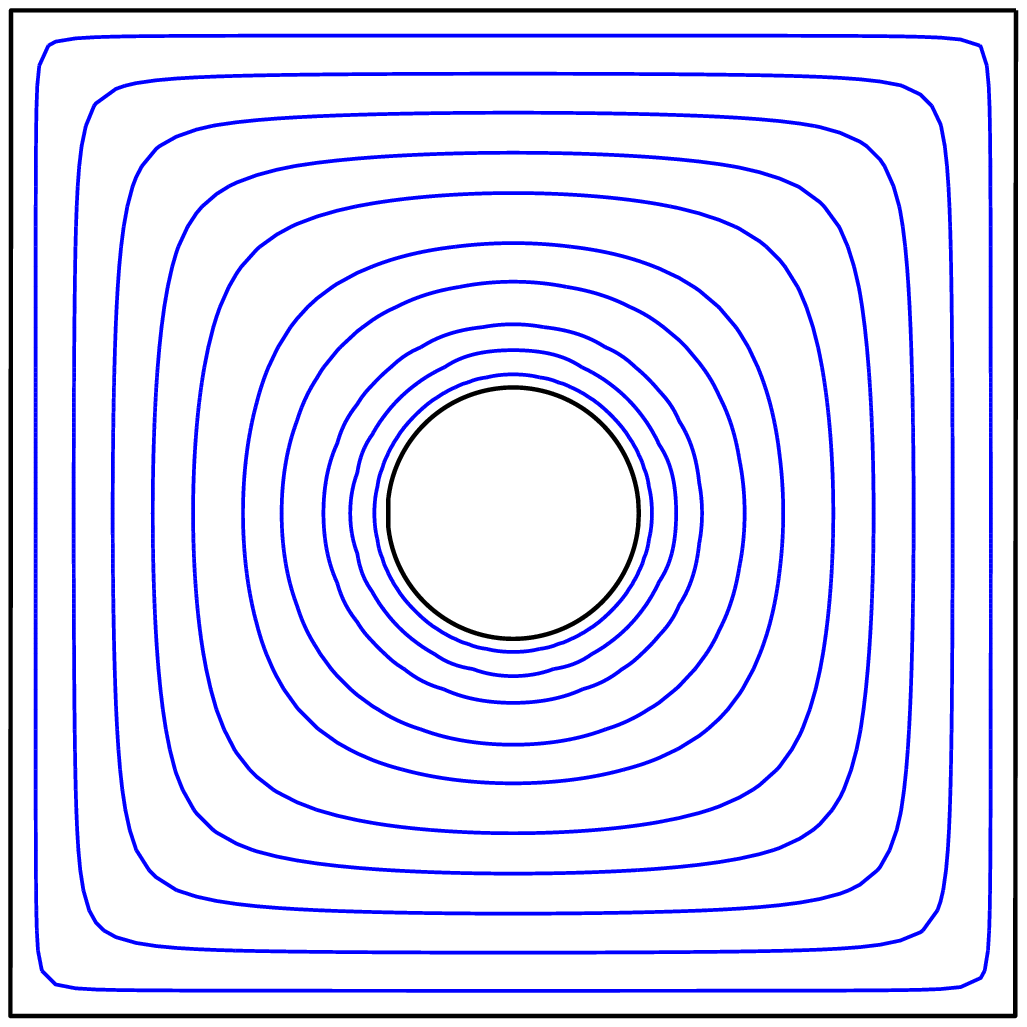}}
}
\caption{The contour maps of the function $R(G,\alpha)$ for $\theta_1=\pi/2$ (left) and $\theta_1=0$ (right).  Critical streamlines are shown in red color.}
\label{fig:an-sq-cr}
\end{figure}

In a similar fashion to the case of an annulus, we consider three examples of $G$ by replacing the inner circle, the outer circle, or both by squares. The vertices of the outer square are $\pm1\pm\i$ and the vertices of the inner square are $\pm0.25\pm0.25\i$ and hence, for these three examples, the domain $G$ is symmetric. The contour maps of the function $R(G,\alpha)$ shown in Figure~\ref{fig:an-sq-cr} are quite revealing in several ways. First, in agreement with the previous examples, Mityuk's radius has no critical point for the canonical domain with a radial slit. In the case of a circular slit, it is apparent that the number of critical points depends on the number of corners in the inner and outer boundary components. In particular, we have the existence of sixteen critical points ($n_m=n_s=8$) when both the inner and outer boundary components are squares, and eight critical points ($n_m=n_s=4$) when one of them is a circle and the other is a square. Note also the symmetrical distribution of critical points location in these three examples.

\subsection{Multiply connected domains}

Finally, in this section, we consider two examples of the domain $G$ when $\ell\geq 2$.

\subsubsection{Three circles}\label{sec:3cr}
\begin{figure}[t] %
\centerline{
\scalebox{0.45}{      \includegraphics[trim=0cm 0cm 0cm 0cm,clip]{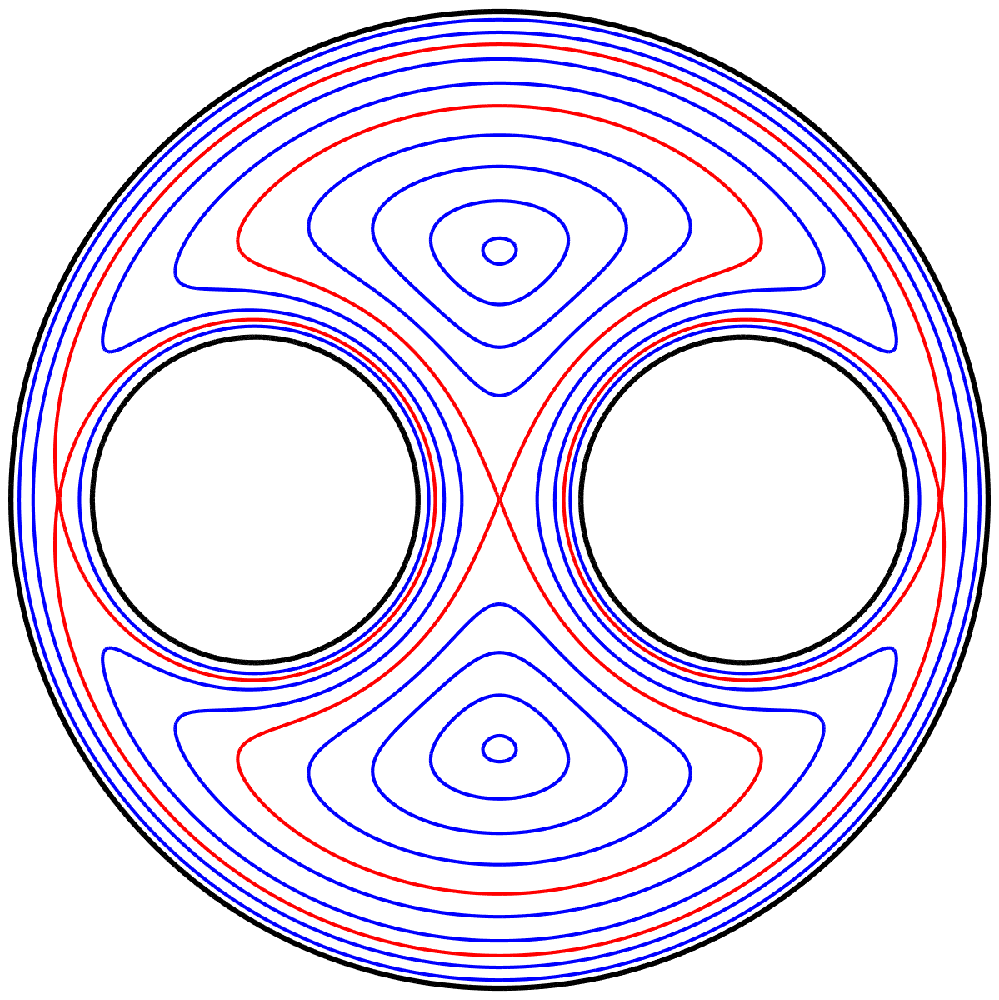}}
\hfill\scalebox{0.45}{\includegraphics[trim=0cm 0cm 0cm 0cm,clip]{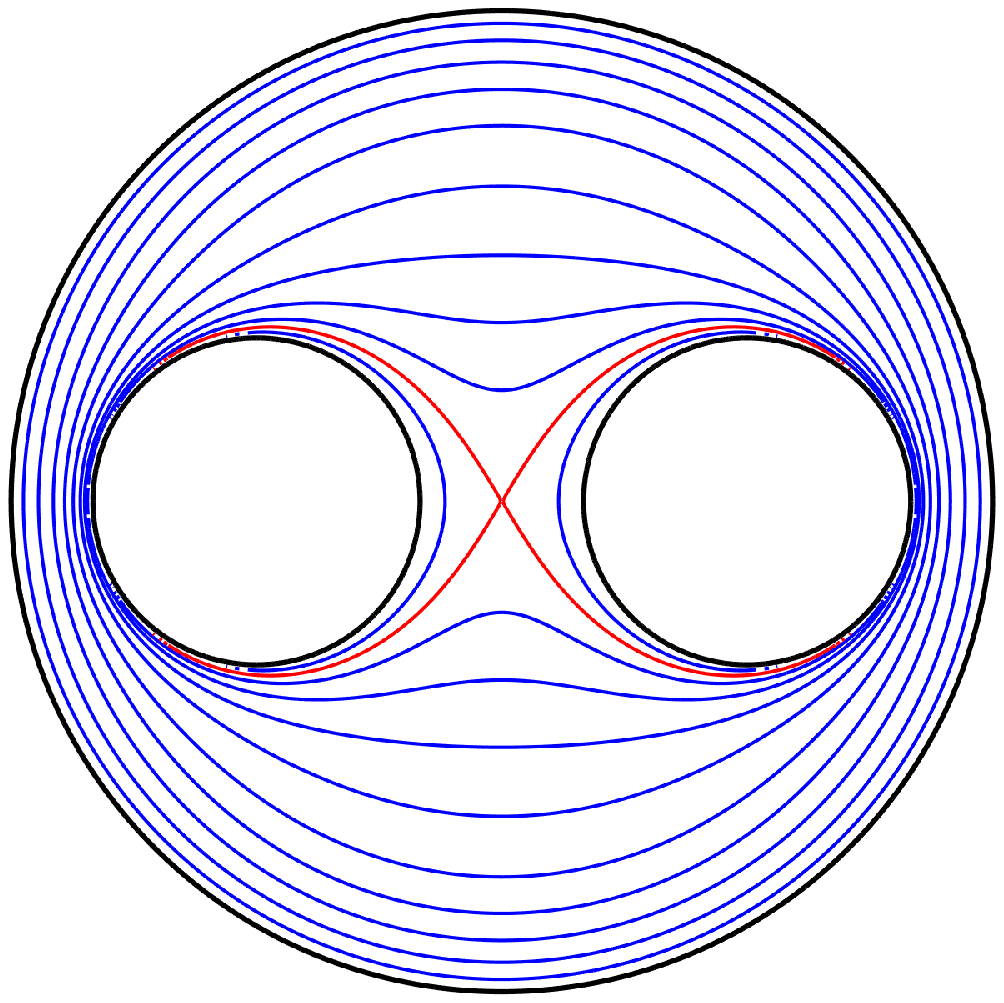}}
\hfill\scalebox{0.45}{\includegraphics[trim=0cm 0cm 0cm 0cm,clip]{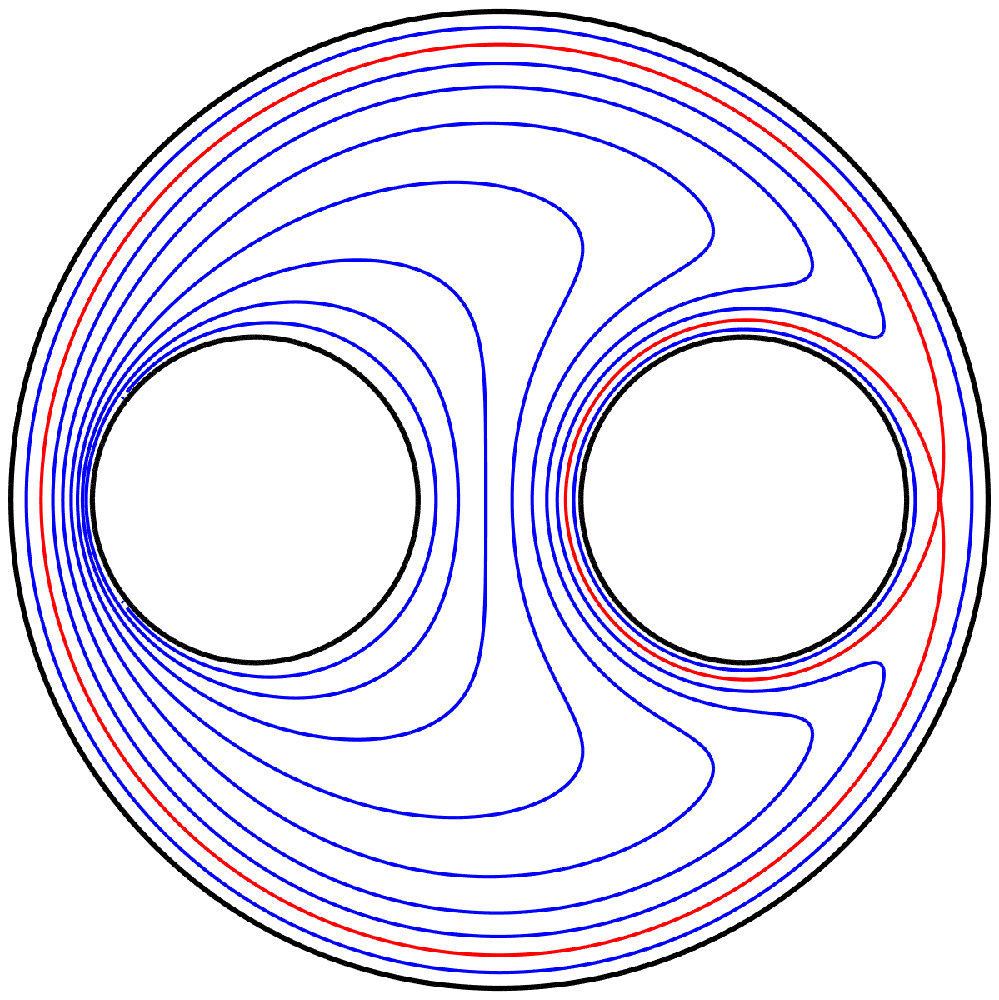}}
}
\caption{The contour maps of the function $R(G,\alpha)$ for: $\theta_1=\theta_2=\pi/2$ (left), $\theta_1=\theta_2=0$ (center), and $\theta_1=\pi/2$, $\theta_2=0$ (right).  Critical streamlines are shown in red color.}
\label{fig:cir3L}
\end{figure}

\begin{figure}[t] %
\centerline{
\scalebox{0.35}{      \includegraphics[trim=0cm 0cm 0cm 0cm,clip]{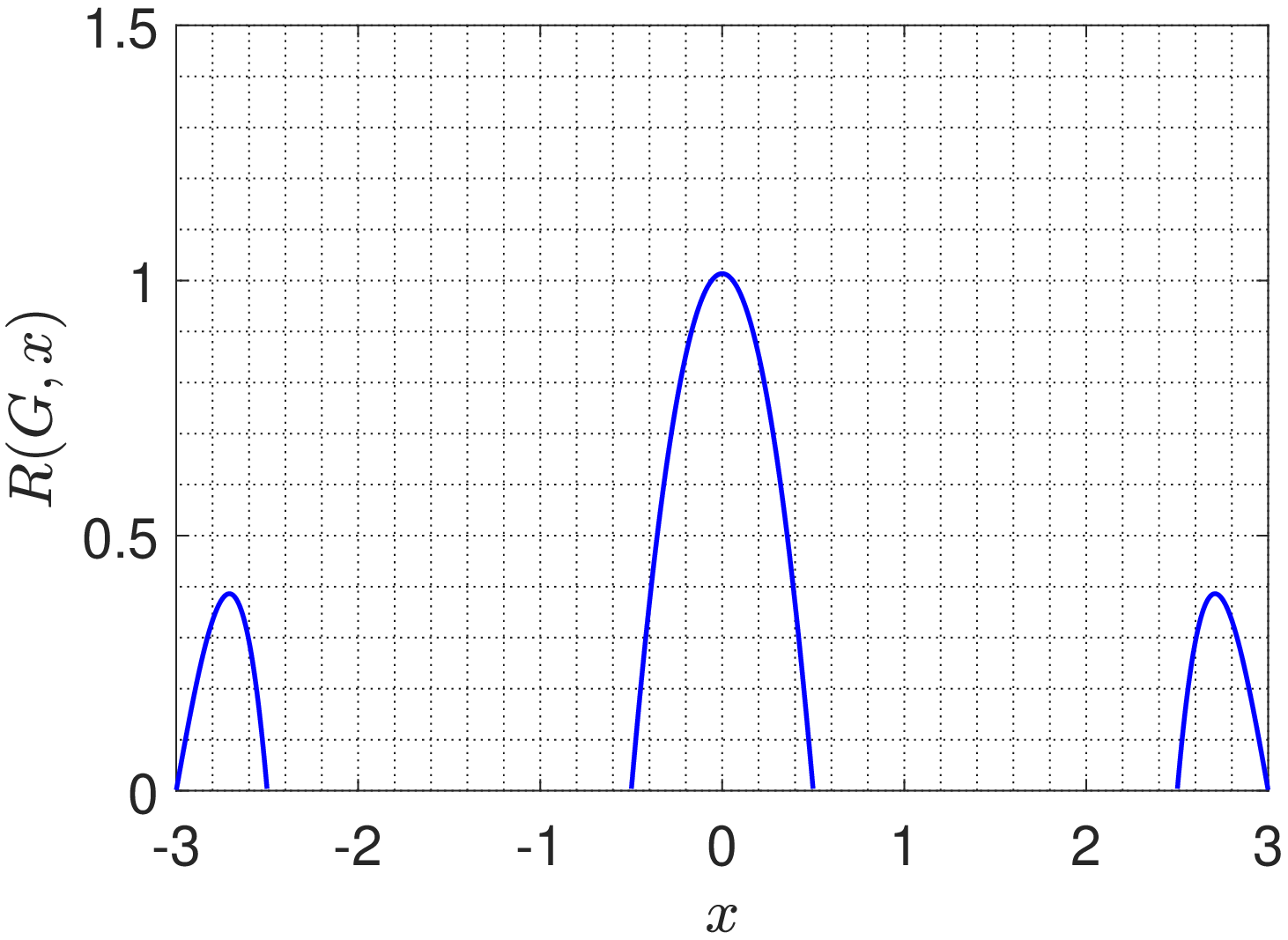}}
\hfill\scalebox{0.35}{\includegraphics[trim=0cm 0cm 0cm 0cm,clip]{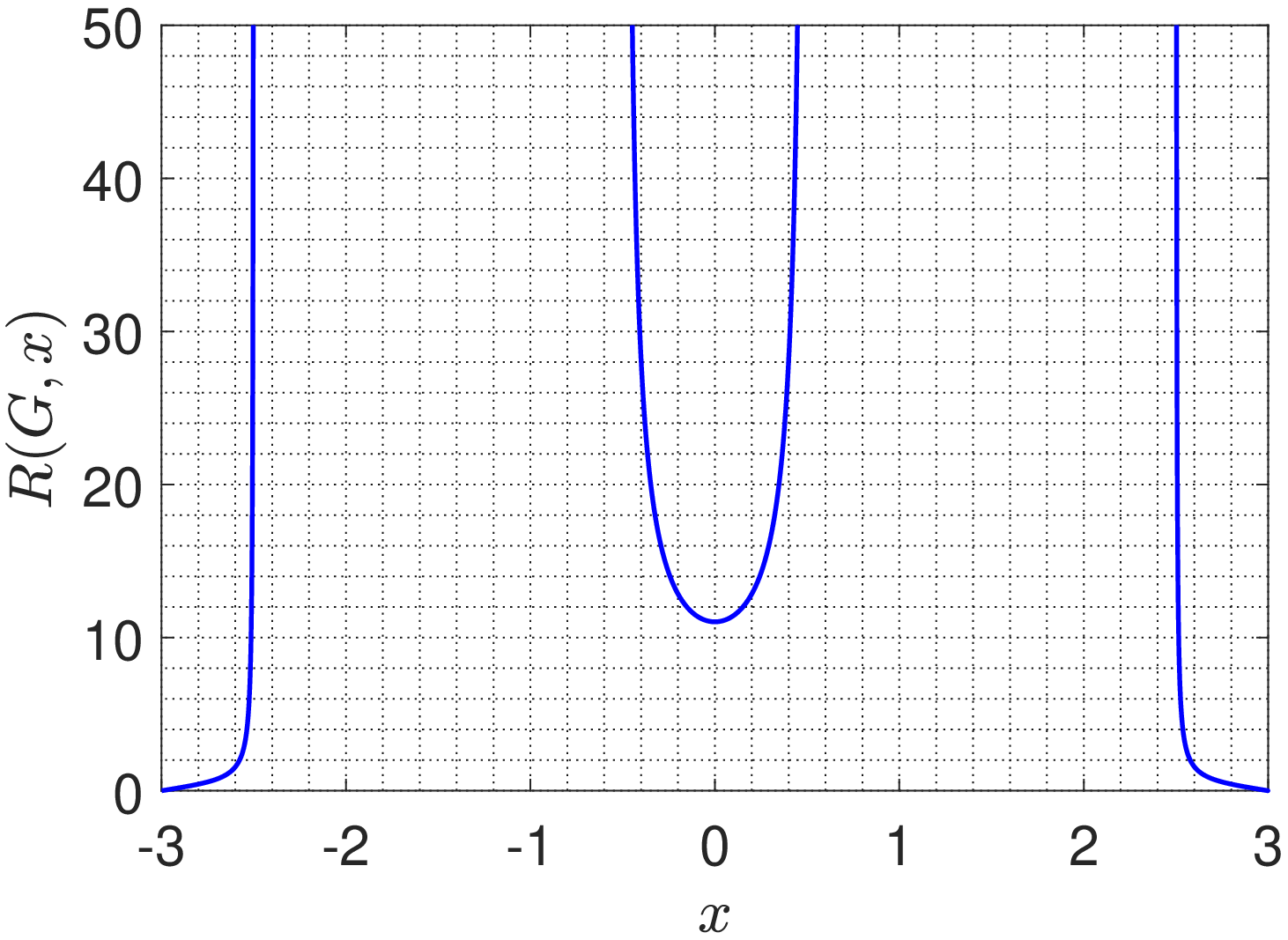}}
\hfill\scalebox{0.35}{\includegraphics[trim=0cm 0cm 0cm 0cm,clip]{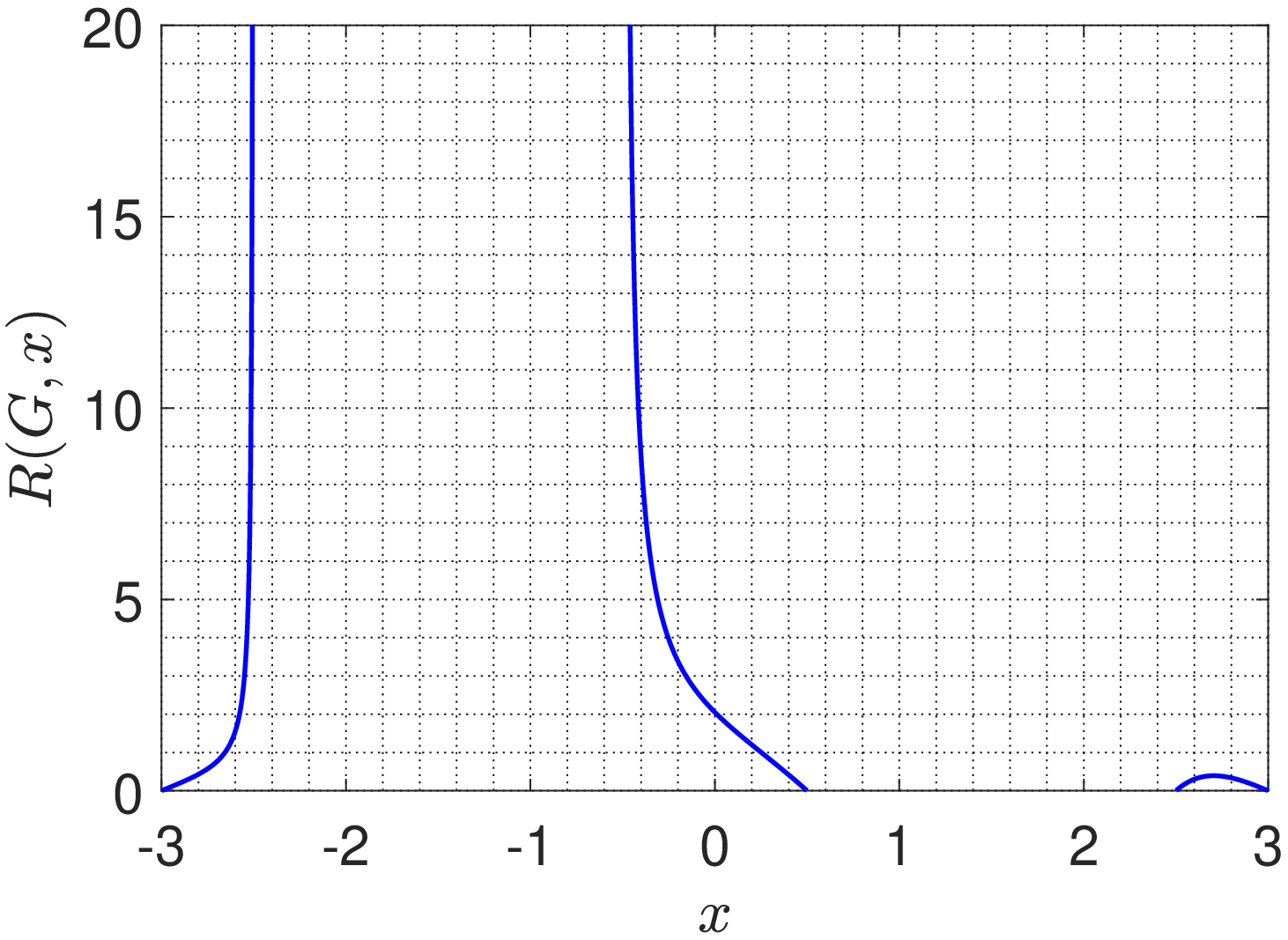}}
}
\caption{The values of $R(G,x)$ for $-3<x<3$ such that $x\in G$ for: $\theta_1=\theta_2=\pi/2$ (left), $\theta_1=\theta_2=0$ (center), and $\theta_1=\pi/2$, $\theta_2=0$ (right). The domain $G$ is the same as in Figure~\ref{fig:cir3L}.}
\label{fig:cir3x}
\end{figure}

Let $G$ be the triply connected domain interior to the circle $|z|=3$ and exterior to the circles $|z-1.5|=1$ and $|z+1.5|=1$.
We include three cases of the canonical domain: two circular slits, two radial slits, and one circular slit and one radial slit. The contour maps of the function $R(G,\alpha)$ are shown in Figure~\ref{fig:cir3L}.
We can see that for the unit disk with two circular slits, Mityuk's radius has five critical points ($n_m=2$ and $n_s=3$), where four of them are symmetric with respect to the saddle point $0$. Remark also that there is only one critical point (saddle) for the two remaining cases of the canonical domain.
The existence of critical points in the cases of two circular slits and a mix of circular and radial slits confirms the theoretical results of~\cite{Kin18,Mit}. On the other hand, the nature of critical points in the case of two radial slits opens a room for theoretical investigation.
Finally, we present in Figure~\ref{fig:cir3x} the graph of $R(G,x)$ with respect to $x\in G$ such that $-3<x<3$. It can be seen that the boundary  behavior of $R(G,x)$  agrees with~\eqref{eq:lim-k} and~\eqref{eq:lim-j} for the three cases.

\subsubsection{Six circles}

\begin{figure}[t] %
\centerline{
\hfill\scalebox{0.45}{      \includegraphics[trim=0cm 0cm 0cm 0cm,clip]{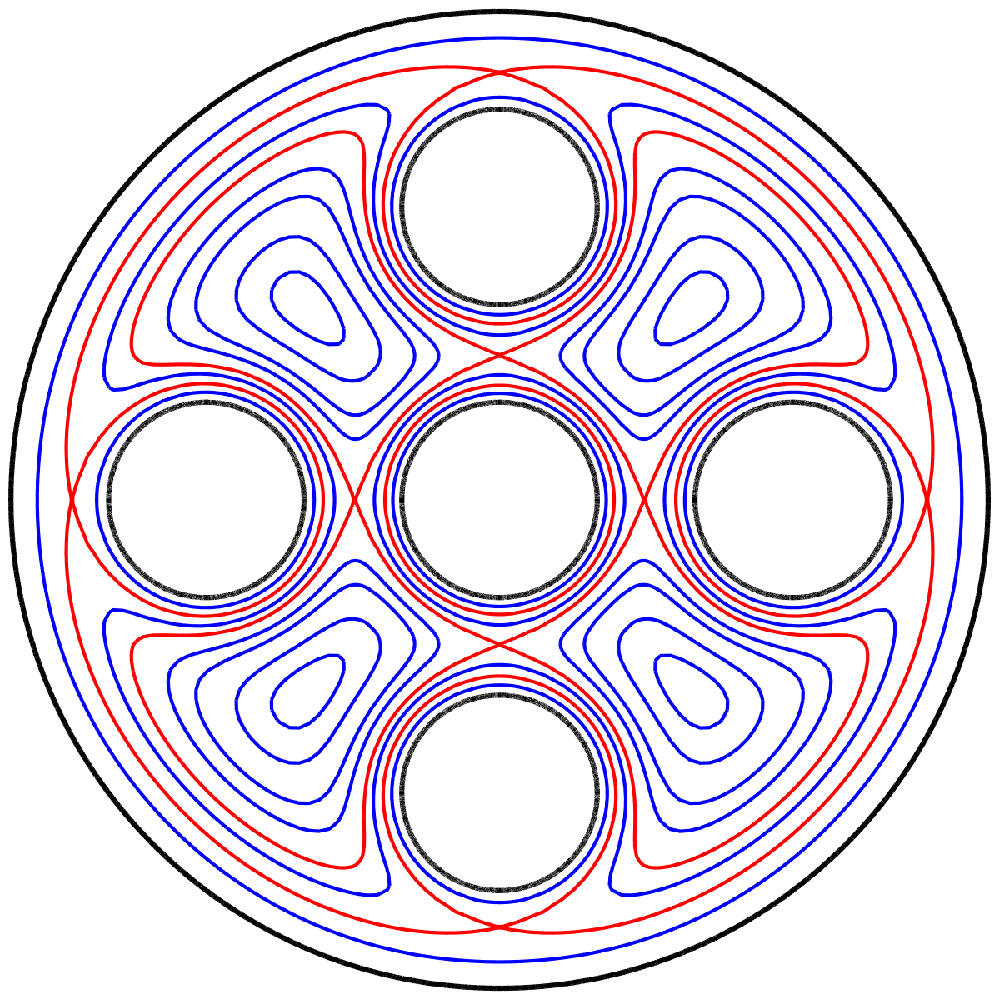}}
\hfill\scalebox{0.45}{\includegraphics[trim=0cm 0cm 0cm 0cm,clip]{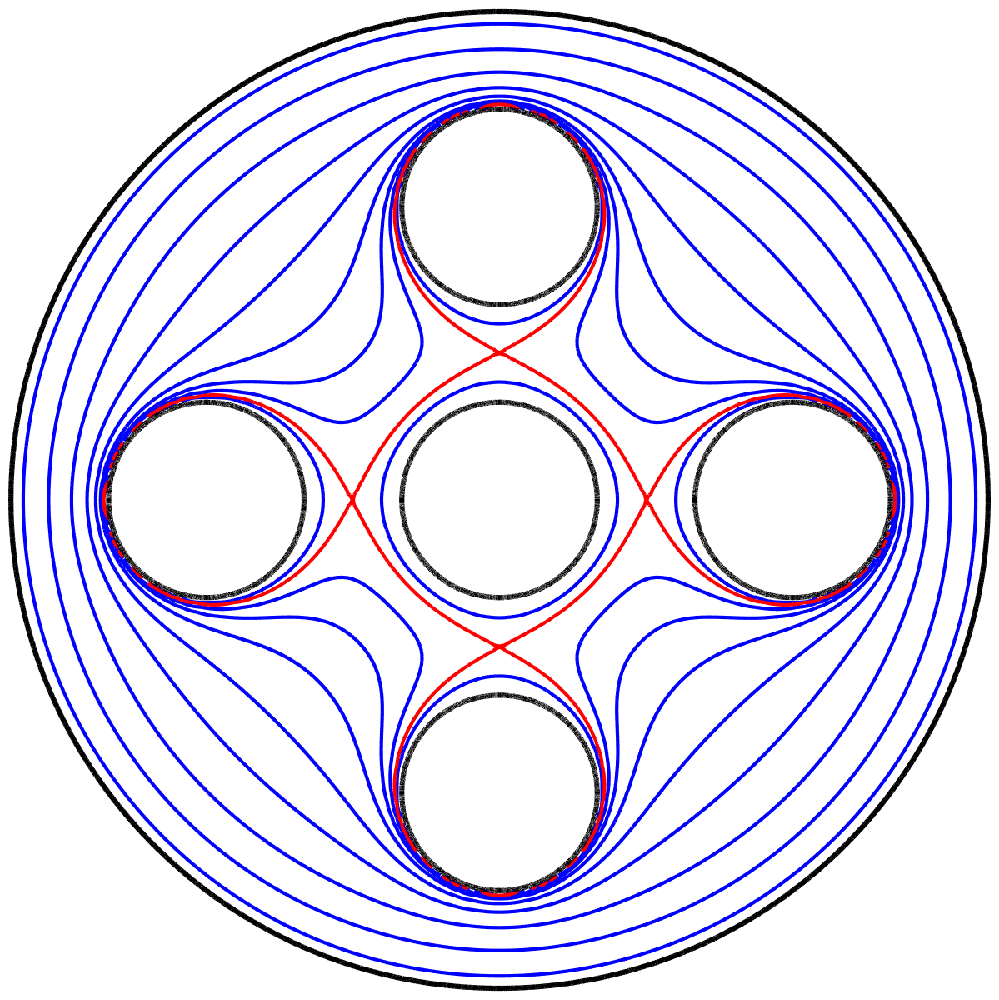}}
\hfill}
\centerline{
\scalebox{0.45}{      \includegraphics[trim=0cm 0cm 0cm 0cm,clip]{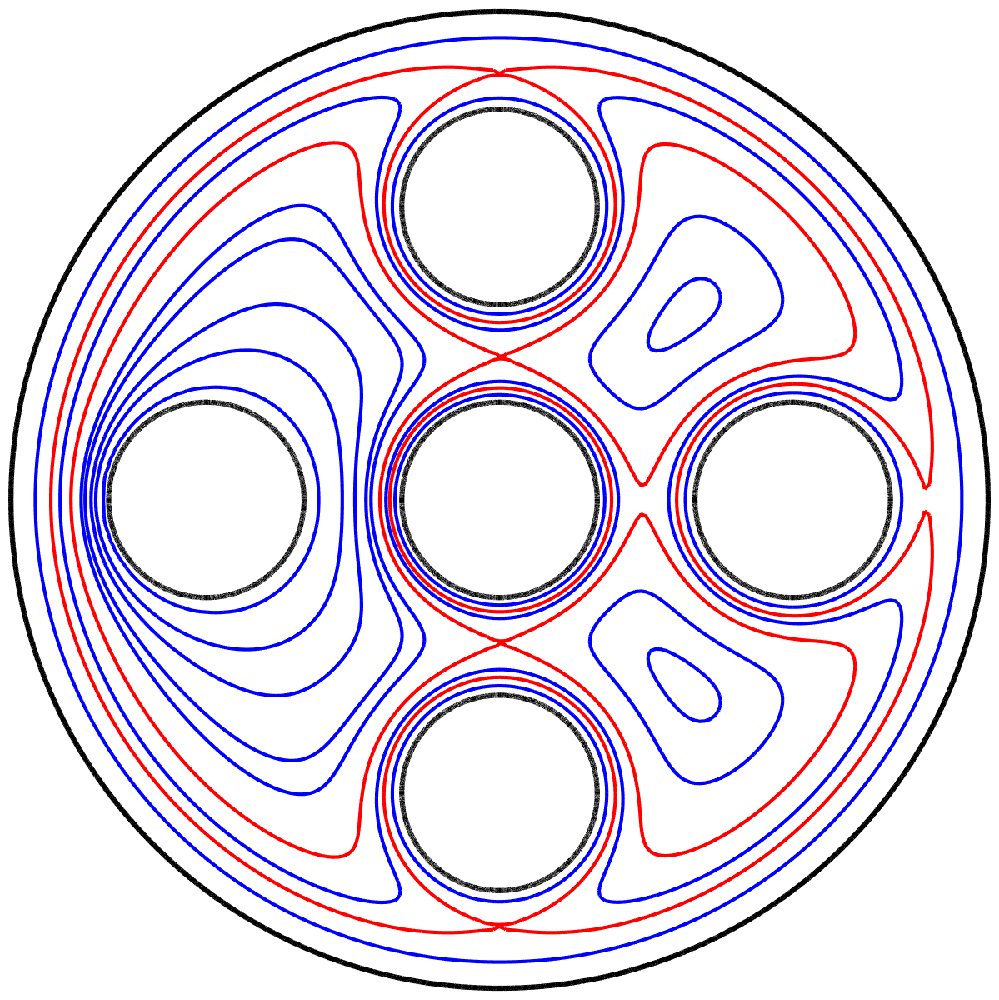}}
\hfill\scalebox{0.45}{\includegraphics[trim=0cm 0cm 0cm 0cm,clip]{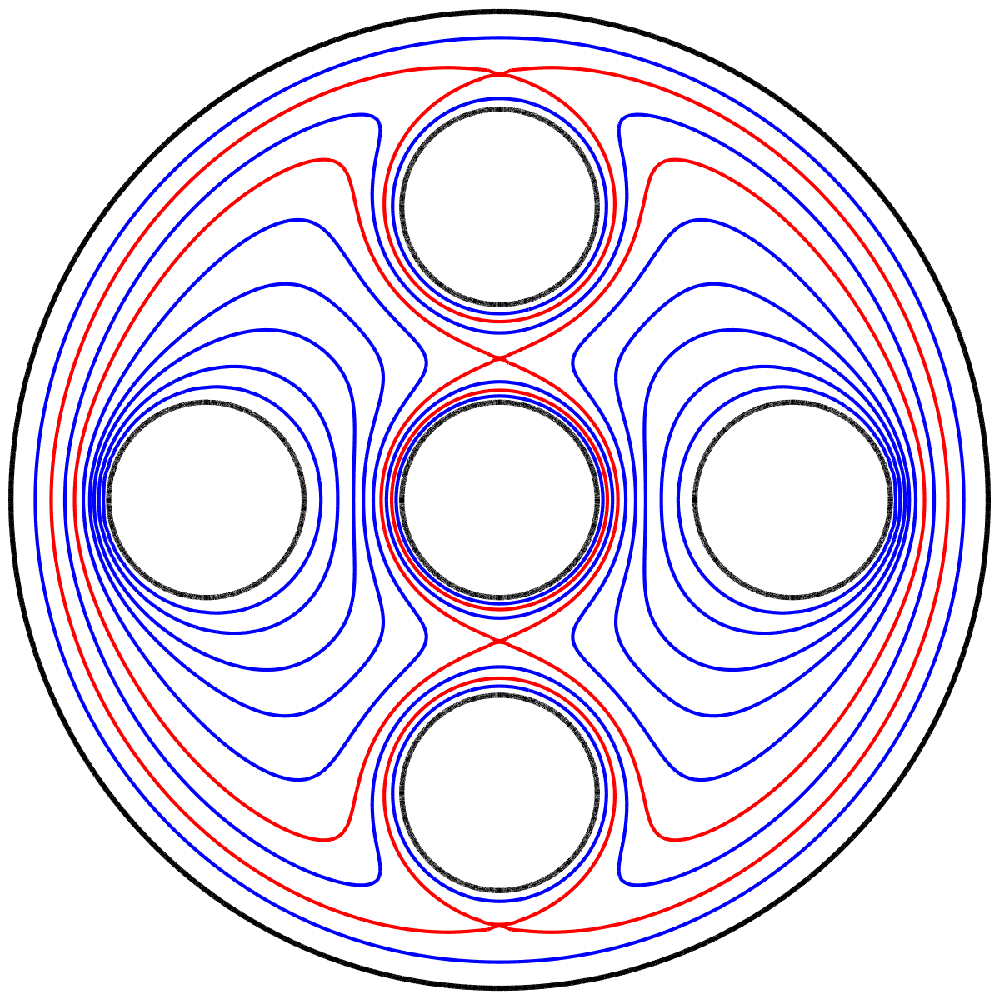}}
\hfill\scalebox{0.45}{\includegraphics[trim=0cm 0cm 0cm 0cm,clip]{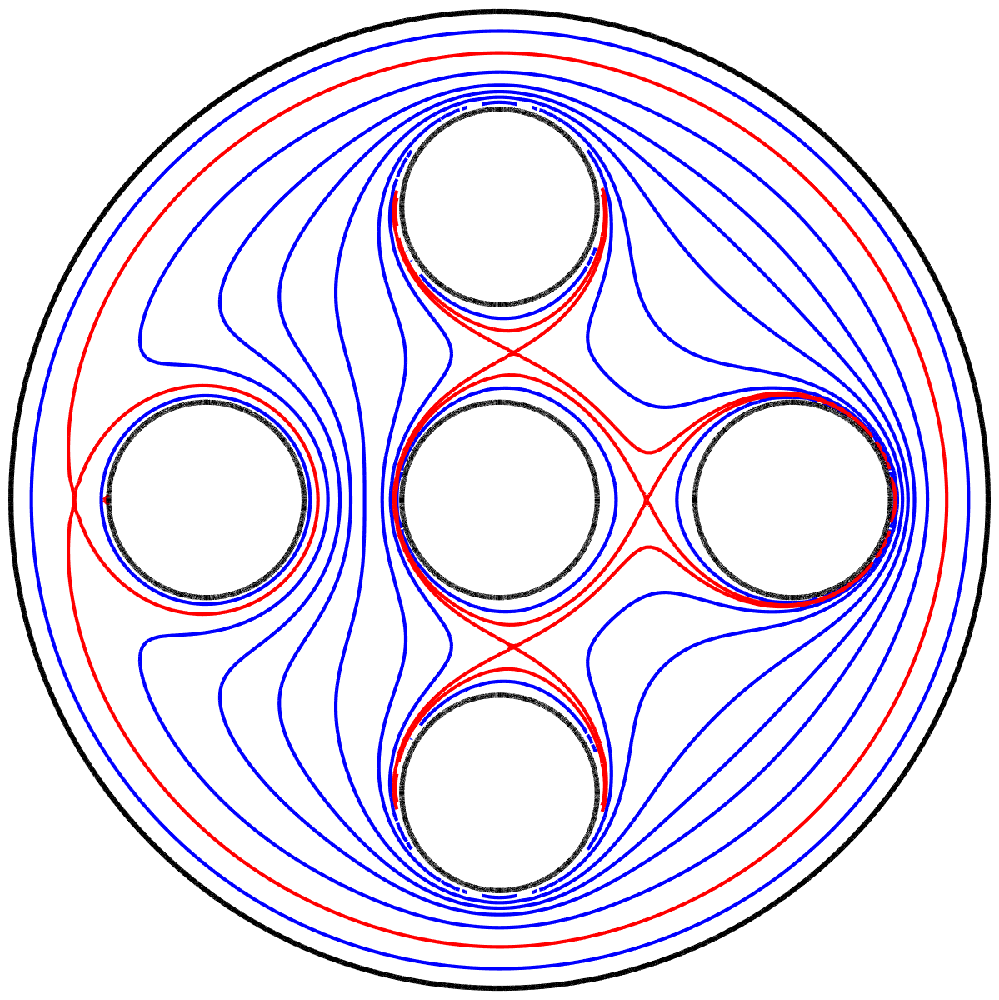}}
}
\caption{The contour maps of the function $R(G,\alpha)$ for: $\theta_1=\cdots=\theta_5=\pi/2$ (first row, left), $\theta_1=\cdots=\theta_5=0$ (first row, right), $\theta_1=\theta_2=\theta_3=\theta_5=\pi/2$, $\theta_4=0$ (second row, left), $\theta_1=\theta_3=\theta_5=\pi/2$, $\theta_2=\theta_4=0$ (second row, center), and $\theta_4=\pi/2$, $\theta_1=\theta_2=\theta_3=\theta_5=0$ (second row, right).  Critical streamlines are shown in red color.}
\label{fig:cir5L}
\end{figure}

Suppose $G$ is of connectivity $6$ being interior to the circle $|z|=1$ and exterior to the five circles with centers $0$, $0.6$, $0.6\i$, $-0.6$ and $-0.6\i$. The radii of all inner circles is $0.2$.
We include three cases of the canonical domain: five circular slits ($\theta_1=\cdots=\theta_5=\pi/2$), five radial slits ($\theta_1=\cdots=\theta_5=0$), and three circular slits and two radial slits ($\theta_1=\theta_3=\theta_5=\pi/2$ and $\theta_2=\theta_4=0$).
The contour maps of the function $R(G,\alpha)$ are shown in Figure~\ref{fig:cir5L}.
We can see that Mityuk's radius has twelve critical points ($n_m=4$ and $n_s=8$) for the case of five circular slits, and only four saddle points in the two other cases of the canonical domain. We remark also the symmetrical distribution of critical points location in all the three cases.

\subsubsection{Seven circles}

\begin{figure}[t] %
\centerline{
\scalebox{0.45}{      \includegraphics[trim=0cm 0cm 0cm 0cm,clip]{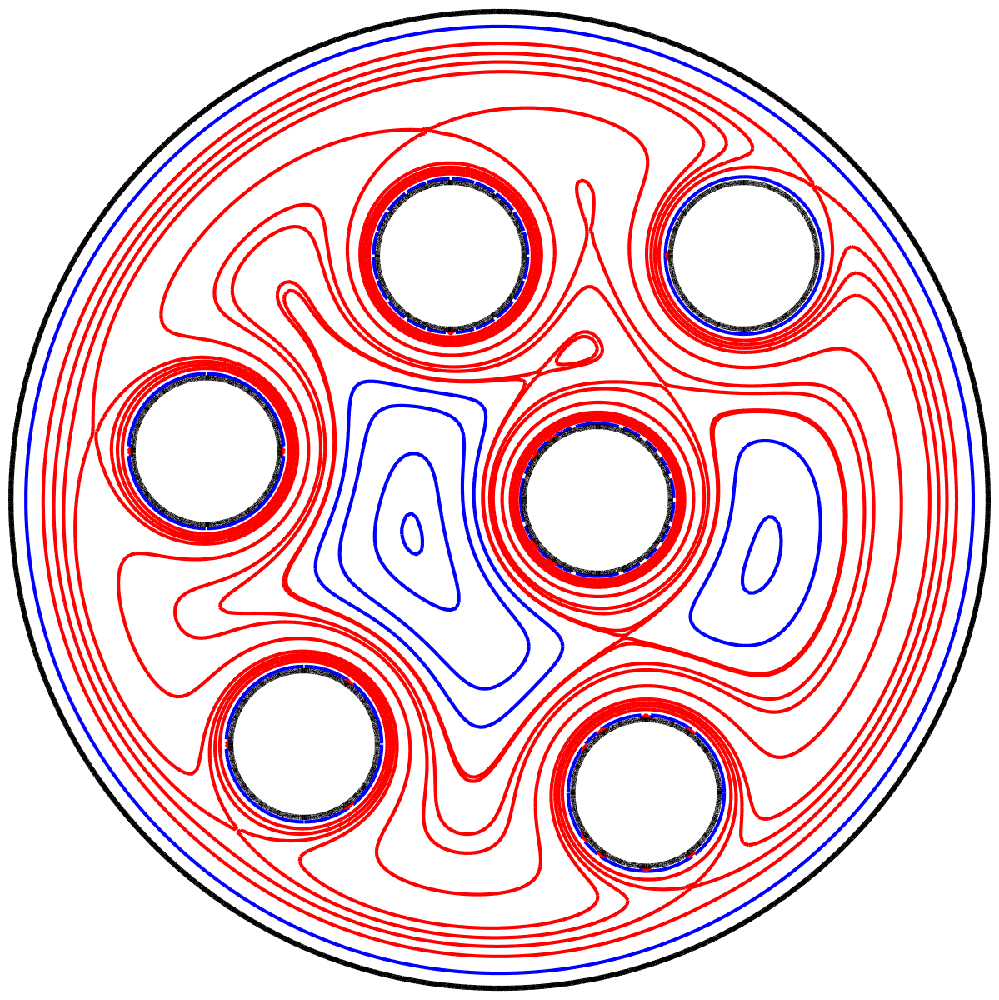}}
\hfill\scalebox{0.45}{\includegraphics[trim=0cm 0cm 0cm 0cm,clip]{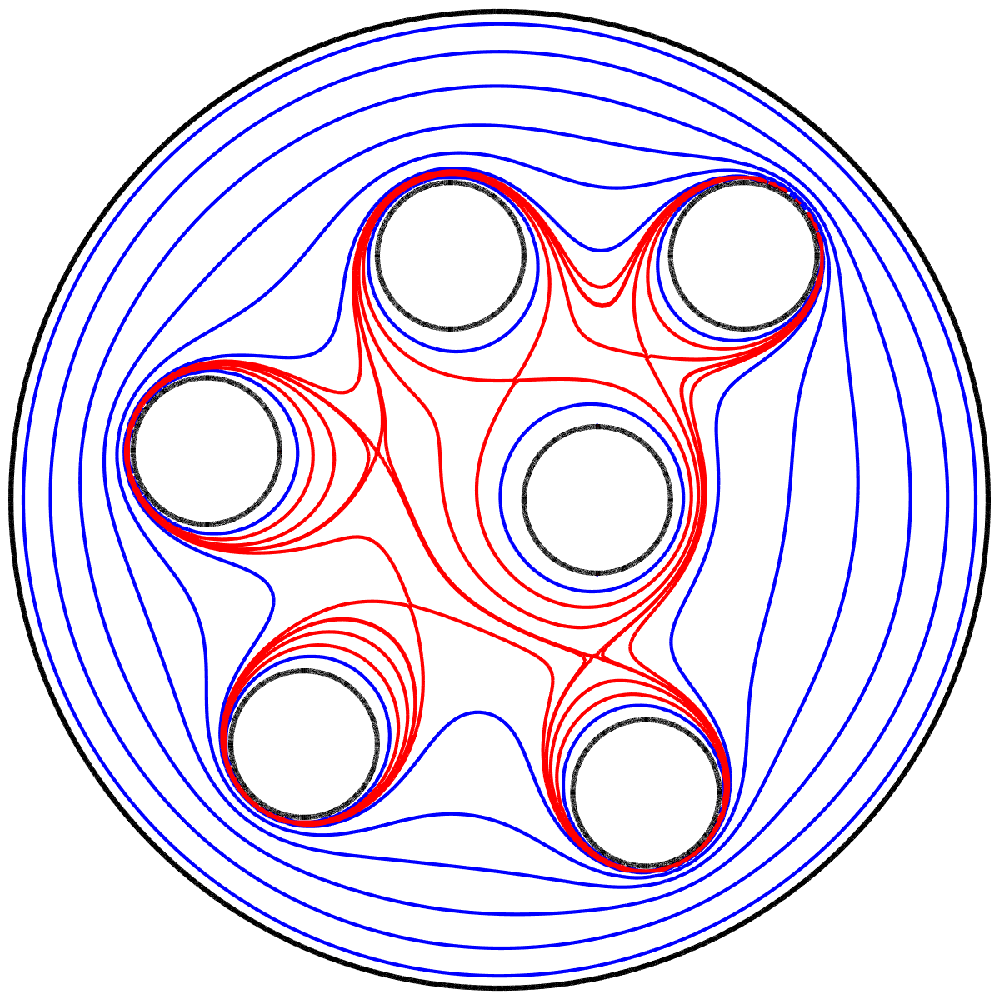}}
\hfill\scalebox{0.45}{\includegraphics[trim=0cm 0cm 0cm 0cm,clip]{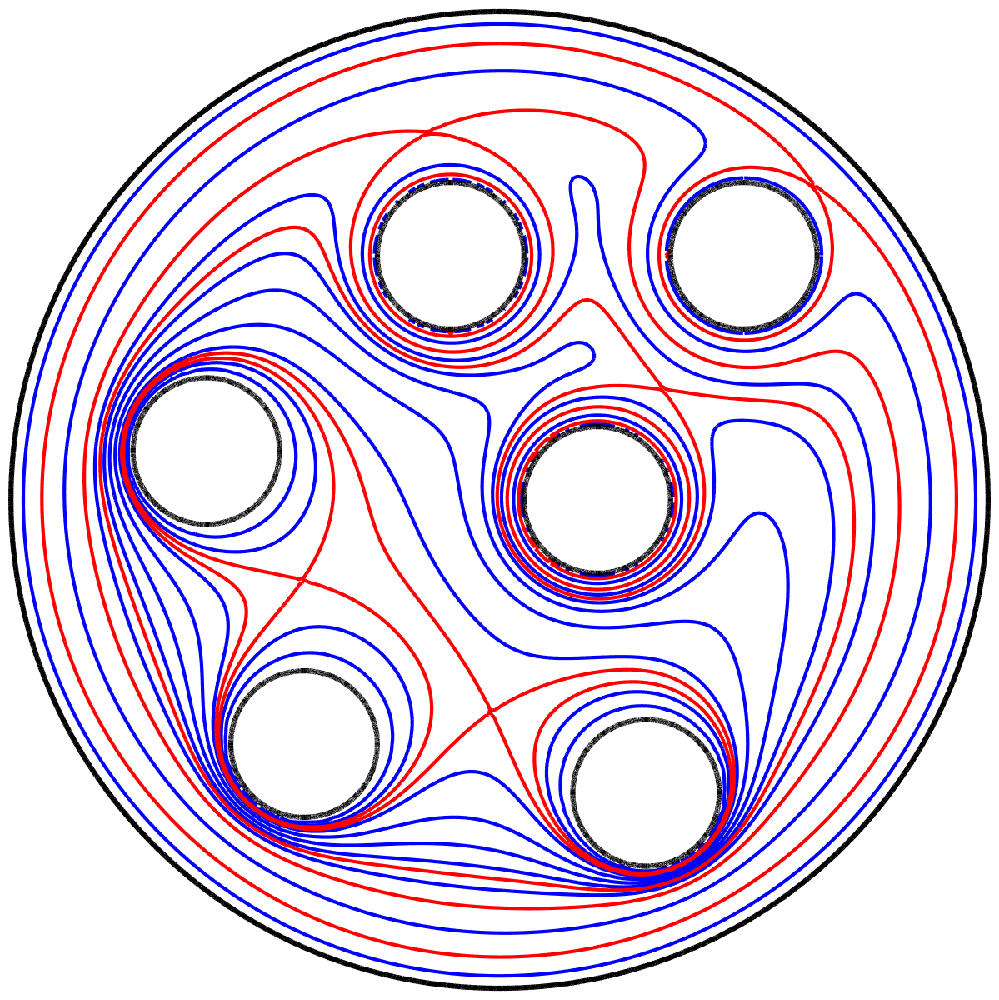}}}
\caption{The contour maps of the function $R(G,\alpha)$ for: $\theta_1=\cdots=\theta_6=\pi/2$ (left), $\theta_1=\cdots=\theta_6=0$ (center), $\theta_1=\theta_2=\theta_3=\pi/2$, $\theta_4=\theta_5=\theta_6=0$ (right).  Critical streamlines are shown in red color.}
\label{fig:cir7L}
\end{figure}

Assume $G$ is of connectivity $7$, interior to the circle $|z|=1$ and exterior to the five circles with centers $0.2$, $0.5+0.5\i$, $-0.1+0.5\i$, $-0.6+0.1\i$, $-0.4-0.5\i$, and $0.3-0.6\i$. The radii of all inner circles is $0.15$.
We include three cases of the canonical domain: six circular slits ($\theta_1=\cdots=\theta_6=\pi/2$), five radial slits ($\theta_1=\cdots=\theta_6=0$), and three circular slits and three radial slits ($\theta_1=\theta_2=\theta_3=\pi/2$ and $\theta_4=\theta_5=\theta_6=0$).
The contour maps of the function $R(G,\alpha)$ are shown in Figure~\ref{fig:cir7L}. It is apparent from this figure that the critical points of Mityuk's radius satisfy the equation $n_m-n_s=-5$ for the three cases of the canonical domain.

\section{Conclusion}
We have introduced a numerical method to compute Mityuk's function and radius of multiply connected domains with respect to the canonical domain consisting of the unit disk with $\ell$ circular/radial slits. Our main interest has been to validate the theoretical results on the existence and nature of critical points of Mityuk's radius as well as the boundary behavior of this function as proven in~\cite{Kin1,Kin18,Kin84}.
For doubly connected domains, we have first considered the special case of an annulus. 
{\color{blue}
In contrast to the case of the unit disk with a radial slit where Mityuk's radius has no critical points~\cite{Kin18}, 
there is an infinite number of critical points located on the circle with center zero and radius equal to the square root of the inner circle radius in the other case of the unit disk with a circular slit~\cite{Kin1}.
The other examples of doubly connected domains exhibit the same property of existence of critical points for a circular slit and nonexistence for a radial slit. Except for the case of an annulus, where there is an infinite number of critical points, equation~\eqref{eq:criticalPts} is also true for doubly connected domains. In particular, the number of maxima is equal to the number of saddle points~\cite{Kin84}.}

For multiply connected domains of connectivity $\ell \geq 2$, we have presented several examples for which Mityuk's function/radius admits critical points for different cases of the canonical domain.
{\color{blue}
The numerical results illustrate that Mityuk's radius has at least one local maximum if the canonical domain contains only circular slits.}
Concerning the boundary behavior of Mityuk's radius, we have remarked that the limit at the external boundary component is always equal to zero even when the boundary is piecewise smooth. On the other hand, the limit values in~\eqref{eq:lim-j} do not hold true when one of the internal curves is a slit mapped to a radial slit.
{\color{blue}
Two final observations on the numerical results presented in this paper that might open a room for theoretical investigation are given as follows:}
\begin{itemize}
	\item The symmetry and geometry of the domain is directly reflected into the location and number of critical points. Mityuk's radius has more critical points for symmetric domains compared to non-symmetric domains of the same connectivity.
	\item For Mityuk's radius $R(G, \alpha)$ of the domain $G$ with respect to the point $\alpha\in G$ and the canonical domain consisting of the unit disk with circular slits, it follows from the limits~\eqref{eq:lim-k} and~\eqref{eq:lim-j} that $R(G,\alpha)$ approaches $0$ as $\alpha$ approaches the boundary $\Gamma=\partial G$. Our numerical tests suggest that
\noindent
\begin{equation} \label{MitLoBound} R(G, \alpha) \ge d(\alpha,\partial
G)
\end{equation}
\noindent
holds for all $\alpha \in G$ and for the all multiply connected domains $G$ considered in this paper. This observation suggests the
following question: ``Is it true that \eqref{MitLoBound} holds for all finitely connected domains $G$ bordered by  Jordan curves?''

\end{itemize}

{\color{blue}
\section*{Acknowledgements} 
The authors are grateful to an anonymous referee for his valuable comments and suggestions, which improved the presentation of this paper. 

}
	



\begin{thebibliography}{99}	


\bibitem{Kin1}
{\sc L.A. Aksent{'}ev, M.I. Kinder, and S.B. Sagitova}, \emph{Solvability of the exterior inverse boundary value problem in the case of a multiply connected domain}, Trudy Sem. Kraev. Zadacham. 20 (1983) 22--34   [in Russian].
		
\bibitem{Atk97}
{\sc K.E. Atkinson}, {\em The Numerical Solution of Integral Equations of the Second Kind}, Cambridge University Press, 1997.

\bibitem{Tre-Trap1}
{\sc A.P. Austin, P. Kravanja, and L.N. Trefethen}, \emph{Numerical algorithms based on analytic function values at roots of unity}, SIAM J. Numer. Anal. 52 (2014) 1795--1821.

\bibitem{CF07}
{\sc D. Crowdy and A. Fokas}, {\em Conformal mappings to a doubly connected polycircular arc domain}, Proc. R. Soc. A 463 (2007) 1885--1907.

\bibitem{CM05}
{\sc D. Crowdy and J. Marshall}, {\em Analytic formulae for the Kirchhoff--Routh path function in multiply connected domains}, Proc. Roy. Soc. A. 461 (2005) 2477--2501.

\bibitem{CM06}
{\sc D. Crowdy and J. Marshall}, {\em Conformal mappings between canonical multiply connected domains}, Comput. Methods Funct. Theory 6(1) (2006) 59--76.

\bibitem{Du} { \sc V. N. Dubinin}, {\em Condenser Capacities and Symmetrization in Geometric Function Theory,} Birkh\"auser, 2014.
		
\bibitem{Eli17}
{\sc A.M. Elizarov, A.V. Kazantsev, and M.I. Kinder}, {\em Strict superharmonicity of Mityuk's function for countably connected domains of simple structure}, Lobachevskii Journal of Mathematics 38(3) (2017) 408--413.

\bibitem{Kin18}
{\sc A.M. Elizarov, A.V. Kazantsev, and M.I. Kinder}, {\em Generalized reduced module of a domain over the unit disc with circular and radial slits}, Lobachevskii Journal of Mathematics 39(5) (2018) 664--672.

\bibitem{gak77}
{\sc F.D. Gakhov}, {\em Boundary Value Problems}, Pergamon Press, Oxford, 1966.

\bibitem{Gol69}
{\sc G. Goluzin}, {\em Geometric Theory of Functions of a Complex Variable},
Amer. Math. Soc., Rhode Island, 1969.


\bibitem{Gre-Gim12}
{\sc L. Greengard and Z. Gimbutas},
{\em FMMLIB2D: A {MATLAB} toolbox for fast multipole method in two dimensions},
Version 1.2, \url{http://www.cims.nyu.edu/cmcl/fmm2dlib/fmm2dlib.html}. Accessed 1 Jan 2018.

\bibitem{Kaz17}
{\sc A.V. Kazantsev}, {\em Sectio Aurea conditions for Mityuk's radius of two-connected domains}, Uch. Zap. Kazan. Univ., Ser. Fiz.-Mat. Nauki 159 (2017) 33--46.		

\bibitem{Kin84}
{\sc M.I. Kinder}, {\em The number of solutions of F. D. Gakhov's equation in the case of a multiply connected domain},
Sov. Math. (Iz. VUZ) 28(8) (1984) 91--95.

\bibitem{Kre}{\sc R. Kress}, 
{\em A Nystr\"om method for boundary integral equations in domains with corners}, 
Numer. Math. 58(2) (1990) 145--161.

\bibitem{LSN}
{\sc J. Liesen, O. S\'ete and M.M.S. Nasser}, {\em A fast and accurate computation of the logarithmic capacity of compact sets},  Comput. Methods Funct. Theory  17 (2017) 689--713.
 		
\bibitem{Mit}
{\sc I.P. Mityuk}, {\em A generalized reduced module and some of its applications},  Izv. Vyssh. Uchebn. Zaved., Mat. 2 (1964) 110--119  [in Russian].
		
\bibitem{Mur-Bul}
{\sc A.H.M. Murid and M.M.S. Nasser}, {\em Eigenproblem of the generalized
{N}eumann kernel}, Bull. Malays. Math. Sci. Soc. 26 (2003) 13--33.

\bibitem{Nas-jmaa11}
{\sc M.M.S. Nasser}, {\em Numerical conformal mapping of multiply connected regions onto the second, third and fourth categories of Koebe's canonical slit domains}, J. Math. Anal. Appl. 382 (2011) 47--56.

\bibitem{Nas-ETNA}
{\sc M.M.S. Nasser}, {\em Fast solution of boundary integral equations with the generalized Neumann kernel}, Electron. Trans. Numer. Anal. 44 (2015) 189--229.

\bibitem{NSML}
{\sc M.M.S. Nasser, T. Sakajo, A.H.M. Murid and L.K. Wei}, {\em A fast computational method for potential flows in multiply connected coastal domains}, Japan J. Indust. Appl. Math. 32 (2015) 205--236.
		
\bibitem{Tre-Trap}
{\sc L.N. Trefethen and J. A. C. Weideman}, \emph{The exponentially convergent trapezoidal rule}, SIAM Review 56 (2014) 385--458.
		
\bibitem{Vas02}
{\sc A. Vasil'ev}, {\em Moduli of Families of Curves for Conformal and Quasiconformal Mappings}, Springer-Verlag, Berlin, 2002.
		
\bibitem{Weg-Mur-Nas}
{\sc R. Wegmann, A.H.M. Murid, and M.M.S. Nasser}, {\em The {R}iemann-{H}ilbert
problem and the generalized {N}eumann kernel}, J. Comput. Appl. Math. 182
(2005) 388--415.
		
\bibitem{Weg-Nas}
{\sc R. Wegmann and M.M.S. Nasser}, {\em The {R}iemann-{H}ilbert problem and
the generalized {N}eumann kernel on multiply connected regions}, J. Comput.
Appl. Math. 214 (2008) 36--57.

\bibitem{Wen92}
{\sc G.C. Wen}, {\em Conformal Mapping and Boundary Value Problems}, Amer. Math. Soc., Providence, 1992.
		
\end{thebibliography}
\end{document}